\documentclass{amsproc}

\usepackage[arrow, curve, frame, matrix, graph, tips]{xy}
\usepackage{amssymb}
\usepackage{amsthm}
\usepackage{mtens}

\newtheorem{theorem}{Theorem}[section]
\newtheorem{lemma}[theorem]{Lemma}
\newtheorem{proposition}[theorem]{Proposition}
\newtheorem{corollary}[theorem]{Corollary}

\theoremstyle{definition}
\newtheorem{definition}[theorem]{Definition}
\newtheorem{example}[theorem]{Example}
\newtheorem{examples}[theorem]{Examples}
\newtheorem{non-example}[theorem]{Non-Example}

\theoremstyle{remark}
\newtheorem{remark}[theorem]{Remark}

\newcommand{\TriTwoCell}[7]{\xymatrix{
{#1} \ar[rr]^-{#5} \ar[dr]_{#4} \save \POS?="dom" \restore
&& {#2} \ar[dl]^{#6} \save \POS?="cod" \restore \\ & {#3}
\POS "dom"; "cod" **@{} ?(.35) \ar@{=>}^{#7} ?(.65)
}}

\newcommand{\TwoDiagRel}[3]{\begin{xy}
(0,0);<2em,0em>:<0em,2em>::
(-1,0)*-!R{\xybox{#1}};
(1,0)*-!L{\xybox{#3}}
**@{} ?*{#2}
\end{xy}}

\newcommand{\LaxSq}[9]{\xymatrix{{#1} \ar[r]^-{#6}
\ar[d]_{#5} \save \POS?="dom" \restore
& {#2} \ar[d]^{#7} \save \POS?="cod" \restore \\
{#3} \ar[r]_-{#8} & {#4}
\POS "dom"; "cod" **@{} ?(.35) \ar@{=>}^{#9} ?(.65)}}

\newcommand{\MNDmorphism}[6]{\xymatrix{{#1} \ar[r]^-{#5} \ar[d]_{#2}
& {#3} \save \POS="dom" \restore \ar[d]^{#4} \\
{#1} \save \POS="cod" \restore \ar[r]_-{#5} & {#3}
\POS "dom"; "cod" **@{} ?(.35) \ar@{=>}^{#6} ?(.65)}}

\newcommand{\PbSq}[8]{\xymatrix{{#1} \ar[d]_{#5}
\save \POS?(.3)="lpb" \restore
\ar[r]^-{#6} \save \POS?(.3)="tpb" \restore
& {#2} \ar[d]^{#7} \save \POS?(.3)="rpb" \restore \\
{#4} \ar[r]_-{#8} \save \POS?(.3)="bpb" \restore & {#3}
\POS "rpb"; "lpb" **@{}; ?!{"bpb";"tpb"}="cpb" **@{}; ? **@{-};
"tpb"; "cpb" **@{}; ? **@{-}}}

\newcommand{\rel}{\SelectTips{cm}{}\object@{/}}

\begin{document}

\title{Algebras of higher operads as enriched categories II}

\author{Michael Batanin}
\address{Department of Mathematics,
Macquarie University}
\email{mbatanin@ics.mq.edu.au}
\thanks{}
\author{Denis-Charles Cisinski}
\address{Departement des Mathematiques,
Universit\'e Paris 13 Villanteuse}
\email{cisinski@math.paris13.fr}
\thanks{}
\author{Mark Weber}
\address{Max Planck Institute for Mathematics, Bonn}
\email{mark.weber.math@gmail.com}
\thanks{}
\maketitle
\begin{abstract}
One of the open problems in higher category theory is the systematic construction of the higher dimensional analogues of the Gray tensor product. In this paper we continue the work of \cite{EnHopI} to adapt the machinery of globular operads \cite{Bat98} to this task. The resulting theory includes the Gray tensor product of 2-categories and the Crans tensor product \cite{Crans99} of Gray categories. Moreover much of the previous work on the globular approach to higher category theory is simplified by our new foundations, and we illustrate this by giving an expedited account of many aspects of Cheng's analysis \cite{ChengCompOp} of Trimble's definition of weak $n$-category. By way of application we obtain an ``Ekmann-Hilton'' result for braided monoidal 2-categories, and give the construction of a tensor product of $A$-infinity algebras.
\end{abstract}
\tableofcontents

\section{Introduction}\label{sec:Intro}

In \cite{Bat98} the problem of how to give an explicit combinatorial definition of weak higher categories was solved, and the development of a conceptual framework for their further analysis was begun. In the aftermath of this, the expository work of other authors, most notably Street \cite{Str98} and Leinster \cite{Lei}, contributed greatly to our understanding of these ideas. The central idea of \cite{Bat98} is that the description of any $n$-dimensional categorical structure $X$, may begin by starting with just the underlying $n$-globular set, that is, the sets and functions
\[ \xymatrix{{X_0} & {X_1} \ar@<1ex>[l]^{t}  \ar@<-1ex>[l]_{s} & {X_2} \ar@<1ex>[l]^{t}  \ar@<-1ex>[l]_{s}
& {X_3} \ar@<1ex>[l]^{t}  \ar@<-1ex>[l]_{s} & {...} \ar@<1ex>[l]^{t}  \ar@<-1ex>[l]_{s} & {X_n} \ar@<1ex>[l]^{t}  \ar@<-1ex>[l]_{s}} \]
satisfying the equations $ss=st$ and $ts=tt$, which embody the the objects (elements of $X_0$), arrows (elements of $X_1$) and higher cells of the structure in question. At this stage no compositions have been defined, and when they are, one has a globular set with extra structure. In this way the problem of defining an n-categorical structure of a given type is that of defining the monad on the category $\PSh {\G}_{{\leq}n}$ of $n$-globular sets whose algebras are these structures.

As explained in the introduction to \cite{EnHopI}, this approach works because the monads concerned have excellent formal properties, which facilitate their explicit description and further analysis. The $n$-operads of \cite{Bat98} can be defined from the point of view of monads: one has the monad $\ca T_{{\leq}n}$ on $\PSh {\G}_{{\leq}n}$ whose algebras are strict $n$-categories, and an $n$-operad consists of another monad $A$ on $\PSh {\G}_{{\leq}n}$ equipped with a cartesian monad morphism $A \rightarrow \ca T_{{\leq}n}$. The algebras of this $n$-operad are just the algebras of $A$.

Strict $n$-categories are easily defined by iterated enrichment: a strict $(n{+}1)$-category is a category enriched in the category of strict $n$-categories via its cartesian product, but are too strict for the intended applications in homotopy theory and geometry. For $n=3$ the strictest structure one can replace an arbitrary weak $3$-category with -- and not lose vital information -- is a Gray category, which is a category enriched in $\Enrich 2$ using the Gray tensor product of 2-categories instead of its cartesian product \cite{GPS95}. This experience leads naturally to the idea of trying to define what the higher dimensional analogues of the Gray tensor product are, so as to set up a similar inductive definition as for strict $n$-categories, but to capture the appropriate semi-strict $n$-categories, which in the appropriate sense, would form the strictest structure one can replace an arbitrary weak $n$-category with and not lose vital information.

Crans in \cite{Crans99} attempted to realise this idea in dimension 4, and one of our main motivations is to obtain a theory that will deliver the sort of tensor product that Crans was trying to define explicitly, but in a conceptual way that one could hope to generalise to still higher dimensions. Our examples(\ref{ex:Gray}) and (\ref{ex:Crans}) embody the progress that we have achieved in this direction in this paper. In \cite{WebFunny} the theory of the present paper is used to show that the \emph{funny tensor product} of categories -- which is what one obtains by considering the Gray tensor product of $2$-categories but ignoring what happens to 2-cells -- generalises to give an analogous symmetric monoidal closed structure on the category of algebras of any higher operad. From these developments it seems that a conceptual understanding of the higher dimensional analogues of the Gray tensor product is within reach.

Fundamentally, we have two kinds of combinatorial objects important for the description and study of higher categorical structures -- $n$-operads and tensor products. In \cite{EnHopI} a description of the relationship between tensor products and $n$-operads was begun, and $(n{+}1)$-operads whose algebras involve no structure at the level objects{\footnotemark{\footnotetext{In \cite{EnHopI} these were called \emph{normalised} $(n{+}1)$-operads. In the present work we shall, for reasons that will become apparent below, refer to these operads as \emph{being over $\Set$}.}}} were canonically related with certain lax tensor products on $\PSh {\G}_{{\leq}n}$. Under this correspondence the algebras of the $(n{+}1)$-operad coincide with categories enriched in the associated lax tensor product.

Sections(\ref{sec:EG-LMC})-(\ref{sec:2-functoriality}) of the present paper continue this development by studying, for a given category $V$, the passage
\[ \textnormal{Lax tensor products on $V$} \mapsto \textnormal{Monads on $\ca GV$} \]
where $\ca GV$ is the category of graphs enriched in $V$, in a systematic way. This analysis culminates in section(\ref{sec:2-functoriality}) where the above assignment is seen as the object part of a 2-functor
\[ \Gamma : \DISTMULT \rightarrow \MND(\CAT/\Set) \]
where $\DISTMULT$ is a sub 2-category of the 2-category of lax monoidal categories, and $\MND$ is as defined by the formal theory of monads \cite{Str72}. From this perspective, one is able to describe in a more efficient and general way, many of the previous developments of higher category theory in the globular style. For instance, in section(\ref{ssec:induction}) we give a short and direct explicit construction of the monads $\ca T_{{\leq}n}$ for strict $n$-categories from which all their key properties are easily witnessed. In sections(\ref{ssec:general-op-mult}) and (\ref{ssec:induction}) we give shorter and more general proofs of some of the main results of \cite{EnHopI}. In section(\ref{ssec:Monmonad-Distlaw}) using a dual version of our 2-functor $\Gamma$ and the formal theory of monads \cite{Str72}, we obtain a satisfying general explanation for how it is that monad distributive laws arise in higher category theory -- see \cite{ChengDist} \cite{ChengCompOp}. In sections(\ref{ssec:TCI}) and (\ref{ssec:TCII}) we apply our theory to simplifying many aspects of \cite{ChengCompOp}.

The correspondence between $(n{+}1)$-operads and certain lax monoidal structures on $\PSh {\G}_{{\leq}n}$ given in \cite{EnHopI}, associates to the 3-operad $G$ for Gray categories, a lax tensor product on the category of 2-globular sets. However the Gray tensor product itself is a tensor product of 2-categories. Any lax monoidal structure on a category $V$ comes with a ``unary'' tensor product, which rather than being trivial as is the usual experience with non-lax tensor products, is in fact a monad on $V$. For the lax tensor product induced by $G$, this is the monad for 2-categories. In section(\ref{sec:lift-mult}) we solve the general problem of lifting a lax monoidal structure, to a tensor product on the category of algebras of the monad defined by its unary part. This result, theorem(\ref{thm:lift-mult}), is the main result of the paper, and provides also the sense in which these lifted tensor products are unique. In practical terms this means that in order to exhibit a given tensor product on some category of higher dimensional structures as arising from our machinery, it suffices to exhibit an operad whose algebras are categories enriched in that tensor product. In this way, one is able see that the usual Gray tensor product and that of Crans, do so arise.

Moreover applying this lifting to the lax tensor products on $\PSh {\G}_{{\leq}n}$ associated to general $(n{+}1)$-operads (over $\Set$), one exhibits the structures definable by $(n{+}1)$-operads as enriched categories whose homs are some $n$-dimensional structure. In this way the globular approach is more closely related to some of the inductive approaches to higher category theory, such as that of Tamsamani \cite{Tam99}.

In section(\ref{ssec:A-infinity}) we describe two applications of the lifting theorem. In theorem(\ref{thm:A-infinity-app}) we construct a tensor product of $A_{\infty}$-algebras. As explained in \cite{Markl} the problem of providing such a tensor product is of relevance to string theory, and it proved resistant because of the negative result \cite{Markl} which shows that no ``genuine'' tensor product can exist. However this result does not rule out the existence of a \emph{lax} tensor product, which is what we were able to provide in theorem(\ref{thm:A-infinity-app}). It is possible to see the identification by Joyal and Street \cite{JS93}, of braided monoidal categories as monoidal categories with a multiplication as an instance of theorem(\ref{thm:A-infinity-app}). Another instance is our second application given in corollary(\ref{cor:coh-bm2c}), namely an analogous result to that of Joyal and Street but for braided monoidal 2-categories.

A weak $n$-category is an algebra of a \emph{contractible} $n$-operad{\footnotemark{\footnotetext{In this work we use the notion of contractibility given in \cite{Lei} rather than the original notion of \cite{Bat98}.}}}. In section(\ref{sec:contractibility}) we recall this notion, give an analogous notion of contractible lax monoidal structure and explain the canonical relationship between them.

In this paper we operate at a more abstract level than in much of the previous work on this subject. In particular, instead of studying monads on the category of $n$-globular sets, or even on presheaf categories, we work with monads defined on some category $\ca GV$ of enriched graphs. As our work shows, the main results and notions of higher category theory in the globular style can be given in this setting. So one could from the very beginning start not with $\Set$ as the category of $0$-categories, but with a nice enough $V$. For all the constructions to go through, such as that of $\ca T_{{\leq}n}$, the correspondence between monads and lax tensor products, their lifting theorem, as well as the very definition of weak $n$-category, it suffices to take $V$ to be a locally c-presentable category in the sense defined in section(\ref{ssec:lcpres}).

Proceeding this way one obtains then the theory of $n$-dimensional structures enriched in $V$. That is to say, the object of $n$-cells between any two $(n{-}1)$-cells of such a structure would be an object of $V$ rather than a mere set. Some alternative choices of $V$ which could perhaps be of interest are: (1) the ordinal $[1]=\{0<1\}$ (for the theory of locally ordered higher dimensional structures), (2) simplicial sets (to obtain a theory of higher dimensional structures which come together with a simplicial enrichment at the highest level), (3) the category of sheaves on a locally connected space, or more generally a locally connected Grothendieck topos, (4) the algebras of any $n$-operad or (5) the category of multicategories (symmetric or not). The point is, the theory as we have developed it is actually \emph{simpler} than before, and the generalisations mentioned here come at \emph{no} extra cost.

\section{Enriched graphs and lax monoidal categories}\label{sec:EG-LMC}

\subsection{Enriched graphs and the reduced suspension of spaces}\label{ssec:enriched-graphs}
Given a topological space $X$ and points $a$ and $b$ therein, one may define the topological space $X(a,b)$ of paths in $X$ from $a$ to $b$ at a high degree of generality. In recalling the details let us denote by $\Top$ a category of ``spaces'' which is complete, cocomplete and cartesian closed. We shall write $1$ for the terminal object. We shall furthermore assume that $\Top$ comes equipped with a bipointed object $I$ playing the role of the interval. A conventional choice for $\Top$ is the category of compactly generated Hausdorff spaces with its usual interval, although there are many other alternatives which would do just as well from the point of view of homotopy theory.

Let us denote by $\sigma{X}$ the reduced suspension of $X$, which can be defined as the pushout
\[ \xygraph{!{0;(1.5,0):(0,.667)::} {X{+}X}="tl" [r] {I{\times}X}="tr" [d] {\sigma{X}.}="br" [l] {1{+}1}="bl" "tl"(:"tr":"br",:"bl":"br") "br" [u(.3)l(.3)] (:@{-}[r(.15)],:@{-}[d(.15)])} \]
Writing $\Top_{\bullet}$ for the category of bipointed spaces, that is to say the coslice $1{+}1/\Top$, the above definition exhibits the reduced suspension construction as a functor
\[ \sigma : \Top \rightarrow \Top_{\bullet}. \]
In a sense this functor is the mother of homotopy theory -- applying it successively to the inclusion of the empty space into the point, one obtains the inclusions of the $(n{-}1)$-sphere into the $n$-disk for all $n \in \N$, and its right adjoint
\[ h : \Top_{\bullet} \rightarrow \Top \]
is the functor which sends the bipointed space $(a,X,b)$, to the space $X(a,b)$ of paths in $X$ from $a$ to $b$. This adjunction $\sigma \ladj h$ is easy to verify directly using the above elementary definition of $\sigma(X)$ as a pushout, and the pullback square
\[ \xygraph{!{0;(1.5,0):(0,.667)::} {X(a,b)}="tl" [r] {X^I}="tr" [d] {X^{1{+}1}}="br" [l] {1}="bl" "tl"(:"tr":"br"^{X^i},:"bl":"br"_-{(a,b)}) "tl" [d(.3)r(.3)] (:@{-}[l(.15)],:@{-}[u(.15)])} \]
where $i$ is the inclusion of the boundary of $I$. The collection of spaces $X(a,b)$ is our first example of an enriched graph in the sense of
\begin{definition}\label{def:enriched-graph}
Let $V$ be a category. A \emph{graph $X$ enriched in $V$} consists of an underlying set $X_0$ whose elements are called \emph{objects}, together with an object $X(a,b)$ of $V$ for each ordered pair $(a,b)$ of objects of $X$. The object $X(a,b)$ will sometimes be called the \emph{hom} from $a$ to $b$. A morphism $f:X{\rightarrow}Y$ of $V$-enriched graphs consists of a function $f_0:X_0{\rightarrow}Y_0$ together with a morphism $f_{a,b}:X(a,b){\rightarrow}Y(fa,fb)$ for each $(a,b)$. The category of $V$-graphs and their morphisms is denoted as $\ca GV$, and we denote by $\ca G$ the obvious 2-functor
\[ \begin{array}{lccr} {\ca G : \CAT \rightarrow \CAT} &&& {V \mapsto \ca GV} \end{array} \]
with object map as indicated.
\end{definition}
The 2-functor $\ca G$ is the mother of higher category theory in the globular style -- applying it successively to the inclusion of the empty category into the point (ie the terminal category), one obtains the inclusion of the category of $(n{-}1)$-globular sets into the category of $n$-globular sets. In the case $n>0$ this is the inclusion with object map
\[ \begin{array}{lcr}
{\xygraph{{X_0}="l" [r] {...}="m" [r] {X_{n{-}1}}="r" "r":@<1ex>"m":@<1ex>"l" "r":@<-1ex>"m":@<-1ex>"l"}} & \mapsto &
{\xygraph{{X_0}="l" [r] {...}="m" [r] {X_{n{-}1}}="r" [r] {\emptyset}="rr" "rr":@<1ex>"r":@<1ex>"m":@<1ex>"l" "rr":@<-1ex>"r":@<-1ex>"m":@<-1ex>"l"}}
\end{array} \]
and when $n{=}0$ this is the functor $1{\rightarrow}\Set$ which picks out the empty set. Thus there is exactly one $(-1)$-globular set which may be identified with the empty set.

It is often better to think of $\ca G$ as taking values in $\CAT/\Set$. By applying the endofunctor $\ca G$ to the unique functor $V{\rightarrow}1$ for each $V$, produces the forgetful functor $\ca GV{\rightarrow}\Set$ which sends an enriched graph to its underlying set of objects. This manifestation
\[ \ca G : \CAT \rightarrow \CAT/\Set \]
has a left adjoint which we shall denote as $(-)_{\bullet}$ for reasons that are about to become clear. The functor $(-)_{\bullet}$ is a variation of the Grothendieck construction. To a given functor $f:A{\rightarrow}\Set$ it associates the category $A_{\bullet}$ with objects triples $(x,a,y)$ where $a$ is an object of $A$, and $(x,y)$ is an ordered pair of objects of $fa$. Maps are just maps in $A$ which preserve these base points in the obvious sense.

It is interesting to look at the unit and counit of this 2-adjunction. Given a category $V$, $(\ca GV)_{\bullet}$ is the category of bipointed enriched graphs in $V$. The counit $\varepsilon_V:(\ca GV)_{\bullet}{\rightarrow}V$ sends $(a,X,b)$ to the hom $X(a,b)$. When $V$ has an initial object $\varepsilon_V$ has a left adjoint given by $X \mapsto (X)$. Given a functor $f:A{\rightarrow}\Set$ the unit $\eta_f:A{\rightarrow}\ca G(A_{\bullet})$ sends $a \in A$ to the enriched graph whose objects are elements of $fa$, and the hom $a(x,y)$ is given by the bipointed object $(x,a,y)$.

Consider the case where $0 \in A$ and $f$ is the representable $f=A(0,-)$. Then $A_{\bullet}$ may be regarded as the category of endo-cospans of the object $0$, that is to say the category of diagrams
\[ 0 \rightarrow a \leftarrow 0 \]
and a point of $a \in A$ is now just a map $0{\rightarrow}a$. When $A$ is also cocomplete one can compute a left adjoint to $\eta_A$. To do this note that a graph $X$ enriched in $A_{\bullet}$ gives rise to a functor
\[ \overline{X} : X^{(2)}_0 \rightarrow A \]
where $X_0$ is the set of objects of $X$. For any set $Z$, $Z^{(2)}$ is defined as the following category. It has two kinds of objects: an object being either an element of $Z$, or an ordered pair of elements of $Z$. There are two kinds of non-identity maps
\[ x \rightarrow (x,y) \leftarrow y \]
where $(x,y)$ is an ordered pair from $Z$, and $Z^{(2)}$ is free on the graph just described. A more conceptual way to see this category is as the category of elements of the graph
\[ \xygraph{{Z{\times}Z}="l" [r] {Z}="r" "l":@<1ex>"r" "l":@<-1ex>"r"} \]
where the source and target maps are the product projections, as a presheaf on the category
\[ \xygraph{{\G_{\leq{1}}} [r(.75)] {=} [r(1.25)] {\xybox{\xygraph{0 [r] 1 "0":@<1ex>"1":@<1ex>@{<-}"0"}}} *\frm{-}}  \]
and so there is a discrete fibration $Z^{(2)}{\rightarrow}\G_{{\leq}1}$. The functor $\overline{X}$ sends singletons to $0 \in A$, and a pair $(x,y)$ to the head of the hom $X(x,y)$. The arrow map of $\overline{X}$ encodes the bipointings of the homs. One may then easily verify
\begin{proposition}
Let $0 \in A$, $f=A(0,-)$ and $A$ be cocomplete. Then $\eta_f$ has left adjoint given on objects by $X \mapsto \colim(\overline{X})$.
\end{proposition}
In the exposition thus far we have focussed on building an analogy between the reduced suspension of a space and the graphs enriched in a category. Now we shall bring these constructions together. As we have seen already to each space $X$ one can associate a canonical topologically enriched graph whose homs are the path spaces of $X$. Denoting this enriched graph as $PX$, the assignment $X \mapsto PX$ is the object map of the composite right adjoint in
\[ \xygraph{!{0;(2,0):} {\Top}="l" [r] {\ca G(\Top_{\bullet})}="m" [r] {\ca G\Top}="r" "l":@<-1.5ex>"m"_-{\eta}|{}="b1":@<-1.5ex>"l"|{}="t1" "t1":@{}"b1"|{\perp} "m":@<-1.5ex>"r"_-{\ca Gh}|{}="b2":@<-1.5ex>"m"_-{\ca G\sigma}|{}="t2" "t2":@{}"b2"|{\perp}}. \]
As explained by Cheng \cite{ChengCompOp}, this functor $P = \ca G(h)\eta$ is a key ingredient of the Trimble definition of weak $n$-category.

\subsection{Some exactness properties of the endofunctor $\ca G$}\label{ssec:G-exactness}
The important properties of $\ca G$ are apparent because of the close connection between $\ca G$ and the $\Fam$ construction. A very mild reformulation of the notion of $V$-graph is the following: a $V$-graph $X$ consists of a set $X_0$ together with an $(X_0{\times}X_0)$-indexed family of objects of $V$. Together with the analogous reformulation of the maps of $\ca GV$, this means that we have a pullback square
\[ \PbSq {\ca GV} {\Fam{V}} {\Set} {\Set} {(-)_0=\ca Gt_V} {} {\Fam(t_V)} {(-)^2} \]
in $\CAT$, and thus a cartesian 2-natural transformation $\ca G{\implies}\Fam$. From \cite{Fam2fun} theorem(7.4) we conclude
\begin{proposition}\label{prop:GFam2fun}
$\ca G$ is a familial 2-functor.
\end{proposition}
\noindent In particular it follows from the theory of \cite{Fam2fun} that $\ca G$ preserves conical connected limits as well as all the notions of ``Grothendieck fibration'' which one can define internal to a finitely complete 2-category. Moreover the obstruction maps for comma objects are right adjoints. See \cite{Fam2fun} for more details on this part of 2-category theory. We shall not use these observations very much in what follows. More important for us is
\begin{lemma}\label{lem:G-EM-object}
$\ca G$ preserves Eilenberg-Moore objects.
\end{lemma}
\noindent Given a monad $T$ on a category $V$, we shall write $V^T$ for the category of $T$-algebras and morphisms thereof, and $U^T:V^T{\rightarrow}V$ for the forgetful functor. We shall denote a typical object of $V^T$ as a pair $(X,x)$, where $X$ is the underlying object in $V$ and $x:TX{\rightarrow}X$ is the $T$-algebra structure. From \cite{Str72} the 2-cell $TU^T{\implies}U^T$, whose component at $(X,x)$ is $x$ itself has a universal property exhibiting $V^T$ as a kind of 2-categorical limit called an \emph{Eilenberg-Moore object}. See \cite{Str72} or \cite{LS00} for more details on this general notion. The direct proof that for any monad $T$ on a category $V$, the obstruction map $\ca G(V^T){\rightarrow}\ca G(V)^{\ca G(T)}$ is an isomorphism comes down to the obvious fact that for any $V$-graph $B$, a $\ca GT$-algebra structure on $B$ is the same thing as a $T$-algebra structure on the homs of $B$, and similarly for algebra morphisms.

\subsection{Multitensors}\label{ssec:LMC} Let us recall the notions of lax monoidal category and category enriched therein from \cite{EnHopI}. For a category $V$, the free strict monoidal category $MV$ on $V$ has a very simple description. An object of $MV$ is a finite sequence $(Z_1,...,Z_n)$ of objects of $V$. A map is a sequence of maps of $V$ -- there are no maps between sequences of objects of different lengths. The unit $\eta_V:V{\rightarrow}MV$ of the 2-monad $M$ is the inclusion of sequences of length $1$. The multiplication $\mu_V:M^2V{\rightarrow}MV$ is given by concatenation.

A \emph{lax monoidal category} is a lax algebra for the 2-monad $M$. Explicitly it consists of an underlying category $V$, a functor $E:MV{\rightarrow}V$, and maps
\[ \begin{array}{lcr} {u_Z:Z \rightarrow E(Z)} && {\sigma_{Z_{ij}}:\opE\limits_i\opE\limits_j Z_{ij} \rightarrow \opE\limits_{ij} Z_{ij}} \end{array} \]
for all $Z$, $Z_{ij}$ from $V$ which are natural in their arguments, and such that
\[ \xy (0,0); (10,0):
(0,0)*{\xybox{\xymatrix @C=1em {{\opE\limits_iZ_i} \ar[r]^-{u\opE\limits_i} \ar[d]_{1} \save \POS?(.4)="domeq" \restore & {E_1\opE\limits_iZ_i} \ar[dl]^{\sigma} \save \POS?(.4)="codeq" \restore \\ {\opE\limits_iZ_i} \POS "domeq"; "codeq" **@{}; ?*{=}}}};
(4,0)*{\xybox{\xymatrix @C=1.5em {{\opE\limits_i\opE\limits_j\opE\limits_kZ_{ijk}} \ar[r]^-{{\sigma}\opE\limits_k} \ar[d]_{\opE\limits_i\sigma} \save \POS?="domeq" \restore & {\opE\limits_{ij}\opE\limits_kZ_{ijk}} \ar[d]^{\sigma} \save \POS?="codeq" \restore \\
{\opE\limits_i\opE\limits_{jk}Z_{ijk}} \ar[r]_-{\sigma} & {\opE\limits_{ijk}Z_{ijk}} \POS "domeq"; "codeq" **@{}; ?*{=}}}};
(8,0)*{\xybox{{\xymatrix @C=1em {{\opE\limits_iE_1Z_i} \ar[dr]_{\sigma} \save \POS?(.4)="domeq" \restore & {\opE\limits_iZ_i} \ar[d]^{1} \save \POS?(.4)="codeq" \restore  \ar[l]_-{\opE\limits_iu} \\ & {\opE\limits_iZ_i}} \POS "domeq"; "codeq" **@{}; ?*{=}}}}
\endxy \]
in $V$. As in \cite{EnHopI} we use either of the expressions
\[ \begin{array}{lcr} {\opE\limits_{1{\leq}i{\leq}n} X_i} && {\opE\limits_i X_i} \end{array} \]
as a convenient yet precise short-hand for $E(X_1,...,X_n)$, and we refer to the endofunctor of $V$ obtained by observing the effect of $E$ on singleton sequences as $E_1$. The data $(E,u,\sigma)$ is called a \emph{multitensor} on $V$, and $u$ and $\sigma$ are referred to as the unit and substitution of the multitensor respectively.

Given a multitensor $(E,u,\sigma)$ on $V$, a \emph{category enriched in $E$} consists of $X \in \ca GV$ together with maps
\[ \kappa_{x_i} : \opE\limits_i X(x_{i-1},x_i) \rightarrow X(x_0,x_n) \]
for all $n \in \N$ and sequences $(x_0,...,x_n)$ of objects of $X$, such that
\[ \xygraph{
{\xybox{\xygraph{!{0;(2,0):(0,.6)::}
{X(x_0,x_1)} (:[r]{E_1X(x_0,x_1)}^-{u} :[d]{X(x_0,x_1)}="bot"^{\kappa}, :"bot"_{\id})}}}
[r(5)][d(.15)]
{\xybox{\xygraph{!{0;(2.75,0):(0,.5)::}
{\opE\limits_i\opE\limits_jX(x_{(ij)-1},x_{ij})} (:[r]{\opE\limits_{ij}X(x_{(ij)-1},x_{ij})}^-{\sigma} :[d]{X(x_0,x_{mn_m})}="bot"^{\kappa},:[d]{\opE\limits_iX(x_{(i1)-1},x_{in_i})}_{\opE\limits_i\kappa} :"bot"_-{\kappa})}}}} \]
commute, where $1{\leq}i{\leq}m$, $1{\leq}j{\leq}n_i$ and $x_{(11)-1}{=}x_0$. Since a choice of $i$ and $j$ references an element of the ordinal $n_{\bullet}$, the predecessor $(ij){-}1$ of the pair $(ij)$ is well-defined when $i$ and $j$ are not both $1$. With the obvious notion of $E$-functor (see \cite{EnHopI}), one has a category $\Enrich E$ of $E$-categories and $E$-functors together with a forgetful functor
\[ U^E : \Enrich E \rightarrow \ca GV. \]
The notation we use makes transparent the analogy between multitensors and monads, and categories enriched in multitensors and algebras for a monad. In particular the unit and subtitution for $E$ provide $E_1$ with the unit and multiplication of a monad structure. Moreover, any object of the form $\opE\limits_i Z_i$ is canonically an $E_1$-algebra, as is the hom of any $E$-category, and the substitution maps of $E$ are $E_1$-algebra morphisms (see \cite{EnHopI} lemma(2.7)). Thus in a sense, any multitensor $(E,u,\sigma)$ on a category $V$ is aspiring to be a multitensor on the category $V^{E_1}$ of $E_1$-algebras, but of course there is no meaningful way to regard $u$ as living in $V^{E_1}$, except in the boring situation when $E_1$ is the identity monad, that is, when $u$ is an identity natural transformation. The multitensors with $u$ the identity are called \emph{normal}.
\begin{definition}\label{def:lift}
Let $(E,u,\sigma)$ be a multitensor on a category $V$. A \emph{lift} of $(E,u,\sigma)$ is a normal multitensor $(E',\id,\sigma')$ on $V^{E_1}$ together with an isomorphism $\Enrich E \iso \Enrich {E'}$ which commutes with the forgetful functors into $\ca G(V^{E_1})$.
\end{definition}
\noindent In \cite{EnHopI} we explained how to associate normalised $(n+1)$-operads and $n$-multitensors, which are multitensors on the category of $n$-globular sets. In the present paper we shall explain why any $n$-multitensor has a canonical lift.

\section{Multitensors from monads}\label{sec:Monad->Mult}

\subsection{Monads over $\Set$}\label{ssec:DefNMonad}
At an abstract level much of this paper is about the interplay between the theory of monads on categories of enriched graphs, and the theory of multitensors. It is time to be more precise about which monads on $\ca GV$ we are interested in.
\begin{definition}\label{def:NMnd}
Let $V$ be a category. A monad \emph{over $\Set$} on $\ca GV$ is a monad on
\[ (-)_0 : \ca GV \rightarrow \Set \]
in the 2-category $\CAT/\Set$.
\end{definition}
\noindent That is, a monad $(T,\eta,\mu)$ on $\ca GV$ is over $\Set$ when the functor $T$ doesn't affect the object sets, in other words $TX_0=X_0$ for all $X \in \ca GV$ and similarly for maps, and moreover the components of $\eta$ and $\mu$ are identities on objects. In this section we will describe how such a monad, in the case where $V$ has an initial object denoted as $\emptyset$, induces a multitensor on $V$ denoted $(\overline{T},\overline{\eta},\overline{\mu})$.

Let us describe this multitensor explicitly. First we note that $\emptyset \in V$ enables us to regard any sequence of objects $(Z_1,...,Z_n)$ of $V$ as a $V$-graph. The object set is $\{0,...,n\}$,
\[ (Z_1,...,Z_n)(i-1,i) = Z_i \]
for $1{\leq}i{\leq}n$, and all the other homs are equal to $\emptyset$. Then we define
\begin{equation}\label{eq:Tbar} \overline{T}(Z_1,...,Z_n) := T(Z_1,...,Z_n)(0,n) \end{equation}
and the unit as
\begin{equation}\label{eq:etabar} \overline{\eta}_Z := \{\eta_{(Z)}\}_{0,1}. \end{equation}
Before defining $\overline{\mu}$ we require some preliminaries. Given objects $Z_i$ of $V$ where $1{\leq}i{\leq}n$, and $1{\leq}a{\leq}b{\leq}n$ denote by
\[ s_{a,b} : (Z_a,...,Z_b) \rightarrow (Z_1,...,Z_n) \]
the obvious subsequence inclusion in $\ca GV$: the object map preserves successor and $0 \mapsto (a-1)$, and the hom maps are identities. Now given objects $Z_{ij}$ of $V$ where $1{\leq}i{\leq}k$ and $1{\leq}j{\leq}n_i$, one has a map
\[ \tilde{\tau}_{Z_{ij}} : (\Tbar\limits_{1{\leq}j{\leq}n_1}Z_{1j},...,\Tbar\limits_{1{\leq}j{\leq}n_k}Z_{kj}) \rightarrow T(Z_{11},......,Z_{kn_k}) \]
given on objects by $0 \mapsto 0$ and $i \mapsto (i,n_i)$ for $1{\leq}i{\leq}k$, and the hom map between $(i-1)$ and $i$ is $Ts_{i1,in_i}$. With these definitions in hand we can now define the components of $\overline{\mu}$ as
\begin{equation}\label{eq:mubar} \xymatrix{{\Tbar\limits_i\Tbar\limits_jZ_{ij}} \ar[rr]^-{\{T\tilde{\tau}\}_{0,k}} && {T(Z_{11},...,Z_{kn_k})(0,n_{\bullet})} \ar[rr]^-{\mu_{0,n_{\bullet}}} && {\Tbar\limits_{ij}Z_{ij}}}. \end{equation}
From now until the end of (\ref{ssec:Tbar}) we shall be occupied with the proof of
\begin{theorem}\label{thm:Tbar}
Let $V$ be a category with an initial object $\emptyset$ and $(T,\eta,\mu)$ be a monad over $\Set$ on $\ca GV$. Then $(\overline{T},\overline{\eta},\overline{\mu})$ as defined in (\ref{eq:Tbar})-(\ref{eq:mubar}) defines a multitensor on $V$.
\end{theorem}
\noindent In principle one could supply a proof of this result immediately by just slogging through a direct verification of the axioms. Instead we shall take a more conceptual approach, and along the way encounter various ideas that are of independent interest. For most of the time we will assume a little more of $V$: that it has finite coproducts, to enable our more conceptual approach. In the end though, we will see that only the initial object is necessary.

We break up the construction of $(\overline{T},\overline{\eta},\overline{\mu})$ into three steps. First in (\ref{ssec:NMonad->MonMonad}), we describe how $\ca GV_{\bullet}$ acquires a monoidal structure and $(T,\eta,\mu)$ induces a monoidal monad $(T_{\bullet},\eta_{\bullet},\mu_{\bullet})$ on $\ca GV_{\bullet}$. Then we see that this monoidal monad induces a lax monoidal structure on $\ca GV_{\bullet}$, which in turn can be transferred across an adjunction to obtain $(\overline{T},\overline{\eta},\overline{\mu})$. These last two steps are very general: they work at the level of the theory of lax algebras for an arbitrary 2-monad (which in our case is $M$). So in (\ref{ssec:laxalg-const1}) and (\ref{ssec:laxalg-const2}) we describe these general constructions, and in (\ref{ssec:Tbar}) we finish the proof of theorem(\ref{thm:Tbar}). Finally in (\ref{ssec:path-like}) we present a condition on $T$ which ensures that $T$-algebras and $\overline{T}$-categories may be identified.

\subsection{Monoidal monads from monads over $\Set$}\label{ssec:NMonad->MonMonad} From now until just before the end of (\ref{ssec:Tbar}) we shall assume that $V$ has finite coproducts. We now describe some consequences of this. First, the functor $(-)_0:\ca GV{\rightarrow}\Set$ which sends an enriched graph to its set of objects becomes representable. We shall denote by $0$ the $V$-graph which represents $(-)_0$. It has one object and its unique hom is $\emptyset$.

The second consequence is that $\ca GV_{\bullet}$ inherits a natural monoidal structure and any normalised monad $(T,\eta,\mu)$ on $\ca GV$ can then be regarded as a monoidal monad $(T_{\bullet},\eta_{\bullet},\mu_{\bullet})$. The explanation for this begins with the observation that the representability of the underlying set functor $(-)_0$ enables a useful reformulation of the category $\ca GV_{\bullet}$ as the category of endocospans of $0$ as in section(\ref{ssec:enriched-graphs}). The usefulness of this is that such cospans can be composed, thus endowing $\ca GV_{\bullet}$ with a canonical monoidal structure.

The presence of $\emptyset$ in $V$ enables one to compute coproducts in $\ca GV$. The coproduct $X$ of a family $(X_i:i{\in}I)$ of $V$-graphs has object set given as the disjoint union of the object sets of the $X_i$, $X(x,y)=X_i(x,y)$ when $x$ and $y$ are objects of $X_i$, and all the other homs are $\emptyset$. With finite coproducts available one can also compute pushouts under $0$, that is the pushout $P$ of maps
\[ \xymatrix{X & 0 \ar[r]^-{y} \ar[l]_-{x} & Y} \]
in $\ca GV$ is described as follows. The object set of $P$ is the disjoint union of the object sets of $X$ and $Y$ modulo the identification of $x$ and $y$, and let us write $z$ for this special element of $P$. The homs of $P$ are inherited from $X$ and $Y$ in almost the same way as for coproducts. That is if $a$ and $b$ are either both objects of $X$ or both objects of $Y$ and they are not both $z$, then their hom $P(a,b)$ is taken as in $X$ or $Y$. The hom $P(z,z)$ is the coproduct $X(x,x)+Y(y,y)$. Otherwise this hom is $\emptyset$. Note that in the special case where the homs $X(x,x)$ and $Y(y,y)$ are both $\emptyset$, one only requires the initial object in $V$ to compute this pushout.

Given a sequence
\[ ((a_1,X_1,b_1),...,(a_n,X_n,b_n)) \]
of doubly-pointed $V$-graphs, one defines their \emph{join}
\[ (a_1,(a_1,X_1,b_1)*...*(a_n,X_n,b_n),b_n) \]
where $(a_1,X_1,b_1)*...*(a_n,X_n,b_n)$ is the colimit of
\begin{equation}\label{eq:join} \xymatrix{{X_1} & 0 \ar[l]_-{b_1} \ar[r]^-{a_2} & {...} & 0 \ar[l]_-{b_{n-1}} \ar[r]^-{a_n} & {X_n}} \end{equation}
in $\ca GV$ which can be formed via iterated pushouts under $0$.
This defines a monoidal structure on $\ca GV_{\bullet}$ whose tensor product we denote as $*:M\ca GV_{\bullet}{\rightarrow}\ca GV_{\bullet}$.

Given a functor $T:\ca GV{\rightarrow}\ca GW$ over $\Set$, one has a functor
\[ \begin{array}{lcr} {T_{\bullet} : \ca GV_{\bullet} \rightarrow \ca GW_{\bullet}} && {(a,X,b) \mapsto (a,TX,b)} \end{array} \]
whose object map is indicated on the right in the previous display. When both $V$ and $W$ have an initial object, one defines for each sequence of doubly-pointed $V$-graphs a map
\[ \tau_{X_i} : T_{\bullet}(a_1,X_1,b_1)*...*T_{\bullet}(a_n,X_n,b_n) \rightarrow T_{\bullet}((a_1,X_1,b_1)*...*(a_n,X_n,b_n)) \]
as follows. Write
\[ \begin{array}{c} {c_i:X_i \rightarrow (a_1,X_1,b_1)*...*(a_n,X_n,b_n)} \\ {c'_i:TX_i \rightarrow (a_1,TX_1,b_1)*...*(a_n,TX_n,b_n)} \end{array} \]
for the components of the colimit cocones (\ref{eq:join}). Using the unique map $0{\rightarrow}T0$, there is a unique map $\tau_{X_i}$ such that $Tc_i=\tau_{X_i}{c'_i}$. By the unique characterisation of these maps, they assemble together to provide the coherence 2-cell
\[ \LaxSq {M\ca GV_{\bullet}} {\ca GV_{\bullet}} {M\ca GW_{\bullet}} {\ca GW_{\bullet}} {MT_{\bullet}} {*} {T_{\bullet}} {*} {\tau} \]
for a lax monoidal functor, and for
\[ \xymatrix{{\ca GU} \ar[r]^-{S} & {\ca GV} \ar[r]^-{T} & {\ca GW}} \]
over $\Set$ one has $T_{\bullet}S_{\bullet}=(TS)_{\bullet}$ as monoidal functors. Moreover any natural transformation $\phi:S{\rightarrow}T$ over $\Set$ defines a monoidal natural tranformation $\phi_{\bullet}:S_{\bullet}{\rightarrow}T_{\bullet}$. In fact, denoting by $\ca N$ the 2-category whose objects are categories with initial objects, a 1-cell $T:V{\rightarrow}W$ in $\ca N$ is a functor $T:\ca GV{\rightarrow}\ca GW$ over $\Set$, and a 2-cell between these is a natural tranformation also over $\Set$, we have defined a 2-functor
\[ J : \ca N \rightarrow \LaxAlg M. \]
Applying $J$ to monads gives
\begin{proposition}\label{prop:NMnd->MonMnd}
If $V$ has finite coproducts and $(T,\eta,\mu)$ is a monad over $\Set$ on $\ca GV$, then $(T_{\bullet},\eta_{\bullet},\mu_{\bullet})$ is a monoidal monad on $\ca GV_{\bullet}$, whose monoidal structure is given by pushout-composition of cospans.
\end{proposition}
Finally we note that the functor $\varepsilon:\ca GV_{\bullet}{\rightarrow}V$ has a left adjoint $L$ which we shall now describe. Given $Z \in V$ the underlying $V$-graph of $LZ$, which we shall denote as $(Z)$, has object set $\{0,1\}$ and the distinguished pair is $(0,1)$. As for the homs, $(Z)(0,1)$ is just $Z$ itself, and all the other homs are $\emptyset$. Formally it is the composite
\[ \xymatrix{{MV} \ar[r]^-{ML} & {M\ca GV_{\bullet}} \ar[r]^-{*} & {\ca GV_{\bullet}} \ar[r]^-{U} & {\ca GV}} \]
where $U$ is the obvious forgetful functor, which enables us to view a sequence of objects of $V$ as a $V$-graph, as in (\ref{ssec:DefNMonad}) above.

\subsection{A general lax algebra construction}\label{ssec:laxalg-const1}
Now and in (\ref{ssec:laxalg-const2}) let $(S,\eta,\mu)$ be a 2-monad on a 2-category $\ca K$.
Suppose that we are given a monad \emph{in} $\LaxAlg S$. Let us write $(V,E,\iota,\sigma)$ for the underlying lax $S$-algebra, $(F,\phi)$ for the lax $S$-algebra endomorphism of $V$, and $i$ and $m$ for the unit and multiplication respectively. Then one obtains another lax $S$-algebra structure on $V$ with one cell part given as the composite
\[ \xymatrix{{SV} \ar[r]^-{E} & V \ar[r]^-{F} & V} \]
and the 2-cell data as follows
\[ \xygraph{{\xybox{\xygraph{!{0;(1.5,0):(0,0.7)::}
{V}="tl" [r] {SV}="tr" [d] {V}="br" [l] {V}="bl"
"tl":"tr"^-{\eta}:"br"^{E}:"bl"^{F} "tl":"br"|{1} "tl":"bl"_{1}
"tl" [d(.35)r(.6)] :@{=>}[r(.2)]^{\iota} "tl" [d(.7)r(.2)] :@{=>}[r(.2)]^{i}}}}
[r(3)]
{\xybox{\xygraph{{S^2V}="l1" [d] {SV}="l2" [d] {SV}="l3" [r] {V}="m" [r] {V}="r3" [u] {V}="r2" [u] {SV}="r1"
"l1":"l2"_{SE}:"l3"_{SF}:"m"_-{E}:"r3"_-{F} "l1":"r1"^-{\mu}:"r2"^{E}:"r3"^{F} "l2":"r2"^-{E}:"m"_{F}
"l1" [d(.5)r(.85)] :@{=>}[r(.3)]^{\sigma} "l2" [d(.5)r(.5)] :@{=>}[r(.3)]^{\phi} "m" [u(.2)r(.5)] :@{=>}[r(.3)]^{m}}}}} \]
The verification of the lax algebra axioms is an easy exercise that is left to the reader.

\subsection{Another general lax algebra construction}\label{ssec:laxalg-const2}
Now suppose we are given a lax S-algebra $(V,E,\iota,\sigma)$
together with an adjunction
\[ \xygraph{!{0;(1.2,0):} V [r] W "V":@<-1.2ex>"W"_-{R}:@<-1.2ex>"V"_-{L} "V":@{}"W"|-{\perp}} \]
with unit $u$ and counit $c$. One can then induce a lax $S$-algebra structure on $W$. The one-cell part is given as the composite
\[ \xymatrix{{SW} \ar[r]^-{SL} & {SV} \ar[r]^-{E} & V \ar[r]^-{R} & W} \]
and the 2-cell data as follows
\[ \xygraph{{\xybox{
\xygraph{{W}="l" [r] {V}="m" [r] {SV}="r" [l(.5)d] {V}="br" [l] {W}="bl" "m" [u] {SW}="t"
"l":"t"^-{\eta}:"r"^-{SL}:"br"^-{E}|{}="E":"bl"^-{R} "l":"m"^-{L}:"br"_{1}|{}="oneV" "m":"r"^-{\eta} "l":"bl"_{1}|{}="oneW"
"m":@{}"t"|(.4)*{=}
"oneW":@{}"oneV"|(.35){}="d1"|(.65){}="c1" "d1":@{=>}"c1"^-{u}
"oneV":@{}"E"|(.2){}="d2"|(.8){}="c2" "d2":@{=>}"c2"^-{\iota}}}}
[r(4)] {\xybox{
\xygraph{!{0;(1,0):(0,.8)::}
{S^2W}="l1" [dl] {S^2V}="l2" [d] {SV}="l3" [d] {SW}="l4" [r] {SV}="ml" [r] {V}="mr" [r] {W}="r4" [u] {V}="r3" [u] {SV}="r2" [ul] {SW}="r1"
"l1":"l2"_{S^2L}:"l3"_{SE}:"l4"_{SR}:"ml"_{SL}:"mr"_{E}:"r4"_{R} "l1":"r1"^-{\mu}:"r2"^{SL}:"r3"^{E}:"r4"^{R} "l3":"ml"^{1}|{}="oneMV" "mr":"r3"^{1}|{}="oneV" "l2":"r2"^-{\mu}
"l1" [d(2)r(.35)] :@{=>}[r(.3)]^{\sigma} "l3" [d(.7)r(.15)] :@{=>}[r(.25)]^{Sc}
"r4" [u(.25)l(.35)] {=} "l1" [d(.5)r(.5)] {=}}}}} \]
The reader will easily verify that the lax $S$-algebra axioms for $W$ follow from those of $V$ and the triangle identities of the adjunction.

\subsection{Multitensors from monads over $\Set$}\label{ssec:Tbar} Let us now put together (\ref{ssec:NMonad->MonMonad})-(\ref{ssec:laxalg-const2}).
Given a monad $(T,\eta,\mu)$ on $\ca GV$ over $\Set$ such that $V$ has finite coproducts, we obtained the monoidal monad $(T_{\bullet},\eta_{\bullet},\mu_{\bullet})$ on $\ca GV_{\bullet}$ in proposition(\ref{prop:NMnd->MonMnd}). In other words $(T_{\bullet},\eta_{\bullet},\mu_{\bullet})$ is a monad in $\LaxAlg M$. Applying (\ref{ssec:laxalg-const1}) for $S=M$ gives us a multitensor on $\ca GV_{\bullet}$, and then applying (\ref{ssec:laxalg-const1}) to this last multitensor and the adjunction
\[ \xygraph{!{0;(1.2,0):} {\ca GV_{\bullet}}="V" [r] {V}="W" "V":@<-1.2ex>"W"_-{\varepsilon}:@<-1.2ex>"V"_-{L} "V":@{}"W"|-{\perp}} \]
gives us a multitensor on $V$ which we denote as $(\overline{T},\overline{\eta},\overline{\mu})$.

We shall now unpack this multitensor to see that it does indeed agree with that of theorem(\ref{thm:Tbar}). The one cell part $\overline{T}$ is the composite
\[ \xymatrix{{MV} \ar[r]^-{ML} & {M\ca GV_{\bullet}} \ar[r]^-{*} & {\ca GV_{\bullet}} \ar[r]^-{T_{\bullet}} & {\ca GV_{\bullet}} \ar[r]^-{\varepsilon} & V} \]
which agrees with equation(\ref{eq:Tbar}) and one also easily reconciles equation(\ref{eq:etabar}) for the unit. As for the substitution unpacking the $\overline{\mu}$ of our conceptual approach gives the following composite
\[ \xygraph{!{0;(1,0):(0,.85)::}
{M^2V}="p11" [r] {MV}="p12" [dr] {M\ca GV_{\bullet}}="p22" [l(3)] {M^2\ca GV_{\bullet}}="p21" [dl] {M\ca GV_{\bullet}}="p31" [d] {M\ca GV_{\bullet}}="p41" [d] {MV}="p51" [r(1.25)d] {M\ca GV_{\bullet}}="p52" [r(1.25)d] {\ca GV_{\bullet}}="p53" [r(1.25)u] {\ca GV_{\bullet}}="p54" [r(1.25)u] {V}="p55" [u] {\ca GV_{\bullet}}="p42" [u] {\ca GV_{\bullet}}="p32"
"p11":"p12"^-{\mu}|{}="mv1" "p21":"p22"^-{\mu}|{}="mv2" "p31":"p32"^-{*}|{}="mv3" "p51":"p52"_-{ML}:"p53"_-{*}:"p54"_-{T_{\bullet}}:"p55"_-{\varepsilon}
"p11":"p21"_-{M^2L}:"p31"_-{M*}:"p41"_-{MT_{\bullet}}:"p51"_-{M\varepsilon} "p12":"p22"^{ML}:"p32"^{*}:"p42"^{
T_{\bullet}}:"p55"^{\varepsilon}
"p41":@/^{2pc}/"p52"^{1}|{}="mh1" "p32":@/_{2.5pc}/"p53"_{T_{\bullet}}|{}="mh2" "p42":@/_{.5pc}/"p54"_{1}|{}="mh3"
"mv1":@{}"mv2"|*{=}:@{}"mv3"|*{\iso}
"p51":@{}"mh1"|{}="mid1":@{}"mh2"|{}="mid2":@{}"mh3"|{}="mid3":@{}"p55"|*{=}
"mid1" [l(.15)] :@{=>}[r(.3)]^{Mc} "mid2" [l(.15)] :@{=>}[r(.3)]^{\tau} "mid3" [l(.15)] :@{=>}[r(.3)]^{\mu_{\bullet}}} \]
which we shall now unpack further. The counit $c$ of the adjunction $L \ladj \varepsilon$ is described as follows. For a given doubly-pointed $\ca V$-graph $(a,X,b)$, the corresponding counit component
\[ c_{(a,X,b)} : (0,(X(a,b)),1) \rightarrow (a,X,b) \]
is specified by insisting that the hom map between $0$ and $1$ is $1_{X(a,b)}$. Define $\tilde{\mu}$ by
\[ \xygraph{!{0;(4,0):(0,.35)::}
{T(\Tbar\limits_jX_{1j},...,\Tbar\limits_jX_{kj})}="bl" [u] {T\{(0,T(X_{11},...,X_{1n_1}),n_1)*...*(0,T(X_{k1},...,X_{kn_k}),n_k)\}}="tl" [d(.4)r] {T^2(X_{11},......,X_{kn_k})}="tr" [d(.6)] {T(X_{11},......,X_{kn_k})}="br"
"bl"(:"tl"^{T\{c*...*c\}}:"tr"_-{T\tau}:"br"^{\mu},:"br"_-{\tilde{\mu}}) "bl":@{}"tr"|*{=}} \]
and then $\overline{\mu}$ is the effect of $\tilde{\mu}$ on the hom between $0$ and $k$. But one may easily verify that the composite $\tau\{c*...*c\}$ is just $\tilde{\tau}$ described in (\ref{ssec:DefNMonad}).

This completes the proof of theorem(\ref{thm:Tbar}) for the case where $V$ has finite coproducts. The general case is obtained by observing that only the joins of doubly-pointed $V$-graphs $(a,X,b)$ such that $X(a,a)$ and $X(b,b)=\emptyset$ are actually used in the construction of the multitensor and in its axioms, and these only require an initial object in $V$.

\subsection{Path-like monads}\label{ssec:path-like}
We shall now give a condition on a normalised monad $T$ which ensures that categories enriched in $\overline{T}$ are the same thing as $T$-algebras. Let $(a,X,b) \in \ca GV_{\bullet}$ and consider a sequence $x=(x_0,...,x_n)$ of objects of $X$ such that $x_0=a$ and $x_n=b$. Define the $V$-graph
\[ x^*X := (X(x_0,x_1),X(x_1,x_2),...,X(x_{n-1},x_n)) \]
and we have a map
\[ \overline{x} : (0,x^*X,n) \rightarrow (a,X,b) \]
given on objects by $i \mapsto x_i$ for $0{\leq}i{\leq}n$, and the effect on the hom between $(i-1)$ and $i$ is the identity for $1{\leq}i{\leq}n$.
\begin{definition}\label{def:path-like}
Let $V$ be a category with an initial object and $W$ be a category with all small coproducts. A functor $T:\ca GV{\rightarrow}\ca GW$ over $\Set$ is \emph{path-like} when for all $(a,X,b) \in \ca GV_{\bullet}$, the maps
\[ T\overline{x}_{0,n} : Tx^*X(0,n) \rightarrow TX(a,b) \]
for all $n \in \N$ and sequences $x=(x_0,...,x_n)$ such that $x_0=a$ and $x_n=b$, form a coproduct cocone in $W$. A normalised monad $(T,\eta,\mu)$ on $\ca GW$ is \emph{path-like} when $T$ is path-like in the sense just defined. 
\end{definition}
\begin{example}\label{ex:cat-monad-path-like}
Let $V=W=\Set$ and $T$ be the free category endofunctor of $\Graph$. For any graph $X$ and $a,b \in X_0$, the hom $TX(a,b)$ is by definition the set of paths in $X$ from $a$ to $b$. Each path determines a sequence $x=(x_0,...,x_n)$ of objects of $X$ such that $x_0=a$ and $x_n=b$, by reading off the objects of $X$ as they are visited by the given path. Conversely for a sequence $x=(x_0,...,x_n)$ of objects of $X$ such that $x_0=a$ and $x_n=b$, $T\overline{x}$ identifies the elements $Tx^*X(0,n)$ with those paths in $X$ from $a$ to $b$ whose associated sequence is $x$. Thus $T$ is path-like.
\end{example}
\begin{proposition}\label{prop:pl-alg<->cat}
Let $V$ have small coproducts and $(T,\eta,\mu)$ be a path-like monad on $\ca GV$ over $\Set$. Then $\ca G(V)^T \iso \Enrich {\overline{T}}$.
\end{proposition}
\begin{proof}
Let $X$ be a $V$-graph. To give an identity on objects map $a:TX{\rightarrow}X$ is to give maps $a_{y,z}:TX(y,z){\rightarrow}X(y,z)$. By path-likeness these amount to giving for each $n \in \N$ and $x=(x_0,...,x_n)$ such that $x_0=y$ and $x_n=z$, a map
\[ a_x : \Tbar\limits_iX(x_{i-1},x_i) \rightarrow X(y,z) \]
since $\Tbar\limits_iX(x_{i-1},x_i)=Tx^*X(0,n)$, that is $a_x=a_{y,z}T\overline{x}_{0,n}$. When $n=1$, for a given $y,z \in X_0$, $x$ can only be the sequence $(y,z)$. The naturality square for $\eta$ at $\overline{x}$ implies that $\{\eta_X\}_{y,z}=T\overline{x}_{0,1}\{\eta_{(X(y,z))}\}_{0,1}$, and the definition of $\overline{(\,\,)}$ says that $\{\eta_{(X(y,z))}\}_{0,1}=\overline{\eta}_{X(y,z)}$. Thus to say that a map $a:TX{\rightarrow}X$ satisfies the unit law of a $T$-algebra is to say that $a$ is the identity on objects and that the $a_x$ described above satisfy the unit axioms of a $\overline{T}$-category.

To say that $a$ satisfies the associative law is to say that for all $y,z \in X_0$,
\begin{equation}\label{eq:assoc1} \xymatrix{{T^2X(y,z)} \ar[r]^-{\{\mu_X\}_{y,z}} \ar[d]_{Tx_{y,z}} & {TX(y,z)} \ar[d]^{a_{y,z}} \\ {TX(y,z)} \ar[r]_-{a_{y,z}} & {X(y,z)}} \end{equation}
commutes. Given $x=(x_0,...,x_n)$ from $X$ with $x_0=y$ and $x_n=z$, and $w=(w_0,...,w_k)$ from $Tx^*X$ with $w_0=0$ and $w_k=n$, consider the composite map
\begin{equation}\label{eq:copr} \xymatrix{{Tw^*Tx^*X(0,k)} \ar[r]^-{T\overline{w}_{0,k}} & {T^2x^*X(0,n)} \ar[r]^-{T\overline{x}_{0,n}} & {T^2X(y,z)}} \end{equation}
and note that by path-likeness, and since the coproduct of coproducts is a coproduct, all such maps for $x$ and $w$ such that $x_0=y$ and $x_n=z$ form a coproduct cocone. Precomposing (\ref{eq:assoc1}) with (\ref{eq:copr}) gives the commutativity of
\begin{equation}\label{eq:assoc2} \xymatrix{{\Tbar\limits_i\Tbar\limits_jX(x_{ij-1},x_{ij})} \ar[r]^-{\overline{\mu}} \ar[d]_{\Tbar\limits_ia} & {\Tbar\limits_{ij}X(x_{ij-1},x_{ij})} \ar[d]^{a_x} \\ {\Tbar\limits_iX(x_{w_{i-1}},x_{w_i})} \ar[r]_-{a_w} & {X(y,z)}} \end{equation}
and conversely by the previous sentence if these squares commute for all $x$ and $w$, then one recovers the commutativity of (\ref{eq:assoc1}). This completes the description of the object part of $\ca G(V)^T \iso \Enrich {\overline{T}}$.

Let $(X,a)$ and $(X',a')$ be $T$-algebras and $F:X{\rightarrow}X'$ be a $V$-graph morphism. To say that $F$ is a $T$-algebra map is a condition on the maps $F_{y,z}:X(y,z){\rightarrow}X'(Fy,Fz)$ for all $y,z \in X_0$, and one uses path-likeness in the obvious way to see that this is equivalent to saying that the $F_{y,z}$ are the hom maps for a $\overline{T}$-functor.
\end{proof}
\noindent The proof of proposition(\ref{prop:pl-alg<->cat}) is not new: exactly the same argument was used in the second half of the proof of theorem(7.6) of \cite{EnHopI}, although in that case the setting was far less general. The real novelty is the generality of definition(\ref{def:path-like}) which is crucial for section(\ref{sec:lift-mult}).

\section{Monads from multitensors}\label{sec:Mult->Monad}

\subsection{Distributive multitensors}\label{ssec:DistMult}
The general way of obtaining a monad from a multitensor, which is the topic of this section, applies to multitensors which conform to
\begin{definition}\label{def:Dist-Mult}
Let $V$ be a category with small coproducts. Then a multitensor $(E,u,\sigma)$ is \emph{distributive} when the functor $E$ preserves coproducts in each variable. That is to say, for each $n \in \N$, the functor
\[ E_n : V^n \rightarrow V \]
obtained by observing $E$'s
effect on sequences of length $n$, preserves coproducts in each of its $n$ variables.
\end{definition}
\noindent The first step in associating a monad to a distributive multitensor is to identify the bicategory $\Dist$ which has the property that monads in $\Dist$ are exactly distributive multitensors in the sense of definition(\ref{def:Dist-Mult}). There is also a useful reformulation of the notion of monad on $\ca GV$ over $\Set$: as a monad in another 2-category $\Kl {\ca G_{\bullet}}$ where $\ca G_{\bullet}$ denotes the 2-comonad on $\CAT$ induced by the adjunction $(-)_{\bullet} \ladj \ca G$ of section(\ref{ssec:enriched-graphs}). Our monad-from-multitensor construction is then achieved by means of a pseudo functor (homomorphism of bicategories)
\[ \Gamma : \Dist \rightarrow \Kl {\ca G_{\bullet}} \]
which as a pseudo functor sends monads to monads. In fact for a distributive multitensor $E$, $\Gamma{E}$ is path-like and so $\Enrich E \iso \ca G(V)^{\Gamma{E}}$.

\subsection{The bicategory $\Dist$}\label{ssec:Dist}
The objects of $\Dist$ are categories with coproducts. A morphism $E:V{\rightarrow}W$ in $\Dist$ is a functor $E:MV{\rightarrow}W$ which preserves coproducts in each variable. A 2-cell between $E$ and $E'$ is simply a natural transformation between these functors $MV{\rightarrow}W$. Vertical composition of 2-cells is as for natural transformations.

The horizontal composite of $F:MV{\rightarrow}W$ and $E:MW{\rightarrow}X$, denoted $E \comp F$, is defined as a left kan extension
\[ \LaxSq {M^2V} {MV} {MW} {X} {MF} {\mu} {E \comp F} {E} {l_{E,F}} \]
of $EM(F)$ along $\mu$. Computing this explicitly gives the formula
\[ (E \comp F)(X_1,...,X_n) = \coprod_{n_1+...+n_k=n} \opE\limits_i\opF\limits_j X_{ij} \]
where $1{\leq}i{\leq}k$ and $1{\leq}j{\leq}n_i$ on the right hand side of this formula, and we denote by
\[ \xymatrix{{\opE\limits_{i}\opF\limits_{j}X_{ij}} \ar[r]^-{\opc\limits_{ij}} & {\opEoF\limits_{ij}X_{ij}}} \]
and also by $c_{(n_1,...,n_k)}$, the corresponding coproduct inclusion. The definition of horizontal composition is clearly functorial with respect to vertical composition of 2-cells.
\begin{proposition}\label{prop:Dist}
$\Dist$ is a bicategory and a monad in $\Dist$ is exactly a distributive multitensor.
\end{proposition}
\begin{proof}
It remains to identify the coherences and check the coherence axioms. The notation we have used here for the coproduct inclusions matches that used in section(3) of \cite{EnHopI}. The proof of the first part of proposition(3.3) of \cite{EnHopI} interpretted in our present more general setting, is the proof that $\Dist$ is a bicategory. The characterisation of distributive multitensors as monads in $\Dist$ is immediate from the definitions and our abstract definition of horizontal composition in terms of kan extensions.
\end{proof}
\noindent What in \cite{EnHopI} was called $\Dist(V)$ is here the hom $\Dist(V,V)$.

\subsection{Reformulating monads over $\Set$}\label{ssec:KlDefNMonad}
Because of the adjunction $(-)_{\bullet} \ladj \ca G$ and the definition of the comonad $\ca G_{\bullet}$, a normalised monad on $\ca GV$ is the same thing as a monad on $V$ in the Kleisli 2-category{\footnotemark{\footnotetext{In $\CAT$-enriched sense.}}} $\Kl {\ca G_{\bullet}}$ of the (2-)comonad $\ca G_{\bullet}$. The reason this is sometimes useful is that it expresses how our monads over $\Set$ only involve information at the level of homs. The validity of this reformulation is most plainly seen by realising that the factorisation of $\ca G$ as an identity on objects 2-functor followed by a 2-fully-faithful 2-functor, can be realised as
\[ \xymatrix{{\CAT} \ar[r]^-{R} & {\Kl {\ca G_{\bullet}}} \ar[r]^-{J} & {\CAT/\Set}} \]
where $R$ is the right adjoint part of the Kleisli adjunction for $\ca G_{\bullet}$. Thus applying $J$ sends a monad in $\Kl {\ca G_{\bullet}}$ on $V$ to a monad on $\ca GV$ over $\Set$, and the definition of $\ca G$ and the 2-fully-faithfulness of $J$ ensures that any such monad arises uniquely in this way.

\subsection{The pseudo-functor $\Gamma$}\label{ssec:Gamma}
It is the identity on objects. Given a distributive $E:MV{\rightarrow}W$, $\Gamma{E}$ is defined as a left kan extension
\[ \LaxSq {M\ca G_{\bullet}V} {\ca G_{\bullet}V} {MV} {W} {M\varepsilon} {*} {\Gamma{E}} {E} {\gamma_E} \]
of $EM(\varepsilon)$ along $*$. We shall now explain why this left kan extension exists in general, give a more explicit formula for $\Gamma{E}$ in corollary(\ref{cor:explicit-gamma}), and then with this in hand it will become clear why $\Gamma$ is a pseudo-functor. Of course one could just define $\Gamma{E}$ via the formula in corollary(\ref{cor:explicit-gamma}). We chose instead to give the above more abstract definition, because it will enable us to attain a more natural understanding of why $\Gamma$ produces \emph{path-like} monads from distributive multitensors in (\ref{ssec:DistMult->Monad}).

Conceptually, the reason why the left kan extension $\gamma_E$ involves only coproducts in $W$ is that $*$ is a local left adjoint. A functor $F:A{\rightarrow}B$ is a \emph{local left adjoint} when $\op F$ is a local right adjoint. This is equivalent to asking that for all $a \in A$ the induced functor $F^a:a/A{\rightarrow}Fa/B$ between coslices is a left adjoint. We shall explain in lemma(\ref{lem:spancomp-lra}) why and how pullback-composition of spans can be seen as a local right adjoint, and the statement that $*$ is a local left adjoint is just the dual of this because $*$ is defined as pushout-composition of cospans in $\ca GV$. In lemma(\ref{lem:lla-ple}) we give a formula for computing the left kan extension along a local left adjoint, and this will then be applied to give our promised more explicit description of $\Gamma$'s one-cell map.
\begin{lemma}\label{lem:spancomp-lra}
Let $A$ be a category with finite products and let $a$ be an object of $A$ such that $A/a$ also has finite products. Then the one-cell part
\[ M(A/a{\times}a) \rightarrow A/a{\times}a \]
of the monoidal structure on $A/a{\times}a$ given by span composition is a local right adjoint.
\end{lemma}
\begin{proof}
Note that for all $n \in \N$ the slices $A/(a^n)$, where $a^n$ is the $n$-fold cartesian product of $a$, have finite products, so the statement of the lemma makes sense and all the limits we mention in this proof exist.
In general a functor out of a coproduct of categories is a local right adjoint iff its composite with each coproduct inclusion is a local right adjoint. Thus it suffices to show that $n$-fold composition of spans
\[ (A/a{\times}a)^n \rightarrow A/a{\times}a \]
is a local right adjoint for all $n \in \N$. The case $n=0$ may be exhibited as a composite
\[ \xymatrix{1 \ar[r]^-{1_a} & {A/a} \ar[r]^-{\Delta_!} & {A/a{\times}a}} \]
of local right adjoints and thus is a local right adjoint. The case $n=1$ may be regarded as the identity. It suffices to verify the case $n=2$ because with this in hand an easy induction will give the general case. The pullback composite of spans as shown on the left
\[ \xygraph{*\xybox{\xygraph{!{0;(.65,0):}
d ([dl] b ([dl] {a}="a1", [dr] {a}="a2"), [dr] c [dr] {a}="a3")
"d":"b"_{p}:"a1"_{w} "d":"c"^{q}:"a3"^{z} "b":"a2"^{x} "c":"a2"_{y}
"d" [d(.3)l(.3)] [d(.1)r(.1)]:@{-}[d(.2)r(.2)]:@{-}[u(.2)r(.2)] }}
[r(4)]
*\xybox{\xygraph{!{0;(1.2,0):(0,.7)::}
d [r] {b{\times}c}="bc" [d] {a^4}="a4" [l] {a^3}="a3" [d] {a^2}="a2"
"d"(:@{.>}"a3"_{(wp,xp=yq,zq)}(:@{.>}"a2"_{(\pi_1,\pi_3)},:"a4"_-{a{\times}\Delta{\times}a}), :"bc"^-{(p,q)}:"a4"^{(w,x){\times}(y,z)}
"d" [d(.3)] [r(.1)]:@{-}[r(.2)]:@{-}[u(.2)]}}} \]
may be constructed as the dotted composite shown on the right in the previous display, and so binary composition of endospans of $a$ is encoded by the composite functor
\[ \xymatrix @C=4em {{(A/a{\times}a)^2} \ar[r]^-{\textnormal{prod}_{(a^2,a^2)}} & {A/(a^4)} \ar[r]^-{(a{\times}\Delta{\times}a)^*} & {A/(a^3)} \ar[r]^-{(\pi_1,\pi_3)_!} & {A/a{\times}a}} \]
where $\textnormal{prod}$ denotes the (right adjoint) functor $A{\times}A{\rightarrow}A$ which sends a pair to its cartesian product (and $\textnormal{prod}_{(a^2,a^2)}$ is the slice of $\textnormal{prod}$ over the pair $(a^2,a^2)$). The constituent functors of this last composite are clearly all local right adjoints.
\end{proof}
\noindent Recall \cite{WebGen} \cite{Fam2fun} that local right adjoints can be characterised in terms of generic factorisations. The dual characterisation is as follows. A functor $F:A{\rightarrow}B$ is a local left adjoint iff for all $b \in B$ the components of the comma category $F/b$ have terminal objects. A map $f:Fa{\rightarrow}b$ which is terminal in its component of $F/b$ is said to be \emph{cogeneric}, and by definition any $f$ can be factored as
\[ \xymatrix{{Fa} \ar[r]^-{Fh} & {Fc} \ar[r]^-{g} & b} \]
where $g$ is cogeneric. This factorisation is unique up to unique isomorphism and is called the \emph{cogeneric factorisation} of $f$.
\begin{lemma}\label{lem:lla-ple}
Let $F:A{\rightarrow}B$ be a local left adjoint. Suppose that $A$ has a set $C$ of connected components, each of which has an initial object, and that $B$ is locally small. For a functor $G:A{\rightarrow}D$ where $D$ has coproducts, the left kan extension $L:B{\rightarrow}D$ of $G$ along $F$ exists and is given by the formula
\[ Lb = \coprod\limits_{c{\in}C} \coprod\limits_{f:F0_c{\rightarrow}b} Ga_f \]
where $0_c$ denotes the initial object of the component $c \in C$, and
\[ \xymatrix{{F0_c} \ar[r]^-{Fh_f} & {Fa_f} \ar[r]^-{g_f} & b} \]
is a chosen cogeneric factorisation of $f:F0_c{\rightarrow}b$.
\end{lemma}
\begin{proof}
By the general formula for computing left kan extensions as colimits, it suffices to identify the given formula with the colimit of
\[ \xymatrix{{F/b} \ar[r]^-{p} & A \ar[r]^-{G} & D} \]
where $p$ is the obvious projection. This follows since the components of $F/b$ are indexed by pairs $(c,f)$, where $c \in C$ and $f:F0_c{\rightarrow}b$, and the component of $F/b$ corresponding to $(c,f)$ has terminal object given by $g_f$, which is mapped by $p$ to $a_f$.
\end{proof}
\noindent In the case of $*:M\ca GV_{\bullet}{\rightarrow}\ca GV_{\bullet}$ note the initial objects of the components of $M\ca GV_{\bullet}$ are of the form
\[ 0_n = (\underbrace{(0,(\emptyset),1),...,(0,(\emptyset),1)}_{n}) \]
and a map $*0_n{\rightarrow}(a,X,b)$ is a map
\[ (0,(\emptyset,...,\emptyset),n) \rightarrow (a,X,b) \]
which amounts to a sequence of elements of $X$ of length $(n+1)$ starting at $a$ and finishing at $b$. The reader will easily verify that
\[ \xymatrix{{(\emptyset,...,\emptyset)} \ar[r] & {x^*X} \ar[r]^-{\overline{x}} & {X}} \]
is a cogeneric factorisation of the map associated to $x=(x_0,...,x_n)$. So we have given a conceptual explanation of $x^*X$ and $\overline{x}$ which were used in (\ref{ssec:path-like}), as well as completed the proof of
\begin{corollary}\label{cor:explicit-gamma}
For $V$ and $W$ with coproducts and $E:MV{\rightarrow}W$, the defining left kan extension of $\Gamma{E}$ exists and we have the formula
\[ \Gamma{E}(a,X,b) = \coprod_{\begin{array}{c} {(x_0,...,x_n)} \\ {x_0=a, \, x_n=b}\end{array}} \opE\limits_i X(x_{i-1},x_i) \]
\end{corollary}
We shall identify $\Gamma$ with the composite $J\Gamma$, which amounts to identifying $\Gamma{E}:\ca GV_{\bullet}{\rightarrow}W$ with its mate $\ca GV{\rightarrow}\ca GW$ by the adjunction $(-)_{\bullet} \ladj \ca G$. Then for $X \in \ca GV$ and $a,b \in X_0$, we have for each $n \in \N$ and for each sequence $x=(x_0,...,x_n)$ such that $x_0=a$ and $x_n=b$ a coproduct inclusion
\[ c_x : \opE\limits_i X(x_{i-1},x_i) \rightarrow \Gamma{E}X(a,b) \]
by corollary(\ref{cor:explicit-gamma}). The components of the coherence natural transformations for $\Gamma$ will be identities on objects. From the definition of the unit $I:X{\rightarrow}X$ in $\Dist$, one has that for $a,b \in X_0$ the coproduct inclusion
\[ c_{(a,b)} : X(a,b) \rightarrow \Gamma{I}X(a,b) \]
is an isomorphism. Thus we have an isomorphism $\gamma_0:1{\rightarrow}\Gamma{I}$. Given distributive $F:MU{\rightarrow}V$ and $E:MV{\rightarrow}W$, $X \in \ca GU$ and $a,b \in X_0$, we define the hom maps of
\[ \gamma_{2,X} : \Gamma(E)\Gamma(F)X \rightarrow \Gamma(E \comp F)X \]
for $(a,b) \in X_0$, as the unique isomorphism such that for all $x_{ij} \in X_0$ where $1{\leq}i{\leq}k$ and $1{\leq}j{\leq}n_i$, $x_{11}=a$ and $x_{kn_k}=b$, the diagram
\[ \xygraph{!{(0,0);(0,1.5):(0,2)::}
{\opE\limits_i\opF\limits_jX(x_{ij-1},x_{ij})}="a" ([d]{(\opEoF\limits_{ij})X(x_{ij-1},x_{ij})}="c" [d][r(.5)]{\Gamma(E{\comp}F)X(a,b)}="e",
[r]{\opE\limits_i\Gamma(F)X(x_{i-1},x_i)}="b" [d]{\Gamma(E)\Gamma(F)X(a,b)}="d")
"a":"b"^-{\opE\limits_i\opc\limits_j}:"d"^-{c_{x_{i\bullet}}}:"e"^-{\gamma_2}
"a":"c"_-{\opc\limits_{ij}}:"e"_-{c_{x_{ij}}}} \]
in $W$ commutes, where $x_0=x_{11}$ and $x_i=x_{in_i}$ for $i{>}0$. We have selected the notation so as to match up with the development of \cite{EnHopI} section(4), and the proof of the first part of proposition(4.1) of \emph{loc. sit.} interpretted in the present context gives
\begin{proposition}\label{prop:gamma-psfunctor}
The coherences $(\gamma_0,\gamma_2)$ just defined make $\Gamma$ into a pseudo-functor.
\end{proposition}
\noindent Given a monad $T$ on $\ca GV$ over $\Set$, and a set $Z$, one obtains by restriction a monad $T_Z$ on the category $\ca GV_Z$ of $V$-graphs with fixed object set $Z$. Let us write $\Gamma^{\textnormal{old}}$ for the functor labelled as $\Gamma$ in \cite{EnHopI}. Then for a given distributive multitensor $E$, our present $\Gamma$ and $\Gamma^{\textnormal{old}}$ are related by the formula
\[ \Gamma^{\textnormal{old}}(E) = \Gamma(E)_1 \]
where the $1$ on the right hand side of this equation indicates a singleton. In other words we have just given the ``many-objects version'' of the theory presented in \cite{EnHopI} section(4).

\subsection{Algebras and enriched categories}\label{ssec:DistMult->Monad}
Having just established the machinery to convert distributive multitensors on $V$ to monads on $\ca GV$ over $\Set$, we shall now relate the enriched categories to the algebras. This involves two things: seeing that the normalised monads constructed from distributive multitensors are path-like, and understanding the relationship between $\Gamma$ and the construction of section(\ref{sec:Monad->Mult}) of multitensors from monads.
\begin{lemma}\label{lem:gamma-path-like}
Let $V$ and $W$ have coproducts and $E:MV{\rightarrow}W$ preserve coproducts in each variable. Then $\Gamma{E}:\ca GV{\rightarrow}\ca GW$ is path-like.
\end{lemma}
\begin{proof}
The condition that $T:\ca GV{\rightarrow}\ca GW$ is path-like can be expressed more 2-categorically.
\[ \xygraph{
*{\xybox{\xygraph{!{0;(1.5,0):(0,.67)::}
{*/(a,X,b)}="tl" [r] {1}="tr" [d] {\ca GV_{\bullet}}="mr" [l] {M\ca GV_{\bullet}}="ml" [d] {}="bl" [r] {W}="br"
"tl"(:"ml"_{p}:"mr"^-{*},:"tr":"mr"^{(a,X,b)}:"br"^{T})
"tl" [d(.5)r(.35)] :@{=>}[r(.3)]^{\lambda}}}}
[r(4)d(.1)]
*{\xybox{\xygraph{!{0;(1.5,0):(0,.67)::}
{*/(a,X,b)}="tl" [r] {1}="tr" [d] {\ca GV_{\bullet}}="mr" [l] {M\ca GV_{\bullet}}="ml" [d] {MV}="bl" [r] {W}="br"
"tl"(:"ml"_{p}(:"bl"_{M\varepsilon}:"br"_-{E},:"mr"^-{*}),:"tr":"mr"^{(a,X,b)}:"br"^{\Gamma{E}})
"tl" [d(.5)r(.35)] :@{=>}[r(.3)]^{\lambda} "ml" [d(.5)r(.35)] :@{=>}[r(.3)]^{\gamma_E}}}}} \]
Writing $\lambda$ for 2-cell part of the comma object, $T$ is path-like iff the 2-cell on the left exhibits $T(a,X,b)$ as a colimit. To see this recall from (\ref{ssec:Gamma}) that the set of components of $*/(a,X,b)$ may be regarded as the set of sequences of objects of $X_0$ starting at $a$ and finishing at $b$, that each of these components has a terminal object, and that $\overline{x}:x^*X{\rightarrow}X$ is terminal in the component corresponding to the sequence $x=(x_0,...,x_n)$. The situation for a given $E$ is depicted on the right in the previous display, and by definition this composite 2-cell is a colimit. Thus it suffices to show that the component of the 2-cell $\gamma_E$ at $p\overline{x}$ is invertible for all sequences $x=(x_0,...,x_n)$ from $X_0$. The component of $\gamma_E$ at a general
\[ ((c_1,Y_1,d_1),...,(c_n,Y_n,d_n)) \in M\ca G_{\bullet}V \]
is the coproduct inclusion
\[ \xymatrix{{\opE\limits_i Y_i(c_i,d_i)} \ar[rr]^-{c_w} && {\coprod\limits_{(z_0,...,z_m)} \opE\limits_i Y(z_{i-1},z_i)}} \]
corresponding to the sequence
\[ w = (c_1,d_1=c_2,...,d_{n-1}=c_n,d_n) \]
where $Y=(c_1,Y_1,d_1)*...*(c_n,Y_n,d_n)$. In the case of $p\overline{x} \in M\ca GV_{\bullet}$, $Y=x^*X$, and for summands corresponding to sequences $z$ different from $w$, we will have $Y(z_{i-1},z_i)=\emptyset$ for some $i$. By the distributivity of $E$ those summands will be $\emptyset$, whence $c_w$ will be invertible.
\end{proof}
\noindent Given a distributive multitensor $(E,u,\sigma)$ note that one can apply $\Gamma$ to it and then $\overline{(-)}$ to the result. One has
\[ \overline{\Gamma{E}}(Z_1,...,Z_n) = \Gamma{E}(Z_1,...,Z_n)(0,n) = \coprod_{a_0,...,a_m} \opE_i(Z_1,...,Z_n)(a_{i-1},a_i) \]
where the $a_i$ in the sum are elements of $\{0,...,n\}$ and $a_0=0$ and $a_m=n$. Unless the sequence $(a_0,...,a_n)$ is just an in-order list $(0,...,n)$ of the elements of $\{0,...,n\}$, at least one of the homs $(Z_1,...,Z_n)(a_{i-1},a_i)$ must be $\emptyset$ making that summand $\emptyset$ by the distributivity of $E$. Thus the coproduct inclusion
\[ c_{(0,...,n)} : \opE_iZ_i \rightarrow \overline{\Gamma{E}}(Z_1,...,Z_n) \]
is invertible. Moreover using the explicit description of the multitensor $\overline{\Gamma{E}}$ one may verify that this isomorphism is compatible with the units and substitutions, and so we have
\begin{lemma}\label{lem:gamma-bar}
If $(E,u,\sigma)$ is a distributive multitensor on a category $V$ with coproducts, then one has an isomorphism $E \iso \overline{\Gamma{E}}$ of multitensors.
\end{lemma}
\noindent Together with lemma(\ref{lem:gamma-path-like}) and proposition(\ref{prop:pl-alg<->cat}) this implies
\begin{corollary}\label{cor:gamma-alg-ecat}
If $(E,u,\sigma)$ is a distributive multitensor on a category $V$ with coproducts, then one has $\Enrich E \iso \ca G(V)^{\Gamma{E}}$ commuting with the forgetful functors into $\ca GV$.
\end{corollary}

\subsection{A conceptual view of path-likeness}\label{ssec:pl-conceptual}
We now describe the sense in which $\Gamma$ and $\overline{(-)}$ are adjoint. First let us note that equation(\ref{eq:Tbar}) defining the construction $\overline{(-)}$ in section(\ref{ssec:DefNMonad}) may be seen as providing functors
\[ \overline{(-)}_{V,W} : \CAT/\Set(\ca GV,\ca GW) \rightarrow \CAT(MV,W) \]
for all $V,W$ in $\CAT$. We have abused notation slightly by denoting by $\ca GV$ (resp. $\ca GW$) the category of $V$-enriched graphs together with its forgetful functor into $\Set$. In (\ref{ssec:DefNMonad}) we considered only the case $V=W$ and when $T:\ca GV{\rightarrow}\ca GV$ is part of a monad, but equation(\ref{eq:Tbar}) obviously makes sense in this more general context. In order to relate this with $\Gamma$ we make
\begin{definition}\label{def:ndist}
Let $V$ and $W$ have coproducts. A functor $T:{\ca GV}{\rightarrow}{\ca GW}$ over $\Set$ is \emph{distributive} when $\overline{T}:MV{\rightarrow}W$ preserves coproducts in each variable. We denote by $\NDist(V,W)$ the full subcategory of $\CAT/\Set(\ca GV,\ca GW)$ consisting of the distributive functors from $\ca GV$ to $\ca GW$.
\end{definition}
\begin{proposition}\label{prop:pl-adjoint-char}
Let $V$ and $W$ be categories with coproducts. Then we have an adjunction
\[ \xygraph{!{0;(3,0):}
{\Dist(V,W)}="d" [r] {\NDist(V,W)}="nd" "d":@<1ex>"nd"^-{\Gamma_{V,W}}|{}="t":@<1ex>"d"^-{\overline{(-)}_{V,W}}|{}="b"
"t":@{}"b"|-{\perp}} \]
whose unit is invertible. A distributive $T:{\ca GV}{\rightarrow}{\ca GW}$ is in the image of $\Gamma_{V,W}$ iff it is path-like.
\end{proposition}
\begin{proof}
By lemma(\ref{lem:gamma-bar}) applying $\Gamma$ does indeed produce a distributive functor, so $\Gamma_{V,W}$ is well-defined and one has an isomorphism of $\overline{(-)}_{V,W}\Gamma_{V,W}$ with the identity. For $X \in \ca GV$ and $a,b \in X_0$ one has
\[ \Gamma\overline{T}X(a,b) = \coprod_{a=x_0,...,x_n=b} Tx^*X(0,n) \rightarrow TX(a,b) \]
induced by the hom-maps of $T\overline{x}$, giving $\varepsilon_T:\Gamma\overline{T}{\rightarrow}T$ natural in $T$, and so by \cite{Web2top} lemma(2.6) to establish the adjunction it suffices to show that $\varepsilon$ is inverted by $\overline{(-)}$. To this end note that when $X=(Z_1,...,Z_m)$ for $Z_i \in V$, the above summands are non-initial iff the sequence $x_0,...,x_n$ is the sequence $(0,1,...,m)$, by the distributivity of $T$. The characterisation of path-likeness now follows too, since this condition on a given $T$ is by definition the same as the invertibility of $\varepsilon_T$.
\end{proof}
\noindent An immediate consequence of proposition(\ref{prop:pl-adjoint-char}) and proposition(\ref{prop:GammaE-basic})(\ref{GEb1}) below is the following result. A direct proof is also quite straight forward and is left as an exercise.
\begin{corollary}\label{cor:pl->copr-pres}
Let $V$ and $W$ be categories with coproducts. If $T:{\ca GV}{\rightarrow}{\ca GW}$ over $\Set$ is distributive and path-like, then it preserves coproducts.
\end{corollary}
%

\section{Categorical properties preserved by $\Gamma$}\label{sec:reexpress}

\subsection{}\label{ssec:intro-reexpress}
Let us now regard $\Gamma$ as a pseudo-functor
\[ \Gamma : \Dist \rightarrow \CAT/\Set. \]
That is to say, we take for granted the inclusion of $\Kl {\ca G_{\bullet}}$ in $\CAT/\Set$. In this section we shall give a systematic account of the categorical properties that $\Gamma$ preserves. The machinery we are developing gives an elegant inductive description of the monads $\ca T_{\leq{n}}$ for strict $n$-categories, provides explanations of some of their key properties, and gives a shorter account of the central result of \cite{EnHopI} on the equivalence between $n$-multitensors and $(n+1)$-operads.

\subsection{A review of some categorical notions}\label{ssec:lcpres}
Let $\lambda$ be a regular cardinal. An object $C$ in a category $V$ is \emph{connected} when the representable functor $V(C,-)$ preserves coproducts, and $C$ is \emph{$\lambda$-presentable} when $V(C,-)$ preserves $\lambda$-filtered colimits. The object $C$ is said to be \emph{small} when it is $\lambda$-presentable for some regular cardinal $\lambda$.

A category $V$ is \emph{extensive} when it has coproducts and for all families $(X_i:i \in I)$ of objects of $V$, the functor
\[ \begin{array}{lccr} {\coprod : \prod\limits_i (V/X_i) \rightarrow V/(\coprod\limits_i X_i)}
&&& {(f_i:Y_i{\rightarrow}X_i) \mapsto \coprod\limits_if_i:\coprod\limits_iY_i{\rightarrow}\coprod\limits_iX_i} \end{array} \]
is an equivalence of categories. A more elementary characterisation is that $V$ is extensive iff it has coproducts, pullbacks along coproduct coprojections and given a family of commutative squares
\[ \xymatrix{{X_i} \ar[r]^-{c_i} \ar[d]_{f_i} & X \ar[d]^{f} \\ {Y_i} \ar[r]_-{d_i} & Y} \]
where $i \in I$ such that the $d_i$ form a coproduct cocone, the $c_i$ form a coproduct cocone iff these squares are all pullbacks. It follows that coproducts are disjoint{\footnotemark{\footnotetext{Meaning that coproduct coprojections are mono and the pullback of different coprojections is initial.}}} and the initial object of $V$ is strict{\footnotemark{\footnotetext{Meaning that any map into the initial object is an isomorphism.}}}. Another sufficient condition for extensivity is provided by
\begin{lemma}\label{lem:easy-ext}
If a category $V$ has disjoint coproducts and a strict initial object, and every $X \in V$ is a coproduct of connected objects, then $V$ is extensive.
\end{lemma}
\noindent The proof is left as an easy exercise. Note this condition is not necessary: there are many extensive categories whose objects don't decompose into coproducts of connected objects, for example, take the topos of sheaves on a space which is not locally connected. A category is \emph{lextensive} when it is extensive and has finite limits. There are many examples of lextensive categories: Grothendieck toposes, the category of algebras of any higher operad and the category of topological spaces and continuous maps are all lextensive.

Denoting the terminal object of a lextensive category $V$ by $1$, the representable $V(1,-)$ has a left exact left adjoint
\[ (-) \cdot 1 : \Set \rightarrow V \]
which sends a set $Z$ to the copower $Z{\cdot}1$. This functor enables one to express coproduct decompositions of objects of $V$, \emph{internal to $V$} because to give a map
\[ f : X \rightarrow I \cdot 1 \]
in $V$ is the same thing as giving an $I$-indexed coproduct decomposition of $X$. The lextensive categories in which every object decomposes into a sum of connected objects are characterised by the following well-known result.
\begin{proposition}\label{prop:lext-decompose-charn}
For a lextensive category $V$ the following statements are equivalent:
\begin{enumerate}
\item  Every $X \in V$ can be expressed as a coproduct of connected objects.\label{copcon}
\item  The functor $(-) \cdot 1$ has a left adjoint.\label{pi0}
\end{enumerate}
\end{proposition}
\noindent The most common instance of this is when $V$ is a Grothendieck topos. The toposes $V$ satisfying the equivalent conditions of proposition(\ref{prop:lext-decompose-charn}) are said to be \emph{locally connected}. This terminology is reasonable since for a topological space $X$, one has that $X$ is locally connected as a space iff its associated topos of sheaves is locally connected in this sense.

A set $\ca D$ of objects of $V$ is a \emph{strong generator} when for all maps $f:X{\rightarrow}Y$, if
\[ V(D,f) : V(D,X) \rightarrow V(D,Y) \]
is bijective for all $D \in \ca D$ then $f$ is an isomorphism. A locally small category $V$ is \emph{locally $\lambda$-presentable} when it is cocomplete and has a strong generator consisting of small objects. Finally recall that a functor is \emph{accessible} when it preserves $\lambda$-filtered colimits for some regular cardinal $\lambda$.

The theory of locally presentable categories is one of the high points of classical category theory, and this notion admits many alternative characterisations \cite{GabUlm} \cite{MP} \cite{AR94}. For instance locally presentable categories are exactly those categories which are the $\Set$-valued models for a limit sketch. Grothendieck toposes are locally presentable because each covering sieve in a Grothendieck topology on a category $\C$ gives rise to a cone in $\op {\C}$, and a sheaf is exactly a functor $\op {\C}{\rightarrow}\Set$ which sends these cones to limit cones in $\Set$. That is to say a Grothendieck topos can be seen as the models of a limit sketch which one obtains in an obvious way from any site which presents it. Just as locally presentable categories generalise Grothendieck toposes, the following notion generalises locally connected Grothendieck toposes.
\begin{definition}\label{def:lcc}
A locally small category $V$ is \emph{locally c-presentable} when it is cocomplete and has a strong generator consisting of small connected objects.
\end{definition}
\noindent Just as locally presentable categories have many alternative characterisations we have the following result for locally c-presentable categories. Its proof is obtained by applying the general results of \cite{ABLR} in the case of the doctrine for $\lambda$-small connected categories, which is ``sound'' (see \cite{ABLR}),
and proposition(\ref{prop:lext-decompose-charn}).
\begin{theorem}\label{thm:conn-GabUlm}
For a locally small category $V$ the following statements are equivalent.
\begin{enumerate}
\item  $V$ is locally c-presentable.\label{lc1}
\item  $V$ is cocomplete and has a small dense subcategory consisting of small connected objects.\label{lc2}
\item  $V$ is a full subcategory of a presheaf category for which the inclusion is accessible, coproduct preserving and has a left adjoint.\label{lc4}
\item  $V$ is the category of models for a limit sketch whose distingished cones are connected.\label{lc5}
\item  $V$ is locally presentable and every object of $V$ is a coproduct of connected objects.\label{lc6}
\item  $V$ is locally presentable, extensive and the functor $(-){\cdot}1:\Set{\rightarrow}V$ has a left adjoint.\label{lc7}
\end{enumerate}
\end{theorem}
\begin{examples}\label{ex:lc-groth-toposes}
By theorem(\ref{thm:conn-GabUlm})(\ref{lc7}) a Grothendieck topos is locally connected in the usual sense iff its underlying category is locally c-presentable.
\end{examples}
\noindent Just as with locally presentable categories, locally c-presentable categories are closed under many basic categorical constructions. For instance from theorem(\ref{thm:conn-GabUlm})(\ref{lc6}), one sees immediately that the slices of a locally c-presentable category are locally c-presentable from the corresponding result for locally presentable categories. Another instance of this principle is the following result.
\begin{theorem}\label{thm:acc-monad}
If $V$ is locally c-presentable and $T$ is an accessible coproduct preserving monad on $V$, then $V^T$ is locally c-presentable.
\end{theorem}
\begin{proof}
First we recall that colimits in $V^T$ can be constructed explicitly using colimits in $V$ and the accessibility of $T$ (see for instance \cite{TTT} for a discussion of this). By definition we have a regular cardinal $\lambda$ such that $T$ preserves $\lambda$-filtered colimits and $V$ is locally $\lambda$-presentable. Defining $\Theta_0$ to be the full subcategory of $V$ consisting of the $\lambda$-presentable and connected objects, $(T,\Theta_0)$ is a monad with arities in the sense of \cite{Fam2fun}. One has a canonical isomorphism
\[ \xymatrix @C=3.5em {{V^T} \ar[r]^-{V^T(i,1)} \ar[d]_{U} \save \POS?="d" \restore & {\PSh {\Theta}_T} \ar[d]^{\res_j} \save \POS?="c" \restore \\ V \ar[r]_-{V(i_0,1)} & {\PSh {\Theta}_0} \POS "d";"c" **@{}; ?*{\iso}} \]
in the notation of \cite{Fam2fun}. Thus $V^T(i,1)$ is accessible since $\res_j$ creates colimits, and $T$ and $V(i_0,1)$ are accessible. By the nerve theorem of \cite{Fam2fun} $V^T(i,1)$ is also fully faithful, it has a left adjoint since $V^T$ is cocomplete given by left extending $i$ along the yoneda embedding, and so we have exhibited $V^T$ as conforming to theorem(\ref{thm:conn-GabUlm})(\ref{lc4}).
\end{proof}
\begin{examples}\label{ex:lc-noperad-algebras}
An $n$-operad for $0{\leq}n{\leq}\omega$ in the sense of \cite{Bat98}, gives a finitary coproduct preserving monad on the category $\PSh {\G}_{{\leq}n}$ of $n$-globular sets, and its algebras are just the algebras of the monad. Thus the category of algebras of any $n$-operad is locally c-presentable by theorem(\ref{thm:acc-monad}).
\end{examples}
%

\subsection{What $\ca G$ preserves}\label{ssec:Gamma-pres}
At the object level, to apply $\Gamma$ is to apply $\ca G$, so we shall now collect together many of the categorical properties that $\ca G$ preserves.
For $V$ with an initial object $\emptyset$, we saw in section(\ref{ssec:NMonad->MonMonad}) how to construct coproducts in $\ca GV$ explicitly. From this explicit construction, it is clear that the connected components of a $V$-graph $X$ may be described as follows. Objects $a$ and $b$ of $X$ are in the same connected component iff there exists a sequence $(x_0,...,x_n)$ of objects of $X$ such that for $1{\leq}i{\leq}n$ the hom $X(x_{i-1},x_i)$ is non-initial. Moreover $X$ is clearly the coproduct of its connected components, coproducts are disjoint and the initial object of $\ca GV$, whose $\Set$ of objects is empty, is strict. Thus by lemma(\ref{lem:easy-ext}) we obtain
\begin{proposition}\label{prop:GV-ext}
If $V$ has an initial object then $\ca GV$ is extensive and every object of $\ca GV$ is a coproduct of connected objects.
\end{proposition}
\noindent Given finite limits in $V$ it is straight forward to construct finite limits in $\ca GV$ directly. The terminal $V$-graph has one object and its only hom is the terminal object of $V$. Given maps $f:A{\rightarrow}B$ and $g:C{\rightarrow}B$ in $\ca GV$ their pullback $P$ can be constructed as follows. Objects are pairs $(a,c)$ where $a$ is an object of $A$ and $c$ is an object of $C$ such that $fa=gc$. The hom $P((a_1,c_1),(a_2,c_2))$ is obtained as the pullback of
\[ \xymatrix{{A(a_1,a_2)} \ar[r]^-{f} & {B(fa_1,fa_2)} & {C(c_1,c_2)} \ar[l]_-{g}} \]
in $V$. Thus one has
\begin{proposition}\label{prop:GV-lext}
If $V$ has finite limits then so does $\ca GV$. If in addition $V$ has an initial object, then $\ca GV$ is lextensive and every object of $V$ is a coproduct of connected objects.
\end{proposition}
\noindent As for cocompleteness one has the following result due to Betti, Carboni, Street and Walters.
\begin{proposition}\label{prop:GV-cocomp}\cite{BCSW-VarEnr}
If $V$ is cocomplete then so is $\ca GV$ and $(-)_0:\ca GV{\rightarrow}\Set$ is cocontinuous.
\end{proposition}
\noindent We now turn to local c-presentability. First we require a general lemma which produces a dense subcategory of $\ca GV$ from one in $V$ in a canonical way.
\begin{lemma}\label{lem:GV-dense}
Let $\ca D$ be a full subcategory of $V$ and suppose that $V$ has an initial object. Define an associated full subcategory $\ca D'$ of $\ca GV$ as follows:
\begin{itemize}
\item  $0 \in \ca D'$.
\item  $D \in \ca D \implies (D) \in \ca D'$.
\end{itemize}
If $\ca D$ is dense then so is $\ca D'$. For a regular cardinal $\lambda$, if the objects of $\ca D$ are $\lambda$-presentable then so are those of $\ca D'$.
\end{lemma}
\begin{proof}
Given functions
\[ f_{D'} : \ca GV(D',X) \rightarrow \ca GV(D',Y) \]
natural in $D' \in \ca D'$, we must show that there is a unique $f:X{\rightarrow}Y$ such that $f_{D'}=\ca GV(D',f)$. The object map of $f$ is forced to be $f_0$, and naturality with respect to the maps
\[ \xymatrix{0 \ar[r]^-{0} & (D) & 0 \ar[l]_-{1}} \]
ensures that the functions $f_{D'}$ amount to $f_0$ together with functions
\[ f_{D,a,b} : \ca GV_{\bullet}((0,(D),1),(a,X,b)) \rightarrow \ca GV_{\bullet}((0,(D),1),(f_0a,Y,f_0b)) \]
natural in $D \in \ca D$ for all $a,b \in X_0$. By the adjointness $L \ladj \varepsilon$ these maps are in turn in bijection with maps
\[ f'_{D,a,b} : V(D,X(a,b)) \rightarrow V(D,Y(f_0a,f_0b)) \]
natural in $D \in \ca D$ for all $a,b \in X_0$, and so by the density of $\ca D$ one has unique $f_{a,b}$ in $V$ such that$f'_{D,a,b}=V(D,f_{a,b})$. Thus $f_0$ and the $f_{a,b}$ together form the object and hom maps of the unique desired map $f$.

Since any colimit in $\ca GV$ is preserved by $(-)_0$ \cite{BCSW-VarEnr}, one can easily check directly that $0$ is $\lambda$-presentable for all $\lambda$. Let $D \in \ca D$ be $\lambda$-presentable. One has natural isomorphisms
\[ \begin{array}{c} {\ca GV((D),X) \iso \coprod\limits_{a,b{\in}X_0} \ca G_{\bullet}V((0,(D),1),(a,X,b)) \iso \coprod\limits_{a,b{\in}X_0} V(D,X(a,b))} \end{array} \]
exhibiting $\ca GV((D),-)$ as a coproduct of functors that preserve $\lambda$-filtered colimits, and thus is itself $\lambda$-filtered colimit preserving.
\end{proof}
\begin{corollary}\label{cor:GV-lc}
If $V$ is locally presentable then $\ca GV$ is locally c-presentable.
\end{corollary}
\begin{proof}
Immediate from theorem(\ref{thm:conn-GabUlm})(\ref{lc2}), lemma(\ref{lem:GV-dense}) and proposition(\ref{prop:GV-cocomp}).
\end{proof}
\noindent In \cite{KL-NiceVCat} Kelly and Lack proved that if $V$ is locally presentable then so is $\ca GV$ by an argument almost identical to that given here. The only difference is that in their version of lemma(\ref{lem:GV-dense}), their $\ca D'$ differs from ours only in that they use $0+0$ where we use $0$, and they instead prove that $\ca D'$ is a strong generator given that $\ca D$ is. We have given the above proof because the present form of lemma(\ref{lem:GV-dense}) is more useful to us in section(\ref{ssec:Gamma-pres12-lra}). Next we shall see that $\ca G$ preserves toposes. First a lemma of independent interest.
\begin{lemma}\label{lem:abstract-wreath}
Let $\C$ be a category and $\ca E$ be a lextensive category. Consider the category $\C_+$ constructed from $\C$ as follows. There is an injective on objects fully faithful functor
\[ \begin{array}{lcr} {i_{\C} : \C \rightarrow \C_+} && {C \mapsto C_+} \end{array} \]
and $\C_+$ has an additional object $0$ not in the image of $i_{\C}$. Moreover for each $C \in \C$ one has maps
\[ \begin{array}{lcr} {\sigma_C : 0 \rightarrow C_+} && {\tau_C : 0 \rightarrow C_+} \end{array} \]
and for all $f:C{\rightarrow}D$ one has the equations $f_+\sigma_C=\sigma_D$ and $f_+\tau_C=\tau_D$. Then $\ca G[\op \C,\ca E]$ is equivalent to the full subcategory of $[\C^{\textnormal{op}}_+,\ca E]$ consisting of those $X$ such that $X_0$ is a copower of $1$, the terminal object of $\ca E$.{\footnotemark{\footnotetext{Such $X_0$ in $\ca E$ are said to be \emph{discrete}.}}}
\end{lemma}
\begin{proof}
Let us write $\ca F$ for the full subcategory of $[\C^{\textnormal{op}}_+,\ca E]$ described in the statement of the lemma. We shall describe the functors
\[ \begin{array}{lcr} {(-)_+ : \ca G[\op {\C},\ca E] \rightarrow \ca F} && {(-)^{-} : \ca F \rightarrow \ca G[\op {\C},\ca E]} \end{array} \]
which provide the desired equivalence directly. Given $X \in \ca G[\op {\C},\ca E]$ define $X_+0={X_0}{\cdot}1$ and
\[ \begin{array}{c} {X_+C_+ = \coprod\limits_{a,b \in X_0} X(a,b)(C)} \end{array}. \]
In the obvious way this definition is functorial in $X$ and $C$. Conversely given $Y \in \ca F$ choose a set $Y^{-}_0$ such that $Y0={Y^{-}_0}{\cdot}1$. Such a set is determined uniquely up to isomorphism. Then define the homs of $Y^{-}$ via the pullbacks
\[ \PbSq {Y^{-}(a,b)(C)} {YC_+} {Y0{\times}Y0} {1} {} {} {(Y\sigma,Y\tau)} {(a,b)} \]
in $\ca E$ for all $a,b \in Y^{-}_0$. Note that since $\ca E$ is lextensive and hence distributive, $Y0{\times}Y0$ is itself the coproduct of copies of $1$ indexed by such pairs $(a,b)$. The natural isomorphisms
\[ \begin{array}{lcr} {X^{-}_+(a,b)(C) \iso X(a,b)(C)} && {Y^{-}_+C_+ \iso YC_+} \end{array} \]
come from extensivity.
\end{proof}
\begin{corollary}\label{cor:GV-topos}
\begin{enumerate}
\item If $V$ is a presheaf topos then so is $\ca GV$.\label{G6}
\item If $V$ is a Grothendieck topos then $\ca GV$ is a locally connected Grothendieck topos.\label{G7}
\end{enumerate}
\end{corollary}
\begin{proof}
Applying lemma(\ref{lem:abstract-wreath}) in the case where $\C$ is small and $\ca E=\Set$ one obtains the formula
\[ \ca G \PSh {\C} \catequiv \PSh {\C_+} \]
and thus (\ref{G6}). Since a Grothendieck topos is a left exact localisation of a presheaf category, the 2-functoriality of $\ca G$ together with (\ref{G6}), corollary(\ref{cor:GV-lc}) and example(\ref{ex:lc-groth-toposes}), implies that to establish (\ref{G7}) it suffices to show that $\ca G$ preserves left exact functors between categories with finite limits. This follows immediately from the explicit description of finite limits in $\ca GV$ given in the proof of proposition(\ref{prop:GV-lext}).
\end{proof}
\begin{remark}
The construction $(-)_+$ in the previous proof is easily discovered by thinking about why, as pointed out in section(\ref{ssec:enriched-graphs}), applying $\ca G$ successively to the empty category does one produce the categories of $n$-globular sets for $n \in \N$. The construction $(-)_+$ is just the general construction on categories, which when  applied successively to the empty category produces the categories $\G_{{\leq}n}$ for $n \in \N$, presheaves on which are by definition $n$-globular sets.
\end{remark}
%

\subsection{$\Gamma$'s one and two-cell map}\label{ssec:Gamma-pres12-cart}
First we note that by the explicit description of coproducts, the terminal object and pullbacks of enriched graphs, and the formula of corollary(\ref{cor:explicit-gamma}), one has
\begin{proposition}\label{prop:GammaE-basic}
Let $V$ and $W$ have coproducts and $E:MV{\rightarrow}W$ be distributive.
\begin{enumerate}
\item  $\Gamma{E}$ preserves coproducts.\label{GEb1}
\item  If $E$ preserves the terminal object then so does $\Gamma{E}$.
\item  If $E$ preserves pullbacks then so does $\Gamma{E}$.
\end{enumerate}
\end{proposition}
\noindent The precise conditions under which $\Gamma$ preserves and reflects cartesian natural transformations are identified by proposition(\ref{prop:Gamma-cart}). First we require a lemma which generalises lemma(7.4) of \cite{EnHopI}, whose proof follows easily from the explicit description of pullbacks in $\ca GV$ discussed in (\ref{ssec:Gamma-pres}).
\begin{lemma}\label{lem:pb-hom}
Suppose that $V$ has pullbacks. Given a commutative square (I)
\[ \TwoDiagRel {\xymatrix{W \ar[d]_{f} \save \POS?="d" \restore \ar[r]^-{h} & X \ar[d]^{g} \save \POS?="c" \restore \\ Y \ar[r]_-{k} & Z \POS "d";"c" **@{}; ?*{\textnormal{I}}}} {}
{\xymatrix{{W(a,b)} \ar[d]_{f_{a,b}} \save \POS?="d" \restore \ar[r]^-{h_{a,b}} & {X(ha,hb)} \ar[d]^{g_{ha,hb}} \save \POS?="c" \restore \\ {Y(a,b)} \ar[r]_-{k_{a,b}} & {Z(ha,hb)} \POS "d";"c" **@{}; ?*{\textnormal{II}}}} \]
in $\ca GV$ such that $f_0$ and $g_0$ are identities, one has for each $a,b \in W_0$ commuting squares (II) as in the previous display. The square (I) is a pullback iff for all $a,b \in W_0$, the square (II) is a pullback in $V$.
\end{lemma}
\begin{proposition}\label{prop:Gamma-cart}
Let $V$ have coproducts, $W$ be extensive, and $E$ and $F:MV{\rightarrow}W$ be distributive.
Then $\phi:E{\rightarrow}F$ is cartesian iff $\Gamma{\phi}$ is cartesian.
\end{proposition}
\begin{proof}
Let $f:X{\rightarrow}Y$ be in $\ca GV$, $a,b \in X_0$ and $x_0,...,x_n$ be a sequence of objects of $X$ such that $x_0=a$ and $x_n=b$. For each such $f,a,b$ it suffices by lemma(\ref{lem:pb-hom}), to show that the square (II) in the commutative diagram
\[ \xygraph{!{0;(2,0):(0,.5)::}
{\opE\limits_iX(x_{},x_i)}="i11" [r] {\Gamma{E}X(a,b)}="i12" [d] {\Gamma{E}Y(fa,fb)}="i22" [l] {\opE\limits_iY(fx_{},fx_i)}="i21" "i11" [ul] {\Gamma{F}X(a,b)}="o11"  "i12" [ur] {\opF\limits_iX(x_{},x_i)}="o12" "i21" [dl] {\opF\limits_iY(fx_{},fx_i)}="o21" "i22" [dr] {\Gamma{F}Y(fa,fb)}="o22"
"i11":"i12"^-{c_{x_i}}:"i22"^{\Gamma{E}(f)_{a,b}}|{}="l3" "i11":"i21"_{\opE\limits_if}|{}="l2":"i22"_-{c_{x_i}}
"o11":"o12"^-{c_{x_i}}:"o22"^{\Gamma{F}(f)_{a,b}}|{}="l4" "o11":"o21"_{\opF\limits_if}|{}="l1":"o22"_-{c_{x_i}}
"i11":"o11"^{\phi} "i21":"o21"^{\phi} "i12":"o12"^{\Gamma\phi} "i22":"o22"^{\Gamma\phi}
"l1":@{}"l2"|(.4)*{(III)}:@{}"l3"|*{(I)}:@{}"l4"|(.7)*{(II)}} \]
is a pullback. The square (I) and the largest square are pullbacks since $W$ is extensive and the $c$-maps are coproduct inclusions, and (III) is a pullback since $\phi$ is cartesian. Thus for all such $x_0,...,x_n$ the composite of (I) and (II) is a pullback. The result follows since the $c_{x_i}$ form coproduct cocones and $W$ is extensive. Conversely suppose that we have $f_i:X_i{\rightarrow}Y_i$ in $V$ for $1{\leq}i{\leq}n$. Then by the isomorphisms $E{\iso}\overline{\Gamma(E)}$ and $F{\iso}\overline{\Gamma(F)}$, and their naturality with respect to $\phi$ one may identify the naturality square on the left
\[ \TwoDiagRel {\xymatrix{{\opE\limits_iX_i} \ar[r]^-{\phi_{X_i}} \ar[d]_{\opE\limits_if_i} & {\opF\limits_i} \ar[d]^{\opF\limits_if_i} \\ {\opE\limits_iY_i} \ar[r]_-{\phi_{Y_i}} & {\opF\limits_iY_i}}} {} {\xymatrix{{\Gamma(E)(X_1,...,X_n)(0,n)} \ar[r]^-{\Gamma(\phi)_{0,n}} \ar[d]_{\Gamma(E)(f_1,...,f_n)_{0,n}} & {\Gamma(F)(X_1,...,X_n)(0,n)} \ar[d]_{\Gamma(F)(f_1,...,f_n)_{0,n}} \\ {\Gamma(E)(Y_1,...,Y_n)(0,n)} \ar[r]_-{\Gamma(\phi)_{0,n}} & {\Gamma(F)(Y_1,...,Y_n)(0,n)}}} \]
with that on the right in the previous display, which is cartesian by lemma(\ref{lem:pb-hom}) and since $\Gamma(\phi)$ is cartesian, and so $\phi$ is indeed cartesian.
\end{proof}
\noindent $\Gamma$'s
compatibility with the bicategory structure of $\Dist$ is expressed in
\begin{proposition}\label{prop:Gamma-cart2}
\begin{enumerate}
\item  Let $E:MV{\rightarrow}W$ and $F:MU{\rightarrow}V$ be distributive and $W$ be extensive. If $E$ and $F$ preserve pullbacks then so does $E \comp F$.\label{cart2-1}
\item  Let $E,E':MV{\rightarrow}W$ and $F,F':MU{\rightarrow}V$ be distributive and pullback preserving, and $W$ be extensive. If $\phi:E{\rightarrow}E'$ and $\psi:F{\rightarrow}F'$ are cartesian then so is $\phi \comp \psi$.\label{cart2-2}
\end{enumerate}
\end{proposition}
\begin{proof}
(\ref{cart2-1}). Suppose that for $1{\leq}i{\leq}n$ we have pullback squares
\[ \PbSq {A_i} {B_i} {C_i} {D_i} {h_i} {f_i} {k_i} {g_i} \]
in $U$. Then for each partition $n_1+...+n_k=n$ of $n$ we have a commutative diagram
\[ \xygraph{!{0;(2,0):(0,.5)::}
{(\opEoF\limits_i)A_i}="i11" [r] {(\opEoF\limits_i)B_i}="i12" [d] {(\opEoF\limits_i)D_i}="i22" [l] {(\opEoF\limits_i)C_i}="i21" "i11" [ul] {\opE\limits_i\opF\limits_jA_{ij}}="o11"  "i12" [ur] {\opE\limits_i\opF\limits_jB_{ij}}="o12" "i21" [dl] {\opE\limits_i\opF\limits_jC_{ij}}="o21" "i22" [dr] {\opE\limits_i\opF\limits_jD_{ij}}="o22"
"i11":"i12"^-{(\opEoF\limits_i)f_i}:"i22"^{(\opEoF\limits_i)k_i}|{}="l3" "i11":"i21"_{(\opEoF\limits_i)h_i}|{}="l2":"i22"_-{(\opEoF\limits_i)g_i}
"o11":"o12"^-{(\opE\limits_i\opF\limits_j)f_{ij}}:"o22"^{(\opE\limits_i\opF\limits_j)k_{ij}}|{}="l4" "o11":"o21"_{(\opE\limits_i\opF\limits_j)h_{ij}}|{}="l1":"o22"_-{(\opE\limits_i\opF\limits_j)g_{ij}}
"o11":"i11"^{c_{n_i}} "o21":"i21"^{c_{n_i}} "o12":"i12"^{c_{n_i}} "o22":"i22"^{c_{n_i}}
"l1":@{}"l2"|(.35)*{(I)}:@{}"l3"|*{(II)}:@{}"l4"|(.7)*{(III)}} \]
in $W$ in which the $c$-maps are coproduct inclusions, the coproducts being indexed over the set of all such partitions. For the outer square $1{\leq}i{\leq}k$ and $1{\leq}j{\leq}n_i$. We must show that (II) is a pullback. The squares (I) and (III) are pullbacks because $W$ is extensive. The large square is a pullback because $E$ and $F$ preserve pullbacks. Thus the composite of (I) and (II) is a pullback, and so the result follows by the extensivity of $W$.

(\ref{cart2-2}). Given $f_i:A_i{\rightarrow}B_i$ for $1{\leq}i{\leq}n$ and $n_1+...+n_k=n$ we have a commutative diagram
\[ \xygraph{!{0;(2,0):(0,.5)::}
{(\opEoF\limits_i)A_i}="i11" [r] {(\opEoF\limits_i)B_i}="i12" [d] {(\opEoFpr\limits_i)B_i}="i22" [l] {(\opEoFpr\limits_i)A_i}="i21" "i11" [ul] {\opE\limits_i\opF\limits_jA_{ij}}="o11"  "i12" [ur] {\opE\limits_i\opF\limits_jB_{ij}}="o12" "i21" [dl] {\opEpr\limits_i\opFpr\limits_jA_{ij}}="o21" "i22" [dr] {\opEpr\limits_i\opFpr\limits_jB_{ij}}="o22"
"i11":"i12"^-{(\opEoF\limits_i)f_i}:"i22"^{\psi{\comp}\phi} "i11":"i21"_{\psi{\comp}\phi}:"i22"_-{(\opEoFpr\limits_i)f_i}
"o11":"o12"^-{(\opE\limits_i\opF\limits_j)f_{ij}}:"o22"^{\psi\phi} "o11":"o21"_{\psi\phi}:"o22"_-{(\opEpr\limits_i\opFpr\limits_j)f_{ij}}
"o11":"i11"^{c_{n_i}} "o21":"i21"^{c_{n_i}} "o12":"i12"^{c_{n_i}} "o22":"i22"^{c_{n_i}}} \]
in $W$, to which we apply a similar argument as in (\ref{cart2-1}) to demonstrate that the inner square is cartesian.
\end{proof}
%

\subsection{The general basic correspondence between operads and multitensors}\label{ssec:general-op-mult}
Let $V$ be lextensive. Then by proposition(\ref{prop:Gamma-cart2}) the monoidal structure of $\Dist(V,V)$ restricts to pullback preserving $MV{\rightarrow}V$ and cartesian transformations between them. A multitensor $(E,u,\sigma)$ on $V$ is \emph{cartesian} when $E$ preserves pullbacks and $u$ and $\sigma$ are cartesian. By slicing over $E$ one obtains a monoidal category $\Coll E$, whose objects are cartesian transformations $\alpha:A{\rightarrow}E$. To give $(A,\alpha)$ a monoid structure is to give $A$ the structure of a cartesian multitensor such that $\alpha$ is a cartesian multitensor morphism. Such an $(A,\alpha)$ is called an \emph{$E$-multitensor}. The category of $E$-multitensors is denoted $\Mult E$.
\begin{example}\label{ex:non-sig-operads}
Let us denote by $\prod$ the multitensor on $\Set$ given by finite products. By the lextensivity of $\Set$, $\prod$ is an lra multitensor and thus cartesian. A $\prod$-multitensor is the same thing as a non-symmetric operad in $\Set$. For given a cartesian multitensor map $\varepsilon:E{\rightarrow}\prod$, one obtains the underlying sequence $(E_n \, : \, n \in \N)$ of sets of the corresponding operad as $E_n = \opE\limits_{1{\leq}i{\leq}n} 1$. One uses the cartesianness of the naturality squares corresponding to the maps $(X_1,...,X_n){\rightarrow}(1,...,1)$ to recover $E$ from the $E_n$. Similarly the multitensor structure of $E$ corresponds to the unit and substitution maps making the $E_n$ into an operad.
\end{example}
Given a cartesian normalised monad $T$ on $\ca GV$, one obtains a monoidal category $\NColl T$, whose objects are cartesian transformations $\alpha:A{\rightarrow}T$ over $\Set$. Explicitly, to say that a general collection, which is a cartesian transformation $\alpha:A{\rightarrow}T$, is over $\Set$, is to say that the components of $\alpha$ are identities on objects maps of $V$-graphs. The tensor product of $\NColl T$ is obtained via composition and the monad structure of $T$, and a monoid structure on $(A,\alpha)$ is a cartesian monad structure on $A$ such that $\alpha$ is a cartesian monad morphism. Such an $(A,\alpha)$ is called a \emph{$T$-operad over $\Set$}{\footnotemark{\footnotetext{In \cite{Bat98} and \cite{EnHopI} these were called \emph{normalised} $T$-operads.}}}. The category of $T$-operads over $\Set$ is denoted as $\NOp T$.
\noindent By the results of this subsection $\Gamma$ induces a strong monoidal functor
\[ \Gamma_E : \Coll E \rightarrow \NColl {\Gamma{E}} \]
and we shall now see that this functor is a monoidal equivalence. Applying this equivalence to the monoids in the respective monoidal categories gives the promised general equivalence between multitensors and operads over $\Set$ in corollary(\ref{cor:mult-nop-equiv}) below.
\begin{lemma}\label{lem:transfer-dpl}
Let $V$ be a lextensive category and $T$ be a cartesian monad on $\ca GV$ over $\Set$. Let $\alpha:A{\rightarrow}T$ be a collection over $\Set$. 
\begin{enumerate}
\item  If $T$ is distributive then so is $A$.\label{tdpl1}
\item  If $T$ is path-like then so is $A$.\label{tdpl2}
\end{enumerate}
\end{lemma}
\begin{proof}
(\ref{tdpl1}): given an $n$-tuple $(X_1,...,X_n)$ of objects of $V$ and a coproduct cocone
\[ (c_j : X_{ij} \rightarrow X_i \,\, : \,\, j \in J) \]
where $1{\leq}i{\leq}n$, we must show that the hom-maps
\[ A(X_1,...,c_j,...,X_n)_{0,n} : A(X_1,...,X_{ij},...,X_n)(0,n) \rightarrow A(X_1,...,X_i,...,X_n)(0,n) \]
form a coproduct cocone. For $j \in J$ we have a pullback square
\[ \xygraph{!{0;(6,0):(0,.167)::}
{A(X_1,...,X_{ij},...,X_n)(0,n)}="tl" [r] {A(X_1,...,X_i,...,X_n)(0,n)}="tr" [d] {T(X_1,...,X_i,...,X_n)(0,n)}="br" [l] {T(X_1,...,X_{ij},...,X_n)(0,n)}="bl" "tl"(:"tr"^-{A(X_1,...,c_j,...,X_n)_{0,n}}:"br"^{\alpha},:"bl"_{\alpha}:"br"_-{T(X_1,...,c_j,...,X_n)_{0,n}})
"tl" [r(.1)d(.4)] {\xybox{\xygraph{!{0;(0.2,0):} (:@{-}[u], :@{-}[l])}}}} \]
and by the distributivity of $T$ and lemma(\ref{lem:pb-hom}), the $T(X_1,...,c_j,...,X_n)_{0,n}$ form a coproduct cocone, and thus so do the $A(X_1,...,c_j,...,X_n)_{0,n}$ by the extensivity of $V$.\\
(\ref{tdpl2}): given $X \in \ca GV$, $a,b \in X_0$ and a sequence $(x_0,...,x_n)$ of objects of $X$ such that $x_0=a$ and $x_n=b$, we have the map
\[ A\overline{x}_{0,n} : Ax^*X(0,n) \rightarrow AX(a,b) \]
and we must show that these maps, where the $x_i$ range over all sequences from $a$ to $b$, form a coproduct cocone. By the path-likeness of $T$ we know that the maps
\[ T\overline{x}_{0,n} : Tx^*X(0,n) \rightarrow TX(a,b) \]
form a coproduct cocone, so we can use the cartesianness of $\alpha$, lemma(\ref{lem:pb-hom}) and the extensivity of $V$ to conclude as in (\ref{tdpl1}).
\end{proof}
\begin{proposition}\label{prop:coll-equiv}
Let $V$ be lextensive and $(E,\iota,\sigma)$ be a distributive cartesian multitensor on $V$. Then $\Gamma_E$ is a monoidal equivalence $\Coll E \catequiv \NColl {\Gamma{E}}$.
\end{proposition}
\begin{proof}
The functor $\Gamma_E$ is the result of applying the functor $\Gamma_{V,V}$ of proposition(\ref{prop:pl-adjoint-char}) over $E$. Thus by proposition(\ref{prop:pl-adjoint-char}), proposition(\ref{prop:Gamma-cart}) and lemma(\ref{lem:transfer-dpl}), $\Gamma_E$ is an equivalence.
\end{proof}
\begin{corollary}\label{cor:mult-nop-equiv}
Let $V$ be lextensive and $(E,\iota,\sigma)$ be a distributive cartesian multitensor on $V$. Then applying $\Gamma_E$ gives $\Mult E \catequiv \NOp {\Gamma E}$.
\end{corollary}
%

\subsection{$\Gamma$ and local right adjoint monads}\label{ssec:Gamma-pres12-lra}
Local right adjoint monads, especially defined on presheaf categories, are fundamental to higher category theory. Indeed a deeper understanding of such monads is the key to understanding the relationship between the operadic and homotopical approaches to the subject \cite{Fam2fun}. We will now understand the conditions under which $\Gamma$ preserves local right adjoints. First we require two lemmas.
\begin{lemma}\label{lem:partial-adjoint}
Let $R:V{\rightarrow}W$ be a functor, $V$ be cocomplete, $U$ be a small dense full subcategory of $W$, and $L:U{\rightarrow}V$ be a partial left adjoint to $R$, that is to say, one has isomorphisms $W(S,RX) \iso V(LS,X)$ natural in $S \in U$ and $X \in V$. Defining $\overline{L}:W{\rightarrow}V$ as the left kan extension of $L$ along the inclusion $I:U{\rightarrow}W$, one has $\overline{L} \ladj R$.
\end{lemma}
\begin{proof}
Denoting by $p:I/Y{\rightarrow}U$ the canonical forgetful functor for $Y \in W$ and recalling that $\overline{L}Y = \colim(Lp)$, one obtains the desired natural isomorphism as follows
\[ \begin{array}{rcccl} {V(\overline{L}Y,X)} & {\iso} & {[I/Y,V](Lp,\textnormal{const}(X))} & {\iso} & {\lim_{f{\in}I/Y} V(L(\textnormal{dom}(f)),X)} \\ & {\iso} & {\lim_f W(\textnormal{dom}(f),RX)} & {\iso} &  {\ca B(Y,RX)} \end{array} \]
for all $X \in V$.
\end{proof}
\begin{lemma}\label{lem:lra-dense}
Let $T:V{\rightarrow}W$ be a functor, $V$ be cocomplete and $W$ have a small dense subcategory $U$. Then $T$ is a local right adjoint iff every $f:S{\rightarrow}TX$ with $A \in U$ admits a generic factorisation. If in addition $V$ has a terminal object denoted $1$, then generic factorisations in the case $X=1$ suffice.
\end{lemma}
\begin{proof}
For the first statement ($\implies$) is true by definition so it suffices to prove the converse. The given generic factorisations provide a partial left adjoint $L:I/TX{\rightarrow}V$ to $T_X:V/X{\rightarrow}W/TX$ where $I$ is the inclusion of $U$. Now $I/TX$ is a small dense subcategory of $W/TX$, and so by the previous lemma $L$ extends to a genuine left adjoint to $T_X$. In the case where $V$ has $1$ one requires only generic factorisations in the case $X=1$ by the results of \cite{Fam2fun} section(2).
\end{proof}
\noindent The analogous result for presheaf categories, with the representables forming the chosen small dense subcategory, was discussed in \cite{Fam2fun} section(2).
\begin{proposition}\label{prop:GammaE-lra}
Let $V$ and $W$ be locally c-presentable and $E:MV{\rightarrow}W$ be distributive. If $E:MV{\rightarrow}W$ is a local right adjoint then so is $\Gamma{E}$.
\end{proposition}
\begin{proof}
Let $\ca D$ be a small dense subcategory of $W$ consisting of small connected objects. By lemma(\ref{lem:lra-dense}) and lemma(\ref{lem:GV-dense}) it suffices to exhibit generic factorisations of maps
\[ f:S \rightarrow \Gamma{E}1 \]
where $S$ is either $0$ or $(D)$ for some $D \in \ca D$. In the case where $S$ is $0$ the first arrow in the composite
\[ \xymatrix{0 \ar[r] & {\Gamma{E}0} \ar[r]^-{\Gamma{E}t} & {\Gamma{E}1}} \]
is generic because $0$ is the initial $W$-graph with one object (and $t$ here is the unique map). In the case where $S=(D)$, to give $f$ is to give a map $f':D{\rightarrow}E_n1$ in $V$, by corollary(\ref{cor:explicit-gamma}) since $D$ is connected. Since $E$ is a local right adjoint, $E_n$ is too and so one can generically factor $f'$ to obtain
\[ \xymatrix{D \ar[r]^-{g'_f} & {\opE\limits_iZ_i} \ar[r]^-{\opE\limits_it} & {E_n1}} \]
from which we obtain the generic factorisation
\[ \xymatrix{{(D)} \ar[r]^-{g_f} & {\Gamma{E}Z} \ar[r]^-{\Gamma{E}t} & {\Gamma{E}1}} \]
where $Z=(Z_1,...,Z_n)$, the object map of $g_f$ is given by $0 \mapsto 0$ and $1 \mapsto n$, and the hom map of $g_f$ is $g'_f$ composed with the coproduct inclusion.
\end{proof}
%

\subsection{$\Gamma$ and accessible functors}\label{ssec:Gamma-pres12-access}
First note that while it is a very different thing for $E:MV{\rightarrow}W$ to preserve coproducts compared with preserving coproducts in each variable, the situation is simpler for $\lambda$-filtered colimits, where $\lambda$ is any regular cardinal. Note that
\[ F : V_1 \times ... \times V_n \rightarrow W \]
preserves $\lambda$-filtered colimits in each variable iff $F$ preserves $\lambda$-filtered colimits. For given a connected category $C$, the colimit of a functor $C{\rightarrow}A$ constant at say $X$ is of course $X$, and since $\lambda$-filtered colimits are connected, one can prove $(\impliedby)$ by keeping all but the variable of interest constant. For the converse it is sufficient to prove that $F$ preserves colimits of chains of length less than $\lambda$, and this follows by a straight forward transfinite induction. Since $MV$ is a sum of $V^n$'s,
from the connectedness of $\lambda$-filtered colimits it is clear that $E$ preserves $\lambda$-filtered colimits iff each $E_n:V^n{\rightarrow}W$ does, and so we have proved
\begin{lemma}\label{lem:multi-lambda}
For $E:MV{\rightarrow}W$ the following statements are equivalent for any regular cardinal $\lambda$.
\begin{enumerate}
\item  $E$ preserves $\lambda$-filtered colimits in each variable.
\item  $E$ preserves $\lambda$-filtered colimits.
\item  $E_n:V^n{\rightarrow}W$ preserves $\lambda$-filtered colimits for all $n \in \N$.
\end{enumerate}
\end{lemma}
\noindent As already mentioned, colimits in $\ca GV$ for a cocomplete $V$ were calculated in \cite{BCSW-VarEnr}. Let us spell out transfinite composition in $\ca GV$. Given an ordinal $\lambda$ and a $\lambda$-chain
\begin{equation}\label{chain1} \begin{array}{lcr} {i \leq j \in \lambda} & {\mapsto} & {f_{ij} : X_i \rightarrow X_j} \end{array} \end{equation}
in $\ca GV$ with colimit $X$, one may consider the induced $\lambda$-chain
\begin{equation}\label{chain2} \begin{array}{lcr} {i \leq j \in \lambda} & {\mapsto} & {(f_{ij})_0 \times (f_{ij})_0 : (X_i)_0 \times (X_i)_0 \rightarrow (X_j)_0 \times (X_j)_0} \end{array} \end{equation}
in $\Set$. This will have colimit $X_0 \times X_0$ because $\lambda$-filtered colimits and products commute in $\Set$ and $()_0:V{\rightarrow}\Set$ is cocontinuous. For $(a,b) \in X_0$ let us denote by $D_{a,b}$ the full subcategory of the category of elements of (\ref{chain2}), consisting of those elements which are sent to $(a,b)$ by the universal cocone. We shall call this the \emph{$(a,b)$-component} of the chain (\ref{chain1}). Now $D_{a,b}$ is of course no longer a chain, but one may easily verify that it is $\lambda$-filtered. By the explicit description of colimits in $\Set$, the $D_{a,b}$ are just the connected components of the category of elements of (\ref{chain2}). To pairs $(a',b') \in (X_i)_0 \times (X_i)_0$ which are elements of $D_{a,b}$, one may associate the corresponding hom $X_i(a',b') \in V$, and in this way build a functor $F_{a,b}:D_{a,b}{\rightarrow}V$. The hom $X(a,b)$ is the colimit of this functor.
\begin{proposition}\label{prop:Gamma-accessible}
Let $V$ and $W$ be cocomplete, $E:MV{\rightarrow}W$ be distributive and $\lambda$ be a regular cardinal. If $E$ preserves $\lambda$-filtered colimits in each variable then $\Gamma{E}$ preserves $\lambda$-filtered colimits.
\end{proposition} 
\begin{proof}
It suffices to show $\Gamma{E}$ preserves colimits of $\lambda$-chains. Consider the chain (\ref{chain1}) in $\ca GV$. For all $n \in \N$ one has an $(n+1)$-ary version of (\ref{chain2}), that is involving $(n+1)$-fold instead of binary cartesian products in $\Set$. These of course also commute with $\lambda$-filtered colimits. Similarly one obtains a $\lambda$-filtered category $D_{x_0,...,x_n}$ and a functor
\[ \begin{array}{lcr} {F_{x_0,...,x_n} : D_{x_0,...,x_n} \rightarrow V^n} && {(y_0,...,y_n) \in X_i \mapsto (X_i(y_0,y_1),...,X_i(y_{n-1},y_n))} \end{array} \]
Applying $\Gamma{E}$ does nothing at the object level. Let us write $D'_{a,b}$ for the $(a,b)$-component of the chain obtained by applying $\Gamma{E}$ to (\ref{chain1}), and $F'_{a,b}$ for the corresponding functor into $W$. From the explicit description of $\Gamma$'s effect on homs of corollary(\ref{cor:explicit-gamma}), one sees that $F'_{a,b}$ is the coproduct of the composites
\[ \xymatrix{{D_{x_0,...,x_n}} \ar[r]^-{F_{x_0,...,x_n}} & {V^n} \ar[r]^-{E_n} & W}. \]
over all sequences $(x_0,...,x_n)$ starting at $a$ and finishing at $b$. By lemma(\ref{lem:multi-lambda}) the colimits of the $F_{x_0,...,x_n}$ are preserved by the $E_n$, and so by the explicit description of colimits of $\lambda$-chains in $\ca GW$, the colimit of (\ref{chain1}) is indeed preserved by $\Gamma{E}$.
\end{proof}
%

\subsection{An elegant construction of the strict $n$-category monads}\label{ssec:induction}
Let us recall the construction $(-)^\times$. Given a monad $(T,\eta,\mu)$ on $V$ a category with products, one has a multitensor $T^{\times}$ defined by
\[ \begin{array}{c} {T^{\times}(X_1,...,X_n) = \prod\limits_{1{\leq}i{\leq}n} T(X_i)} \end{array} \]
and the unit and substitution is induced in the obvious way from $\eta$ and $\mu$. When $V$ is lextensive, $T$ is a local right adjoint, and $\eta$ and $\mu$ are cartesian, it follows that $T^{\times}$ is a local right adjoint its unit and multiplication are also cartesian. When $T$ preserves coproducts and the cartesian product for $V$ is distributive, then $T^{\times}$ is a distributive multitensor. If in addition finite limits and filtered colimits commute in $V$ (which happens when, for example $V$ is locally finitely presentable), then $T^{\times}$ is finitary. Moreover by proposition(2.8) of \cite{EnHopI} one has
\begin{equation}\label{eq:Tcross-cat} \Enrich {T^{\times}} \iso \Enrich {V^T} \end{equation}
where the enrichment on the right hand side is with respect to cartesian products. Thus one can consider the following inductively-defined sequence of monads
\begin{itemize}
\item  Put $\ca T_{{\leq}0}$ equal to the identity monad on $\Set$.
\item  Given a monad $\ca T_{{\leq}n}$ on $\ca G^n\Set$, define the monad
$\ca T_{{\leq}n+1} = \Gamma \ca T^{\times}_{{\leq}n}$
on $\ca G^{n+1}\Set$.
\end{itemize}
recalling that $\ca G^n\Set$ is the category of $n$-globular sets.
\begin{theorem}\cite{EnHopI}\label{thm:EnHopI-main-theorem}
For $n \in \N$, $\ca T_{\leq n}$ is the strict $n$-category monad on $n$-globular sets. This monad is coproduct preserving, finitary and local right adjoint.
\end{theorem}
\begin{proof}
By (\ref{eq:Tcross-cat}) and corollary(\ref{cor:gamma-alg-ecat}), one has
\[ \ca G^{n+1}(\Set)^{\ca T_{{\leq}n+1}} \iso \Enrich {\ca T^{\times}_{{\leq}n}} \]
and so by definition $\ca G^n(\Set)^{\ca T_{\leq n}}$ is the category of strict $n$-categories and strict $n$-functors between them. By the remarks at the beginning of this section and corollary(\ref{cor:GV-topos}) $(-)^{\times}$ will produce a distributive, finitary, local right adjoint multitensor on a presheaf category when it is fed a coproduct preserving, finitary, local right adjoint monad on a presheaf category. By corollary(\ref{cor:GV-topos}), proposition(\ref{prop:GammaE-basic}), proposition(\ref{prop:Gamma-accessible}), proposition(\ref{prop:GammaE-lra}) and proposition(\ref{prop:Gamma-cart2}), $\Gamma$ will produce a coproduct preserving, finitary, local right adjoint monad on a presheaf category when it is fed a distributive, finitary, local right adjoint multitensor on a presheaf category. Thus the monads $\ca T_{\leq n}$ are indeed coproduct preserving, finitary and local right adjoint for all $n \in \N$.
\end{proof}
The objects of $\NOp {\ca T_{{\leq}n}}$ -- $n$-operads over $\Set$ -- were in \cite{Bat98} \cite{EnHopI} called ``normalised'' $n$-operads. Many $n$-categorical structures of interest, such as weak $n$-categories, can be defined as algebras of $n$-operads over $\Set$. Objects of $\Mult {\ca T_{{\leq}n}}$ are called $n$-multitensors. These are a nice class of lax monoidal structures on the category of $n$-globular sets. By corollary(\ref{cor:mult-nop-equiv}) and theorem(\ref{thm:EnHopI-main-theorem}) one obtains
\begin{corollary}\cite{EnHopI}\label{cor:nnp1op-nmult}
For all $n \in \N$, applying $\Gamma$ gives $\NOp {\ca T_{{\leq}(n+1)}} \catequiv \Mult {\ca T_{{\leq}n}}$.
\end{corollary}
\noindent That is to say, $\Gamma$ exhibits $(n+1)$-operads over $\Set$ and $n$-multitensors as the same thing, and under this correspondence, the algebras of the operad correspond to the categories enriched in the associated multitensor by corollary(\ref{cor:gamma-alg-ecat}).

\section{The 2-functoriality of the monad-multitensor correspondence}\label{sec:2-functoriality}

\subsection{Motivation}\label{ssec:motivation-2-functoriality}
Up to this point $\Gamma$ has been our notation for the process
\[ \textnormal{Distributive multitensor on $V$} \mapsto \textnormal{Monad on $\ca GV$ over $\Set$} \]
and $\overline{(-)}$ has been our notation for the reverse construction. For the most complete analysis of these constructions one must acknowledge that they are the object maps of 2-functors in two important ways. This 2-functoriality together with the formal theory of monads \cite{Str72} gives a satisfying explanation of how it is that monad distributive laws arise naturally in this subject (see \cite{ChengDist}).

\subsection{2-categories of multitensors and monads}\label{ssec:Dist-Mult}
As the lax-algebras of a 2-monad $M$ (see section(\ref{ssec:LMC})), lax monoidal categories form a 2-category $\LaxAlg M$. See \cite{LackCodesc} for a complete description of the 2-category of lax algebras for an arbitrary 2-monad. Explicitly a lax monoidal functor between lax monoidal categories $(V,E)$ and $(W,F)$ consists of a functor $H:V{\rightarrow}W$, and maps
\[ \psi_{X_i} : \opF\limits_i HX_i \rightarrow H \opE\limits_i X_i \]
natural in the $X_i$ such that
\[ \xygraph{{\xybox{\xygraph{!{0;(.75,0):(0,1.333)::} {HX}="l" [r(2)] {F_1HX}="r" [dl] {HE_1X}="b" "l"(:"r"^-{u_{HX}}:"b"^{\psi_X},:"b"_{Hu_X})}}} [r(5)]
{\xybox{\xygraph{!{0;(2,0):(0,.5)::} {\opF\limits_i\opF\limits_jHX_{ij}}="tl" [r] {\opF\limits_iH\opE\limits_jX_{ij}}="tm" [r] {H\opE\limits_i\opE\limits_jX_{ij}}="tr" [l(.5)d] {H\opE\limits_{ij}X_{ij}}="br" [l] {\opF\limits_{ij}HX_{ij}}="bl" "tl" (:@<1ex>"tm"^-{\opF\limits_i\psi}:@<1ex>"tr"^-{\psi\opE\limits_j}:"br"^{H\sigma},:"bl"_{\sigma{H}}:@<1ex>"br"_-{\psi})}}}} \]
commute for all $X$ and $X_{ij}$ in $V$. A monoidal natural transformation between lax monoidal functors \[ (H,\psi),(K,\kappa):(V,E){\rightarrow}(W,F) \] consists of a natural transformation $\phi:H{\rightarrow}K$ such that
\[ \xygraph{!{0;(2,0):(0,.5)::} {\opF\limits_iHX_i}="tl" [r] {H\opE\limits_iX_i}="tr" [d] {K\opE\limits_iX_i}="br" [l] {\opF\limits_iKX_i}="bl" "tl" (:@<1ex>"tr"^-{\psi}:"br"^{\phi\opE\limits_i},:"bl"_{\opF\limits_i\phi}:@<1ex>"br"_-{\kappa})} \]
commutes for all $X_i$.
\begin{definition}\label{def:2cat-DistMult}
The 2-category $\DISTMULT$ of distributive multitensors, is defined to be the full sub-2-category of $\LaxAlg M$ consisting of the $(V,E)$ such that $V$ has coproducts and $E$ is distributive.
\end{definition}
For any 2-category $\ca K$ recall the 2-category $\MND(\ca K)$ from \cite{Str72} of monads in $\ca K$. Another way to describe this very canonical object is that it is the 2-category of lax algebras of the identity monad on $\ca K$. Explicitly the 2-category $\MND(\CAT)$ has as objects pairs $(V,T)$ where $V$ is a category and $T$ is a monad on $V$. An arrow $(V,T){\rightarrow}(W,S)$ is a pair consisting of a functor $H:V{\rightarrow}W$ and a natural transformation $\psi:SH{\rightarrow}HT$ satisfying the obvious 2 axioms: these are just the ``unary'' analogues of the axioms for a lax monoidal functor written out above.
For example, any lax monoidal functor $(H,\psi)$ as above determines a monad functor $(H,\psi_1):(V,E_1){\rightarrow}(W,F_1)$. A monad transformation between monad functors
\[ (H,\psi),(K,\kappa):(V,T){\rightarrow}(W,S) \]
consists of a natural transformation $\phi:H{\rightarrow}K$ satisfying the obvious axiom. For example a monoidal natural transformation $\phi$ as above is a monad transformation $(H,\psi_1){\rightarrow}(K,\kappa_1)$.

In fact as we are interested in monads over $\Set$, we shall work not with $\MND(\CAT)$ but rather with $\MND(\CAT/\Set)$. An object $(V,T)$ of this latter 2-category is a category $V$ equipped with a functor into $\Set$, together with a monad $T$ on $V$ which ``acts fibre-wise''
with respect to this functor. That is $T$'s
object map doesn't affect the underlying object set, similarly for the arrow map of $T$, and the components of $T$'s
unit and multiplication are identities on objects in the obvious sense. An arrow $(V,T){\rightarrow}(W,S)$ of $\MND(\CAT/\Set)$ is a pair $(H,\psi)$ as in the case of $\MND(\CAT)$, with the added condition that $\psi$'s components are the identities on objects, and similarly the 2-cells of $\MND(\CAT/\Set)$ come with an extra identity-on-object condition.

\subsection{$\Gamma$ as a 2-functor}\label{ssec:Gamma-2-functor}
We shall now exhibit the 2-functor
\[ \Gamma : \DISTMULT \rightarrow \MND(\CAT/\Set) \]
which on objects is given by $(V,E) \mapsto (\ca GV,\Gamma{E})$. Let $(H,\psi):(V,E){\rightarrow}(W,F)$ be a lax monoidal functor between distributive lax monoidal categories. Then for $X \in \ca GV$ and $a,b \in X_0$, we define the hom map $\Gamma(\psi)_{X,a,b}$ to be the composite of
\[ \xygraph{!{0;(5,0):} {\coprod\limits_{a=x_0,...,x_n=b} \opF\limits_iHX(x_{i-1},x_i)}="l" [r] {\coprod\limits_{a=x_0,...,x_n=b} H\opE\limits_iX(x_{i-1},x_i)}="m" "l":@<1ex>"m"^-{\coprod \psi}} \]
and $H$'s coproduct preservation obstruction map.
It follows easily from the definitions that $(\ca GH,\Gamma(\psi))$ as defined here satisfies the axioms of a monad functor. Moreover given a monoidal natural transformation $\phi:(H,\psi){\rightarrow}(K,\kappa)$, it also follows easily from the definitions that \[ \ca G\phi:(\ca GH,\Gamma(\psi)){\rightarrow}(\ca GK,\Gamma(\kappa)) \] is a monad transformation. It is also straight-forward to verify that these assignments are 2-functorial.

\subsection{The image of $\Gamma$}\label{ssec:Gamma-Image}
By proposition(\ref{prop:pl-adjoint-char}) and corollary(\ref{cor:pl->copr-pres}) we understand objects of the image of $\Gamma$ and we collect this information in
\begin{proposition}\label{prop:image-Gamma-objects}
For $V$ a category with coproducts, a monad $(\ca GV, T)$ over $\Set$ is in the image of $\Gamma$ iff $T$ is distributive and path-like. Moreover any such $T$ automatically preserves coproducts. One recovers the distributive multitensor $E$ such that $\Gamma{E}{\iso}T$ as $E=\overline{T}$.
\end{proposition}
\noindent Since the construction $\overline{(-)}$ is itself obviously 2-functorial, the arrows and 2-cells in the image of $\Gamma$ may also be easily characterised.
\begin{proposition}\label{prop:image-Gamma-1-2-cells}
\begin{enumerate}
\item  Let $(V,E)$ and $(W,F)$ be distributive lax monoidal categories. A monad functor of the form
\[ (H,\psi) : (\ca GV,\Gamma{E}) \rightarrow (\ca GW,\Gamma{F}) \]
is in the image of $\Gamma$ iff $H=\ca GH'$ for some $H'$.\label{Gamma-char1}
\item  Let $(H,\psi),(K,\kappa):(V,E){\rightarrow}(W,F)$ be lax monoidal functors between distributive lax monoidal categories. A monad transformation
\[ \phi : (\ca GH,\Gamma{\psi}) \rightarrow (\ca GK,\Gamma{\kappa}) \]
is in the image of $\Gamma$ iff it is of the form $\phi=\ca G\phi'$.\label{Gamma-char2}
\end{enumerate}
\end{proposition}
\begin{proof}
By definition monad functors and transformations in the image of $\Gamma$ have the stated properties, so we must prove the converse.
Given $(\ca GH,\psi)$ such that the components of $\psi$ are the identities on objects, one recovers for $X_1,...,X_n$ from $V$, the corresponding lax monoidal functor coherence map as the hom map from $0$ to $n$ of the component $\psi_{(X_1,...,X_n)}$. That is to say, we apply $\overline{(-)}$ to the appropriate monad functors to prove (\ref{Gamma-char1}), and we do the same to the appropriate monad transformations to obtain (\ref{Gamma-char2}).
\end{proof}
\begin{definition}\label{def:PLMND}
We denote by $\PLMND$ the following 2-category. Its objects are monads $(\ca GV,T)$ over $\Set$ such that $V$ has coproducts and $T$ is distributive and path-like. Its arrows are arrows $(\ca GH,\psi) : (\ca GV,\Gamma{E}) \rightarrow (\ca GW,\Gamma{F})$ of $\MND(\CAT/\Set)$, and its 2-cells are 2-cells $\ca G\phi : (\ca GH,\psi) \rightarrow (\ca GK,\kappa)$ of $\MND(\CAT/\Set)$.
\end{definition}
\noindent Thus from the proof of proposition(\ref{prop:image-Gamma-1-2-cells}) we have
\begin{corollary}\label{cor:2eq-mult-mon}
$\Gamma$ and $\overline{(-)}$ provide a 2-equivalence $\DISTMULT \catequiv \PLMND$.
\end{corollary}
%

\subsection{The dual 2-functoriality of $\Gamma$}\label{ssec:Gamma-another-2-functor}
Lax algebras of a 2-monad organise naturally into \emph{two} different 2-categories depending on whether one takes lax or oplax algebra morphisms. So in particular one has the 2-category $\OpLaxAlg M$ of lax monoidal categories, \emph{op}lax-monoidal functors between them and monoidal natural transformations between those. The coherence $\psi$ for an oplax $(H,\psi):(V,E){\rightarrow}(W,F)$ goes in the other direction, and so its components look like this:
\[ \psi_{X_i} : H \opE\limits_i X_i \rightarrow \opF\limits_i HX_i. \]
The reader should easily be able to write down explicitly the two coherence axioms that this data must satisfy, as well as the condition that must be satisfied by a monoidal natural transformation between oplax monoidal functors. Similarly there is a dual version $\OpMND(\ca K)$ of the 2-category $\MND(\ca K)$ of monads in a given 2-category $\ca K$ discussed above \cite{Str72}. An arrow $(V,T){\rightarrow}(W,S)$ of $\OpMND(\CAT)$ consists of a functor $H:V{\rightarrow}W$ and a natural transformation $\psi:HT{\rightarrow}SH$ satisfying the two obvious axioms. An arrow of $\OpMND(\CAT)$ is called a monad opfunctor. As before $\OpMND(\CAT/\Set)$ differs from $\MND(\CAT/\Set)$ in that all the categories involved come with a functor into $\Set$, and all the functors and natural transformations involved are compatible with these forgetful functors.

We now describe the dual version of the 2-functoriality of $\Gamma$ discussed in sections(\ref{ssec:Gamma-2-functor}) and (\ref{ssec:Gamma-Image}). When defining the one-cell map of $\Gamma$ in section(\ref{ssec:Gamma-2-functor}) we were helped by the fact that the coproduct preservation obstruction went the right way: see the definition of the monad functor $(\ca GH,\Gamma{\psi})$ above. This time however we will not be so lucky, and for this reason we must restrict ourselves in the following definition to coproduct preserving oplax monoidal functors.
\begin{definition}\label{def:OpDistMult}
The 2-category $\OpDISTMULT$ is defined to be the locally full sub-2-category of $\OpLaxAlg M$ consisting of the distributive lax monoidal categories, and the oplax monoidal functors $(H,\psi)$ such that $H$ preserves coproducts. We denote by $\OpPLMND$ the following 2-category. Its objects are monads $(\ca GV,T)$ over $\Set$ such that $V$ has coproducts and $T$ is distributive and path-like. Its arrows are arrows $(\ca GH,\psi) : (\ca GV,\Gamma{E}) \rightarrow (\ca GW,\Gamma{F})$ of $\OpMND(\CAT/\Set)$ such that $H$ preserves coproducts. Its 2-cells are 2-cells $\ca G\phi : (\ca GH,\psi) \rightarrow (\ca GK,\kappa)$ of $\OpMND(\CAT/\Set)$.
\end{definition}
\noindent We now define
\[ \Gamma : \OpDISTMULT \rightarrow \OpMND({\CAT/\Set}) \]
with object map $(V,E) \mapsto (\ca GV,\Gamma{E})$ as before, and the rest of its definition is obtained by modifying the earlier definition of $\Gamma$ in what should now be the obvious way. The proof of the following result is obtained by a similar such modification of the proof of corollary(\ref{cor:2eq-mult-mon}).
\begin{corollary}\label{cor:2eq-mult-mon-dual}
$\Gamma$ and $\overline{(-)}$ provide a 2-equivalence \[ \OpDISTMULT \catequiv \OpPLMND. \]
\end{corollary}
%

\subsection{Monoidal monads and distributive laws}\label{ssec:Monmonad-Distlaw}
As explained in \cite{Str72} the assignment $\ca K \mapsto \MND(\ca K)$ is in fact the object map of a strict 3-functor. Just exploiting 2-functoriality here and corollaries(\ref{cor:2eq-mult-mon}) and (\ref{cor:2eq-mult-mon-dual}) one immediately obtains
\begin{theorem}\label{thm:monmon-distlaw}
$\Gamma$ and $\overline{(-)}$ provide two 2-equivalences of 2-categories:
\begin{enumerate}
\item ${\MND(\DISTMULT) \catequiv \MND(\PLMND)}$.\label{mdl1}
\item ${\MND(\OpDISTMULT) \catequiv \MND(\OpPLMND)}$.\label{mdl2}
\end{enumerate}
\end{theorem}
The meaning of this result is understood by understanding what the objects of the 2-categories involved are, that is to say, what monads are in each of the 2-categories $\DISTMULT$, $\OpDISTMULT$, $\PLMND$ and $\OpPLMND$.

A very beautiful observation of \cite{Str72} is that to give a monad on $(V,T)$ in $\MND(\ca K)$ is to give another monad $S$ on $V$, together with a distributive law $\lambda:TS{\rightarrow}ST$. Similarly to give a monad on $(V,T)$ in $\OpMND(\ca K)$ is to give another monad $S$ on $V$, together with a distributive law $\lambda:ST{\rightarrow}TS$ in the other direction. Thus the 2-categories $\MND(\MND(\ca K))$ and $\MND(\OpMND(\ca K))$ really have the same objects: such an object being a pair of monads on the same category and a distributive law between them. Thus both $\MND(\PLMND)$ and $\MND(\OpPLMND)$ are 2-categories whose objects are monad distributive laws between monads defined on categories of enriched graphs, with some extra conditions.

On the other hand a monad in the 2-category $\LaxAlg M$ of lax monoidal categories and lax monoidal functors is also a well-known thing, and such things are usually called \emph{monoidal monads}. Similarly an opmonoidal monad $T$ on a monoidal category $V$, that is to say a monad on $V$ in $\OpLaxAlg M$, comes with the extra data of coherence maps
\[ T(X_1 \tensor ... \tensor X_n) \rightarrow TX_1 \tensor ... \tensor TX_n \]
that are compatible with the monad structure. If for instance $\tensor$ is just cartesian product, then the product obstruction maps for $T$ endow it with an opmonoidal structure in a unique way.

By definition the objects of $\MND(\DISTMULT)$ are monoidal monads defined on distributive lax monoidal categories, and the objects of $\MND(\OpDISTMULT)$ are coproduct preserving opmonoidal monads defined on distributive lax monoidal categories. Thus the meaning of theorem(\ref{thm:monmon-distlaw}) is that it exhibits these kinds of monoidal and opmonoidal monads as being equivalent to certain kinds of distributive laws. We shall spell this out precisely in corollaries(\ref{cor:monmon-distlaw}) and (\ref{cor:opmonmon-distlaw}) below.

Let $(V,E)$ be a lax monoidal category and $T$ be a monad on $V$. In section(\ref{ssec:laxalg-const1}) we saw that when $T$ is a monoidal monad, that is to say one has coherence maps
\[ \tau_{X_i} : \opE\limits_i TX_i \rightarrow T \opE\limits_i X_i \]
making the underlying endofunctor of $T$ a lax monoidal functor and the unit and multiplication monoidal natural transformations, then one has another multitensor on $V$ given on objects by $T\opE\limits_i X_i$, and with unit and substitution given by the composites
\[ \xygraph{{\xybox{\xygraph{{X}="l" [r] {TX}="m" [r(1.5)] {TE_1X}="r" "l":"m"^-{\eta}:"r"^-{Tu}}}} [r(5)d(.1)]
{\xybox{\xygraph{!{0;(2,0):} {T\opE\limits_iT\opE\limits_j X_{ij}}="l" [r] {T^2\opE\limits_i\opE\limits_j X_{ij}}="m" [r] {T\opE\limits_{ij}X_{ij}}="r" "l":@<1ex>"m"^-{T{\tau}E}:@<1ex>"r"^-{\mu\sigma}}}}} \]
In particular if $E$ is distributive and $T$ preserves coproducts, then this new multitensor $TE$ is also distributive.

If instead $T$ has the structure of an opmonoidal monad, with the coherences
\[ \tau_{X_i} : T\opE\limits_i X_i \rightarrow \opE\limits_i TX_i \]
going in the other direction, then in the same way one can construct a new multitensor $ET$ on $V$ which on objects is defined by given by $\opE\limits_iTX_i$. Once again if $E$ is distributive and $T$ coproduct preserving, then $ET$ is a distributive multitensor. In particular when $E$ is cartesian product, $ET$ is the multitensor $T^{\times}$ of section(\ref{ssec:induction}).

With regards to monoidal monads, unpacking what theorem(\ref{thm:monmon-distlaw})(\ref{mdl1}) says at the object level gives
\begin{corollary}\label{cor:monmon-distlaw}
Let $(V,E)$ be a distributive lax monoidal category and $T$ be a monad on $V$. To give maps $\tau_{X_i} : \opE\limits_i TX_i \rightarrow T\opE\limits_i X_i$
making $T$ into a monoidal monad on $(V,E)$, is the same as giving a monad distributive law $\Gamma(E)\ca G(T){\rightarrow}\ca G(T)\Gamma(E)$ whose components are the identities on objects.
\end{corollary}
\noindent In the case where $T$ preserves coproducts one may readily verify that $\Gamma(TE) \iso \ca G(T)\Gamma(E)$ as monads, and so by corollary(\ref{cor:gamma-alg-ecat}) one understands what the algebras of this composite monad $\ca G(T)\Gamma(E)$ are.
\begin{corollary}\label{cor:algebras-GTGammaE}
If in the situation of corollary(\ref{cor:monmon-distlaw}) $T$ also preserves coproducts, then $\ca GV^{\ca G(T)\Gamma(E)} \iso \Enrich {TE}$.
\end{corollary}
\noindent Similarly, one can unpack what theorem(\ref{thm:monmon-distlaw})(\ref{mdl2}) says at the object level, witness $\Gamma(ET) \iso \Gamma(E)\ca G(T)$ and use corollary(\ref{cor:gamma-alg-ecat}) to conclude
\begin{corollary}\label{cor:opmonmon-distlaw}
Let $(V,E)$ be a distributive lax monoidal category and $T$ be a coproduct preserving monad on $V$. To give maps $\tau_{X_i} : T\opE\limits_i X_i \rightarrow \opE\limits_i TX_i$ making $T$ into an opmonoidal monad on $(V,E)$, is the same as giving a monad distributive law $\ca G(T)\Gamma(E){\rightarrow}\Gamma(E)\ca G(T)$ whose components are the identities on objects. Moreover $\ca GV^{\Gamma(E)\ca G(T)} \iso \Enrich {ET}$.
\end{corollary}
\begin{example}\label{ex:Cheng1}
From the inductive description of $\ca T_{\leq{n}}$ of section(\ref{ssec:induction}) and corollary(\ref{cor:opmonmon-distlaw}) one obtains a distibutive law
\[ \begin{array}{c} {\ca G(\ca T_{\leq{n}})\Gamma(\prod) \rightarrow \Gamma(\prod)\ca G(\ca T_{\leq{n}})} \end{array} \]
for all $n$, between monads on $\ca G^n\Set$, and the composite monad $\Gamma(\prod)\ca G(\ca T_{\leq{n}}) = \ca T_{\leq{(n{+}1)}}$. Thus we have recaptured the decomposition of \cite{ChengDist} of the strict $n$-category monad into a ``distributive series of monads''.
\end{example}
%

\subsection{The Trimble definition \`{a} la Cheng}\label{ssec:TCI}
Pursuing the idea of the previous example, we shall now begin to recover and in some senses generalise Cheng's analysis and description \cite{ChengCompOp} of the Trimble definition of weak $n$-category.

From \cite{EnHopI} example(2.6) non-symmetric operads in the usual sense can be regarded as multitensors. Here we shall identify a non-symmetric operad
\[ \begin{array}{lcr} {(E_n \, \, : \, \, n \in \N)} & {u:I \rightarrow E_1} &
{\sigma:E_k \tensor E_{n_1} \tensor ... \tensor E_{n_k} \rightarrow E_{n_{\bullet}}} \end{array} \]
in a braided monoidal category $V$, with the multitensor
\[ (X_1,...,X_n) \mapsto E_n \tensor X_1 \tensor ... \tensor X_n \]
it generates. Recall that one object $E$-categories for $E$ a non-symmetric operad are precisely algebras of the operad $E$ in the usual sense. If $\tensor$ is cartesian product, then the projections
\[ E_n \times X_1 \times ... \times X_n \rightarrow X_1 \times ... \times X_n  \]
are the components of a cartesian multitensor map $E \rightarrow \prod$. Conversely such a cartesian multitensor map exhibits $E$ as an operad via
\[ E_n := E(\underbrace{1,...,1}_{n}) \]
for all $n \in \N$.

Let $V$ be a distributive category and $T$ a coproduct preserving monad on $V$. Let us denote by $(E,\varepsilon)$ a non-symmetric operad in $V^T$. The ``$\varepsilon$'' is meant to denote the $T$-algebra actions, that is $\varepsilon_n:TE_n \rightarrow E_n$ is the $T$-algebra structure, and so $E$ denotes the underlying operad in $V$. Since $U^T$ preserves products it is the underlying functor of a strong monoidal functor $(V^T,(E,\varepsilon)) \rightarrow (V,E)$ between lax monoidal categories. Since the composites
\[ \xygraph{!{0;(3,0):} {T(E_n \times \prod\limits_iX_i)}="l" [r(1.2)] {TE_n \times \prod\limits_i TX_i}="m" [r] {E_n \times \prod\limits_iTX_i}="r" "l":@<1ex>"m"^-{\textnormal{prod. obstn.}}:@<1ex>"r"^-{\varepsilon_n \times \id}} \]
form the components of an opmonoidal structure for the monad $T$, we find ourselves in the situation of corollary(\ref{cor:opmonmon-distlaw}) and so obtain
\begin{proposition}\label{prop:monad-trimble}
Let $V$ be a distributive category, $T$ a coproduct preserving monad on $V$ and $(E,\varepsilon)$ a non-symmetric operad in $V^T$. Then one has a distributive law $\ca G(T)\Gamma(E) \rightarrow \Gamma(E)\ca G(T)$ between monads on $\ca GV$, and isomorphisms
\[ \Enrich {(E,\varepsilon)} \iso \Enrich {ET} \iso \ca GV^{\Gamma(E)\ca G(T)} \]
of categories over $\ca GV$.
\end{proposition}
\noindent This result has an operadic counterpart.
\begin{proposition}\label{prop:operad-trimble}
Let $V$ be a lextensive category, $T$ a cartesian and coproduct preserving monad on $V$, $\psi:S{\rightarrow}T$ a $T$-operad and $(E,\varepsilon)$ a non-symmetric operad in $V^S$. Then the monad $\Gamma(E)\ca G(S)$, whose algebras by proposition(\ref{prop:monad-trimble}) are $(E,\varepsilon)$-categories, has a canonical structure of a $\Gamma(T^{\times})$-operad.
\end{proposition}
\begin{proof}
With $\Gamma(T^{\times})=\Gamma(\prod)\ca G(T)$ we must exhibit a cartesian monad map $\Gamma(E)\ca G(S) \rightarrow \Gamma(\prod)\ca G(T)$. We have the cartesian multitensor map $\alpha:(E,\varepsilon) \rightarrow \prod$ which exhibits the multitensor $(E,\varepsilon)$ as a non-symmetric operad, thus $U^T\alpha$ is also a cartesian multitensor map, and since $V$ is lextensive $\Gamma$ sends this to a cartesian monad morphism. The required cartesian monad map is thus $\Gamma(U^T\alpha)\ca G(\psi)$.
\end{proof}
Recall the path-space functor $P:\Top \rightarrow \ca G(\Top)$ discussed in section(\ref{ssec:enriched-graphs}). To say that a non-symmetric topological operad $A$ acts on $P$ is to say that $P$ factors as
\[ \xygraph{!{0;(2,0):} {\Top}="l" [r] {\Enrich A}="m" [r] {\ca G(\Top)}="r" "l":"m"^-{P_A}:"r"^-{U^A}} \]
The main example to keep in mind is the version of the little intervals operad recalled in \cite{ChengCompOp} definition(1.1). As this $A$ is a contractible non-symmetric operad, $A$-categories may be regarded as a model of $A$-infinity spaces. Since $P$ is a right adjoint, $P_A$ is also a right adjoint by the Dubuc adjoint triangle theorem.

A product preserving functor
\[ Q : \Top \rightarrow V \]
into a distributive category, may be regarded as the underlying functor of a strong monoidal functor $(\Top,A) \rightarrow (V,QA)$ between lax monoidal categories. Applying $\Gamma$ to this gives us a monad functor
\[ (\ca G(\Top),\Gamma(A)) \rightarrow (\ca GV,\Gamma(QA)) \]
with underlying functor $\ca GQ$, which amounts to giving a lifting $\overline{Q}$ as indicated in the commutative diagram
\[ \xygraph{!{0;(2,0):(0,.5)::} {\Top}="tl" [r] {\Enrich A}="tm" [r] {\Enrich {QA}}="tr" [d] {\ca GV}="br" [l] {\ca G(\Top)}="bl" "tl"(:"tm"^-{P_A}(:"tr"^-{\overline{Q}}:"br"^{U^{QA}},:"bl":"br"_-{\ca G(Q)}),:"bl"_{P})} \]
and so we have produced another product preserving functor
\[ Q^{(+)} : \Top \rightarrow V^{(+)} \]
where $Q^{(+)}=\overline{Q}P_A$ and $V^{(+)}=\Enrich {QA}$. The functor $\overline{Q}$ is product preserving since $\ca G(Q)$ is and $U^{QA}$ creates products. The assignment
\[ (Q,V) \mapsto (Q^{(+)},V^{(+)}) \]
in the case where $A$ is as described in \cite{ChengCompOp} definition(1.1), is the inductive process lying at the heart of the Trimble definition. In this definition one begins with the connected components functor $\pi_0 : \Top \rightarrow \Set$ and defines the category $\Trimble 0$ of ``Trimble 0-categories'' to be $\Set$. The induction is given by
\[ (\Trimble {n{+}1},\pi_{n{+}1}) := (\textnormal{Trm}_{n}^{(+)},\pi_{n}^{(+)}) \]
and so this definition constructs not only a notion of weak $n$-category but the product preserving $\pi_n$'s to be regarded as assigning the fundamental $n$-groupoid to a space.

Applying proposition(\ref{prop:monad-trimble}) to this situation produces the monad on $n$-globular sets whose category of algebras is $\Trimble n$ as well as its decomposition into an iterative series of monads witnessed in \cite{ChengCompOp} section(4.2). Applying proposition(\ref{prop:operad-trimble}) and the inductive description of $\ca T_{\leq{n}}$ of section(\ref{ssec:induction}) exhibits these monads as $n$-operads.

\section{Lifting multitensors}\label{sec:lift-mult}

\subsection{Motivation}\label{ssec:lifting-intro}
Applied to the normalised $3$-operad for Gray categories \cite{Bat98}, the results of the section(\ref{sec:reexpress}) produce a lax monoidal structure $E$ on the category of $2$-globular sets whose enriched categories are exactly Gray categories. For this example it turns out that $E_1$ is $\ca T_{{\leq}2}$, and so providing a lift of $E$ in the sense of definition(\ref{def:lift}) amounts to the construction of a tensor product of $2$-categories whose enriched categories are Gray categories, that is to say, an abstractly constructed Gray tensor product. By the main result of this section theorem(\ref{thm:lift-mult}), \emph{every} $n$-multitensor has a lift which is unique given certain properties.

While the proof of theorem(\ref{thm:lift-mult}) is fairly abstract, and the uniqueness has the practical effect that in the examples we never have to unpack an explicit description of the lifted multitensors provided by the theorem, we provide such an unpacking in section(\ref{ssec:explicit-lifting}) anyway. This enables us to give natural conditions when the construction of the lifted multitensor is simpler. Doing all this requires manipulating some of the transfinite constructions that arise in monad theory, and we give a self-contained review of these in the appendix. 

\subsection{The multitensor lifting theorem}\label{ssec:lifting theorem}
In appendix \ref{sec:Dubuc} we recall an explicit description, for a given monad morphism $\phi:M \to S$ between accessible monads on a locally presentable category $V$, of the left adjoint $\phi_!$ to the canonical forgetful functor $\phi^*:V^S \to V^M$ induced by $\phi$. The key point about $\phi_!$ is that it is constructed via a transfinite process involving only \emph{connected} colimits in $V$. The importance of this is underscored by
\begin{lemma}\label{lem:concol-pathlike}
Let $V$ be a category with an initial object, $W$ be a cocomplete category, $J$ be a small connected category and
\[ F : J \rightarrow [\ca GV,\ca GW] \]
be a functor. Suppose that $F$ sends objects of $J$ to normalised functors, and arrows of $J$ to natural transformations whose components are identities on objects.
\begin{itemize}
\item[(1)]  Then the colimit $K:\ca GV{\rightarrow}\ca GW$ of $F$ may be chosen to be normalised.\label{cpl1}
\end{itemize}
Given such a choice of $K$:
\begin{itemize}
\item[(2)]  If $Fj$ is path-like for all $j \in J$, then $K$ is also path-like.\label{cpl2}
\item[(3)]  If $Fj$ is distributive for all $j \in J$, then $K$ is also distributive.\label{clp3}
\end{itemize}
\end{lemma}
\begin{proof}
Colimits in $[\ca GV,\ca GW]$ are computed componentwise from colimits in $\ca GW$ and so for $X \in \ca GV$ we must describe a universal cocone with components
\[ \kappa_{X,j} : Fj(X) \rightarrow KX. \]
We demand that the $\kappa_{X,j}$ are identities on objects. This is possible since the $Fj(X)_0$ form the constant diagram on $X_0$ by the hypotheses on $F$. For $a,b \in X_0$ we choose an arbitrary colimit cocone
\[ \{\kappa_{X,j}\}_{a,b} : Fj(X)(a,b) \rightarrow KX(a,b) \]
in $W$. One may easily verify directly that since $J$ is connected, the $\kappa_{X,j}$ do indeed define a univeral cocone for all $X$ in order to establish (1). Since the properties of path-likeness and distributivity involve only colimits at the level of the homs as does the construction of $K$ just given, (2) and (3) follow immediately since colimits commute with colimits in general.
\end{proof}
\noindent With these preliminaries in hand we are now ready to present the monad version of the multitensor lifting theorem, and then the lifting theorem itself.
\begin{lemma}\label{lem:mnd-lift-mult}
Let $V$ be a locally presentable category, $R$ be a coproduct preserving monad on $V$, $S$ be an accessible and normalised monad on $\ca GV$, and $\phi:\ca GR{\rightarrow}S$ be a monad morphism whose components are identities on objects. Denote by $T$ the monad induced by $\phi_! \ladj \phi^*$.
\begin{itemize}
\item[(1)]  One may choose $\phi_!$ so that $T$ becomes normalised.
\end{itemize}
Given such a choice of $\phi_!$:
\begin{itemize}
\item[(2)]  If $S$ is path-like then $T$ is path-like.
\item[(3)]  If $S$ is distributive then $T$ is distributive.
\end{itemize}
\end{lemma}
\begin{proof}
Let $\lambda$ be the regular cardinal such that $S$ preserves $\lambda$-filtered colimits. To verify that $(T,\eta^T,\mu^T)$ is a normalised monad one must verify: (i) $T$ is normalised, and (ii) the components of the unit $\eta^T$ are identities on objects. Since $\mu^T$ is a retraction of $\eta^TT$, it will then follow that the components of $\mu^T$ are also identities on objects. But $T$ is normalised iff $\ca GU^RT$ is normalised, and $\ca GU^RT=U^S\phi_!$, so for (i) it suffices to show that one can choose $\phi_!$ making $U^S\phi_!$ normalised. This follows by a transfinite induction using the explicit description of $U^S\phi_!$ of section(\ref{sec:Dubuc}) and lemma(\ref{lem:concol-pathlike}).

For the initial step note that $\phi_{!1}$ can be chosen to be normalised, because $S\eta^{\ca GR}$ is a componentwise-identity on objects natural transformation between normalised functors, since the monads $S$ and $\ca GR$ are normalised. Thus the coequaliser defining $\phi_{!1}$ is a connected colimit involving only normalised functors and componentwise-identity on objects natural transformations, and so $\phi_{!1}$ can be taken to be normalised by lemma(\ref{lem:concol-pathlike}). For the inductive steps the argument is basically the same: at each stage one is taking connected colimits of normalised functors and componentwise identity on objects natural transformations, so that by lemma(\ref{lem:concol-pathlike}) one stays within the subcategory of $[\ca G(V^R),\ca GV]$ consisting of such functors and natural transformations. Moreover using lemma(\ref{lem:concol-pathlike}) $T$ will be path-like if $S$ is. 

As for (ii) it suffices to prove that the components of $\ca G(U^R)\eta^T$ are identities on objects. Writing $q:S{\rightarrow}U^S\phi_{!}$ for transfinite composite constructed as part of the definition of $\phi_!$ (note that $U^S\phi_{!}=\ca G(U^R)T$ by definition) recall from the end of section(\ref{sec:Dubuc}) that one has a commutative square
\[ \xymatrix{{\ca G(R)U^{\ca GR}} \ar[r]^-{\ca G\rho} \ar[d]_{{\phi}U^M} & {U^{\ca GR}} \ar[d]^{\ca G(U^R)\eta^T} \\ {SU^{\ca GR}} \ar[r]_-{q} & {U^S\phi_{!}}} \]
where $\rho$ is the 2-cell datum for $R$'s Eilenberg-Moore object, which we recall is preserved by $\ca G$. Now $\rho$ is componentwise the identity objects since $\eta^{\ca GR}_X$ is and $\rho$ is a retraction of it, $U^S\phi$ is the identity on objects by definition, and $q$ is by construction, so the result follows.
\end{proof}
\noindent Recall from definition(\ref{def:lift}) that a \emph{lift} of $(E,u,\sigma)$ is a normal multitensor $(E',\id,\sigma')$ on $V^{E_1}$ together with an isomorphism $\Enrich E \iso \Enrich {E'}$ which commutes with the forgetful functors into $\ca G(V^{E_1})$. When in addition $E'$ is distributive, we say that it is a \emph{distributive lift} of $E$. Recall from \cite{Str72} that for any category $V$, the functor
\[ \textnormal{Alg} : \op {\Mnd(V)} \rightarrow \CAT/V \]
which sends a monad $T$ on $V$ to the forgetful functor $U^T:V^T{\rightarrow}V$, is fully-faithful.
\begin{theorem}\label{thm:lift-mult}
Let $(E,u,\sigma)$ be a distributive multitensor on $V$ a locally presentable category, and let $E$ be accessible in each variable. Then $E$ has a distributive lift $E'$, which is unique up to isomorphism.
\end{theorem}
\begin{proof}
Write $\ca SE$ for the distributive multitensor on $V$ whose unary part is $E_1$ and whose non-unary parts are constant at $\emptyset$. There is an obvious inclusion $\psi:\ca SE{\rightarrow}E$ of multitensors and one clearly has
\[ \Enrich {\ca SE} \iso \ca G(V^{E_1}) \]
Applying lemma(\ref{lem:mnd-lift-mult}) with $S=\Gamma{E}$, $R=E_1$ and $\phi=\Gamma{\psi}$ one produces a path-like, normalised and distributive monad $T$ on $\ca G(V^ {E_1})$, because $S$ is accessible by proposition(\ref{prop:Gamma-accessible}). Thus one has a distributive multitensor $\overline{T}$ on $V^{E_1}$. Applying proposition(\ref{prop:pl-alg<->cat}) to $\overline{T}$, and corollary(\ref{cor:gamma-alg-ecat}) to $E$, gives
\[ \Enrich {\overline{T}} \iso \Enrich E \]
in view of the monadicity of $\phi^*$. That is to say, $\overline{T}$ is a distributive lift of $E$. As for uniqueness suppose that $(E',\id,\sigma')$ is a distributive lift of $E$. Then by corollary(\ref{cor:gamma-alg-ecat}) and proposition(\ref{prop:pl-adjoint-char}), $\Gamma(E')$ is a distributive monad on $\ca G(V^{E_1})$ and one has
\[ \ca G(V^{E_1})^{\Gamma(E')} \iso \Enrich E \]
commuting with the forgetful functors into $\ca G(V^{E_1})$. By the fully-faithfulness of $\textnormal{Alg}$ recalled above, one has an isomorphism $\Gamma(E'){\iso}T$ of monads, and thus by proposition(\ref{prop:pl-adjoint-char}) an isomorphism $E'{\iso}\overline{T}$ of multitensors.
\end{proof}
Applying this result to any normalised $(n+1)$-operad $A$, exhibits its algebras as categories enriched in the algebras of some $n$-operad. The $n$-operad is $\Gamma(A)_1$, and tensor product over which one enriches is $\Gamma(A)'$. In cases where we already know what our tensor product ought to be, the uniqueness part of theorem(\ref{thm:lift-mult}) ensures that it is. An instance of this is
\begin{example}\label{ex:Gray}
In \cite{Bat98} the normalised 3-operad $G$ whose algebras are Gray categories was constructed. As we have already seen, $\Gamma(G)$ is a lax monoidal structure on $\ca G^2(\Set)$ whose enriched categories are Gray categories, and $\Gamma(G)_1$ is the operad for strict 2-categories. Note that the usual Gray tensor product is symmetric monoidal closed and thus distributive. Thus by theorem(\ref{thm:lift-mult}) $\Gamma(G)'$ is the Gray tensor product. In other words, the general methods of this paper have succeeded in producing the Gray tensor product of $2$-categories from the operad $G$.
\end{example}
\noindent More generally given a distributive tensor product $\tensor$ on the category of algebras of an $n$-operad $B$, and a normalised $(n+1)$-operad $A$ whose algebras are the categories enriched in $B$-algebras, theorem(\ref{thm:lift-mult}) exhibits $\tensor$ as the more generally constructed $\Gamma(A)'$.
\begin{example}\label{ex:Crans}
In \cite{Crans99} Sjoerd Crans explicitly constructed a tensor product on the category of Gray-categories. This explicit construction was extremely complicated. It is possible to exhibit the Crans tensor product as an instance of our general theory, by rewriting his explicit constructions as the construction of the 4-operad whose algebras are teisi in his sense. The multitensor $E$ associated to this 4-operad has $E_1$ equal to the 3-operad for Gray categories. Thus theorem(\ref{thm:lift-mult}) constructs a lax tensor product of Gray categories whose enriched categories are teisi. Since the tensor product explicitly constructed by Crans is distributive, the uniqueness of part of theorem(\ref{thm:lift-mult}) ensures that it is indeed $E'$, since teisi are categories enriched in the Crans tensor product by definition.
\end{example}
\noindent Honestly writing the details of the 4-operad of example(\ref{ex:Crans}) is a formidable task and we have omitted this here. In the end though, such details will not be important, because such a tensor product (or more properly a biclosed version thereof) will only be really useful once it is given a conceptual definition.

\subsection{Applications of the lifting theorem}\label{ssec:A-infinity}
Let $V$ be a symmetric monoidal model category which satisfies the conditions of \cite{BergerMoerdijk} or the monoid axiom of \cite{SS}. In this case the category of pruned $n$-operads of \cite{Bat03} can be equipped with a monoidal model structure \cite{BataninBergerCisinski}. So we can speak of cofibrant $n$-operads in $V$.

For $n=1$ let us fix a particular cofibrant and contractible $1$-operad $A$. The algebras of $A_1$ can be called \emph{$A_{\infty}$-categories enriched in $V$}. Up to homotopy the choice of $A_1$ is not important. So we can speak of \emph{the} category of $A_{\infty}$-categories. For $n=2$ we denote by $A_2$ a cofibrant contractible $2$-operad in $V$. Let  $B= (A_2,u,\sigma)$ be the corresponding multitensor on $\ca GV$. One can always choose $A_2$ in such a way that its unary part is $A_1$. As in \cite{EnHopI} for an arbitary multitensor $E$, one object $E$-categories are called $E$-monoids. Similarly one object $A_{\infty}$-categories are called $A_{\infty}$-monoids.
\begin{theorem}\label{thm:A-infinity-app}
\begin{enumerate}
\item  There is a distributive lift $B'$ of $B$ to the category of $A_{\infty}$-categories.
\item  $B'$ restricts to give a multitensor $C$ on the category of $A_{\infty}$-monoids.
\item  The category of $C$-monoids is equivalent to the category of algebras of $\sym_2(A_2)$ and therefore is Quillen equivalent to the category of the algebras of the little squares operad.
\end{enumerate}
\end{theorem}
\begin{proof}
The first statement is a direct consequence of theorem(\ref{thm:lift-mult}). The second statement follows from the fact that $B'$ is the 
cartesian product on the object level. The last statement follows from the theorem(8.6) of \cite{Bat03}.
\end{proof}
\noindent Applying this result to the case $V=\Cat$ with its folklore model structure one recovers
\begin{corollary}\label{cor:JS-braided}
[Joyal-Street] The category of braided monoidal categories is equivalent to the category of monoidal categories equipped with multiplication.
\end{corollary}
\noindent The previous corollary is proved by Joyal and Street \cite{JS93} by a direct application of a ``categorified'' Eckmann-Hilton argument. The following analogous result for $2$-categories appears to be new.
\begin{corollary}\label{cor:coh-bm2c}
The category of braided monoidal $2$-categories is equivalent to the category of Gray-monoids with multiplication.
\end{corollary}
\begin{proof}
Apply theorem(\ref{thm:A-infinity-app}) with $V=\Enrich 2$ equipped with the Gray tensor product and Lack's folklore model structure for 2-categories \cite{LackFolk2}.
\end{proof}
Thus theorem(\ref{thm:A-infinity-app}) should be considered as an $\infty$-generalisation of the above corollaries. We believe it sheds some light on the problem of defining the tensor product of $A_{\infty}$-algebras initiated by \cite{Saneblidze}. As explained in the introduction, the negative result of \cite{Markl} shows that there is no hope to get an ``honest'' tensor product of such algebras. Thus the multitensor $C$ constructed in theorem(\ref{thm:A-infinity-app}) is genuinely lax, and exhibits laxity as a way around the aforementioned negative result. In future work we will generalise this theorem to arbitrary dimensions.

\subsection{Unpacking $E'$}\label{ssec:explicit-lifting} 
Let us now instantiate the constructions of section(\ref{sec:Dubuc}) to produce a more explicit description of the lifted multitensor $E'$. Beyond mere instantiation this task amounts to reformulating everything in terms of hom maps which live in $V$, because in our case the colimits being formed in $\ca GV$ at each stage of the construction are connected colimits diagrams whose morphisms are all identity on objects. Moreover these fixed object sets are of the form $\{0,...,n\}$ for $n \in \N$.
\\ \\
{\bf Notation}. We shall be manipulating sequences of data and so we describe here some notation that will be convenient. A sequence $(a_1,...,a_n)$ from some set will be denoted more tersely as $(a_i)$ leaving the length unmentioned. Similarly a sequence of sequences
\[ ((a_{11},...,a_{1n_1}),...,(a_{k1},...,a_{kn_k})) \]
of elements from some set will be denoted $(a_{ij})$ -- the variable $i$ ranges over $1{\leq}i{\leq}k$ and the variable $j$ ranges over $1{\leq}j{\leq}n_i$. Triply-nested sequences look like this $(a_{ijk})$, and so on. These conventions are more or less implicit already in the notation we have been using all along for multitensors. See especially section(\ref{ssec:LMC}) and \cite{EnHopI}. We denote by
\[ \con(a_{i_1,...,i_k}) \]
the ordinary sequence obtained from the $k$-tuply nested sequence $(a_{i_1,...,i_k})$ by concatenation. In particular given a sequence $(a_i)$, the set of $(a_{ij})$ such that $\con(a_{ij})=(a_i)$ is just the set of partitions of the original sequence into doubly-nested sequences, and will play an important role below. This is because to give the substitution maps for a multitensor $E$ on $V$, is to give maps
\[ \sigma : \opE\limits_i\opE\limits_j X_{ij} \to \opE\limits_i X_i \]
for all $(X_{ij})$ and $(X_i)$ from $V$ such that $\con(X_{ij})=(X_i)$.
\\ \\
\indent The monad map $\phi:M \to S$ is taken as $\Gamma(\psi):\ca GE_1 \to \Gamma(E)$ where $\psi:E_1 \to E$ is the inclusion of the unary part of the multitensor $E$. Note the notational abuse -- we regard write $E_1$ for the multitensor on $V$ obtained from $E$ by ignoring (ie setting to constant at $\emptyset$) the non-unary parts, but also as the monad on $V$ -- and so $\ca GE_1=\Gamma(E_1)$ as monads. The role of $(X,x)$ in $V^M$ is played by sequences $(X_i,x_i)$ of $E_1$-algebras regarded as objects of $\ca GV^{E_1}$ as in section(\ref{ssec:DefNMonad}).

The transfinite induction produces for each ordinal $m$ and each sequence of $E_1$-algebras as above of length $n$, morphisms
\[ \begin{array}{l} {v^{(m)}_{(X_i,x_i)} : SQ_m(X_i,x_i) \rightarrow Q_{m{+}1}(X_i,x_i)} \\ {q^{(m)}_{(X_i,x_i)}:Q_m(X_i,x_i) \to Q_{m{+}1}(X_i,x_i)} \\ {q^{({<}m)}_{(X_i,x_i)}:S(X_i) \to Q_m(X_i,x_i)} \end{array} \]
in $\ca GV$ which are identities on objects, and thus we shall now evolve this notation so that it only records what's going on in the hom between $0$ and $n$. By the definition of $S$ we have the equation on the left
\[ \begin{array}{lccr} {S(X_i)(0,n) = \opE\limits_iX_i} &&& {Q_m(X_i,x_i)(0,n) = {\opEm\limits_i}(X_i,x_i)} \end{array} \]
and the equation on the right is a definition. Because of these definitions and that of $S$ we have the equation
\[ SQ_m(X_i,x_i)(0,n) = \coprod\limits_{\con(X_{ij},x_{ij}){=}(X_i,x_i)} \opE\limits_i\opEm\limits_j (X_{ij},x_{ij}). \]
The data for the hom maps of the $v^{(m)}$ thus consists of morphisms
\[ \begin{array}{c} {v^{(m)}_{(X_{ij},x_{ij})} : \opE\limits_i\opEm\limits_j (X_{ij},x_{ij}) \to \opEmpone\limits_i (X_i,x_i)} \end{array} \]
in $V$ whenever one has $\con(X_{ij},x_{ij})=(X_i,x_i)$ as sequences of $E_1$-algebras.

To summarise, the output of the transfinite process we are going to describe is, for each ordinal $m$, the following data. For each sequence $(X_i,x_i)$ of $E_1$-algebras, one has an object
\[ \opEm\limits_i (X_i,x_i) \]
and morphisms
\[ \begin{array}{l} {v^{(m)}_{(X_{ij},x_{ij})}} : {\opE\limits_i\opEm\limits_j (X_{ij},x_{ij}) \to \opEmpone\limits_i (X_i,x_i)} \\
{q^{(m)}_{(X_i,x_i)}} : {\opEm\limits_i (X_i,x_i) \to \opEmpone\limits_i (X_i,x_i)} \\
{q^{(<m)}_{(X_i,x_i)}} : {\opE\limits_iX_i \to \opEm\limits_i (X_i,x_i)} \end{array} \]
of $V$ where $\con(X_{ij},x_{ij})=(X_i,x_i)$.
\\ \\
{\bf Initial step}. First we put $\opEzero\limits_i (X_i,x_i) = \opE\limits_i X_i$, $q^{({<}0)}_{(X_i,x_i)_i} = \id$, and then form the coequaliser
\begin{equation}\label{eq:coeq}
\xygraph{!{0;(2,0):} {\opE\limits_iE_1X_i}="l" [r] {\opE\limits_iX_i}="m" [r] {\opEone\limits_i(X_i,x_i)}="r" "l":@<2ex>"m"^-{\sigma} "l":"m"_-{\opE\limits_ix_i}:@<1ex>"r"^-{q^{(0)}_{(X_i,x_i)}}} \end{equation}
in $V$ to define $q^{(0)}$. Put $v^{(0)}=q^{(0)}\sigma$ and $q^{({<}1)}=q^{(0)}$.
\\ \\
{\bf Inductive step}. Assuming that $v^{(m)}$, $q^{(m)}$ and $q^{({<}m{+}1)}$ are given, we have maps
\[ \begin{array}{lcr} {\xybox{\xygraph{!{0;(2,0):} {\opE\limits_i\opE\limits_j\opEm\limits_k}="l" [r] {\opE\limits_i\opEmpone\limits_{jk}}="r" "l":"r"^-{\opE\limits_iv^{(m)}}}}} && {\xybox{\xygraph{!{0;(2,0):} {\opE\limits_i\opE\limits_j\opEm\limits_k}="l" [r] {\opE\limits_{ij}\opEm\limits_k}="m" [r] {\opE\limits_{ij}\opEmpone\limits_k}="r" "l":"m"^-{{\sigma}\opEm\limits_k}:"r"^-{q^{(m)}}}}} \end{array} \]
and these are used to provide the parallel maps in the coequaliser
\[ \xygraph{!{0;(1.5,0):} {\coprod\limits_{\con(X_{ijk},x_{ijk})=(X_i,x_i)} \opE\limits_i\opE\limits_j\opEm\limits_k(X_{ijk},x_{ijk})}="l" [d] {\coprod\limits_{\con(X_{ij},x_{ij})=(X_i,x_i)} \opE\limits_i\opEmpone\limits_j(X_{ij},x_{ij})}="m" [d] {\opEmptwo\limits_i (X_i,x_i)}="r" "l":@<-2ex>"m" "l":@<2ex>"m":"r"^-{(v^{(m{+}1)}_{(X_{ij},x_{ij})})}} \]
which defines the $v^{(m{+}1)}$, the commutative diagram
\[ \xygraph{{\opEmpone\limits_i (X_i,x_i)}="l" [d] {E_1\opEmpone\limits_i (X_i,x_i)}="il" [r(4)] {\coprod\limits_{\con(X_{ij},x_{ij})=(X_i,x_i)} \opE\limits_i\opEmpone\limits_j(X_{ij},x_{ij})}="ir" [ru] {\opEmptwo\limits_i (X_i,x_i)}="r" "l":"il"_-{u}:"ir":"r"^(.35){v^{(m{+}1)}_{(X_i,x_i)}}:@{<-}"l"_-{q^{(m{+}1)}_{(X_i,x_i)}}} \]
in which the unlabelled map is the evident coproduct inclusion defines $q^{(m{+}1)}$, and $q^{({<}m{+}2)}=q^{(m{+}1)}q^{({<}m{+}1)}$.
\\ \\
{\bf Limit step}. Define $\opEm\limits_i (X_i,x_i)$ as the colimit of the sequence given by the objects $\opEr\limits_i (X_i,x_i)$ and morphisms $q^{(r)}$ for $r < m$, and $q_{<{m}}$ for the component of the universal cocone at $r=0$.
\[ \xygraph{!{0;(0,-1.5):(0,-3.45)::} {\colsum\limits_{\con(X_{ijk},x_{ijk})=(X_i,x_i)} \opE\limits_i\opE\limits_j\opEr\limits_k(X_{ijk},x_{ijk})}="tl" [r] {\colsum\limits_{\con(X_{ij},x_{ij})=(X_i,x_i)} \opE\limits_i\opEr\limits_j(X_{ij},x_{ij})}="tm" [r] {\colim_{r{<}m} \opEr\limits_i (X_i,x_i)}="tr" [d] {\opEm\limits_i (X_i,x_i)}="br" [l] {\coprod\limits_{\con(X_{ij},x_{ij})=(X_i,x_i)} \opE\limits_i\opEm\limits_j(X_{ij},x_{ij})}="bm" [l] {\coprod\limits_{\con(X_{ijk},x_{ijk})=(X_i,x_i)} \opE\limits_i\opE\limits_j\opEm\limits_k(X_{ijk},x_{ijk})}="bl" "tl":@<1ex>"tm"^-{\sigma^{(<{m})}}:@<1ex>"tr"^-{v^{(<{m})}} "tl":@<-1ex>"tm"_-{(Ev)^{(<{m})}}:@<-1ex>@{<-}"tr"_-{u^{(<{m})}} "bl":"bm"_-{\mu}:@{<-}"br"_-{uc} "tl":"bl"^{o_{m,2}} "tm":"bm"^{o_{m,1}} "tr":@{=}"br"} \]
As before we write $o_{m,1}$ and $o_{m,2}$ for the obstruction maps, and $c$ denotes the evident coproduct injection. The maps $\sigma^{(<{m})}$, $(Ev)^{(<{m})}$, $v^{(<{m})}$ and $u^{(<{m})}$ are by definition induced by $\sigma{\opEr}$, $(Ev)^{(r)}$, $v^{(r)}$ and $u{\opEr}$ for $r < m$ respectively. Define $v^{(m)}$ as the coequaliser of $o_{m,1}\sigma^{(<{m})}$ and $o_{m,1}(Ev)^{<{m}}$, $q^{(m)}=v^{(m)}(u{\opEm})$ and $q^{(<{m{+}1})}=q^{(m)}q^{(<{m})}$.
\\ \\
Instantiating corollary(\ref{cor:explicit-phi-shreik}) to the present situation gives
\begin{corollary}\label{cor:lifted-obj}
Let $V$ be a locally presentable category, $\lambda$ a regular cardinal, and $E$ a distributive $\lambda$-accessible multitensor on $V$. Then for any ordinal $m$ with $|m| \geq \lambda$ one may take
\[ (\opEm\limits_i (X_i,x_i), a(X_i,x_i)) \]
where the action $a(X_i,x_i)$ is given as the composite
\[ \xygraph{!{0;(3,0):} {E_1\opEm\limits_i (X_i,x_i)}="l" [r] {\opEmpone\limits_i (X_i,x_i)}="m" [r] {\opEm\limits_i (X_i,x_i)}="r" "l":"m"^-{v^{(m)}}:"r"^-{(q^{(m)})^{-1}}} \]
as an explicit description of the object map of the lifted multitensor $E'$ on $V^{E_1}$.
\end{corollary}
In corollaries (\ref{cor:phi-shreik-simple}) and (\ref{cor:vexp-simple}), in which the initial data is a monad map $\phi:M \to S$ between monads on a category $V$ together with an algebra $(X,x)$ for $M$, we noted the simplification of our constructions when $S$ and $S^2$ preserve the coequaliser
\begin{equation}\label{eq:monad-coeq} \xygraph{!{0;(2,0):} {SMX}="l" [r] {SX}="m" [r] {Q_1X}="r" "l":@<-1ex>"m"_-{Sx} "l":@<1ex>"m"^-{\mu^SS(\phi)}:"r"^-{q_0}} \end{equation}
in $V$, which is part of the first step of the inductive construction of $\phi_!$. In the present situation the role of $V$ is played by the category $\ca GV$, the role of $S$ is played by $\Gamma E$, and the role of $(X,x)$ played by a given sequence $(X_i,x_i)$ of $E_1$-algebras, and so the role of the coequaliser (\ref{eq:monad-coeq}) is now played by the coequaliser
\begin{equation}\label{eq:monad-coeq2} \xygraph{!{0;(2,0):} {\Gamma E(E_1X_i)}="l" [r(1.5)] {\Gamma E(X_i)}="m" [r] {Q_1}="r" "l":@<-1ex>"m"_-{Sx} "l":@<1ex>"m"^-{\mu^SS(\phi)}:"r"^-{q^{(0)}}} \end{equation}
in $\ca GV$. Here we have denoted by $Q_1$ the $V$-graph with objects $\{0,...,n\}$ and homs given by
\[ Q_1(i,j) = \left\{\begin{array}{lll} {\emptyset} && {\textnormal{if $i>j$}} \\ {\opEone\limits_{i{<}k{\leq}j}(X_k,x_k)} && {\textnormal{if $i \leq j$.}} \end{array}\right. \]
Taking the hom of (\ref{eq:monad-coeq2}) between $0$ and $n$ gives the coequaliser
\begin{equation}\label{eq:mult-coeq} \xygraph{!{0;(2,0):} {\opE\limits_iE_1X_i}="l" [r] {\opE\limits_iX_i}="m" [r] {\opEone\limits_i(X_i,x_i)}="r" "l":"m"_-{\opE\limits_ix_i} "l":@<2ex>"m"^-{\sigma}:@<1ex>"r"^-{q^{(0)}}} \end{equation}
in $V$ which is part of the first step of the explicit inductive construction of $E'$. We shall refer to (\ref{eq:mult-coeq}) as the \emph{basic coequaliser associated to the sequence $(X_i,x_i)$} of $E_1$-algebras. Note that all coequalisers under discussion here are reflexive coequalisers, with the common section for the basic coequalisers given by the maps $\opE\limits_iu_{X_i}$.

The basic result which expresses why reflexive coequalisers are nice, is the $3{\times}3$-lemma, which we record here for the reader's convenience. A proof can be found in \cite{PTJ-topos77}.
\begin{lemma}\label{lem:3by3}
{\bf $3{\times}3$-lemma}. Given a diagram
\[ \xymatrix @R=3em @C=3em {A \ar@<1ex>[r]^-{f_1} \ar@<-1ex>[r]_{g_1} \ar@<1ex>[d]^{b_1} \ar@<-1ex>[d]_{a_1}
& B \ar[r]^-{h_1} \ar@<1ex>[d]^{b_2} \ar@<-1ex>[d]_{a_2} & C \ar@<1ex>[d]^{b_3} \ar@<-1ex>[d]_{a_3} \\
D \ar@<1ex>[r]^-{f_2} \ar@<-1ex>[r]_-{g_2} & E \ar[r]^-{h_2} & F \ar[d]^{c} \\ && H} \]
in a category such that: (1) the two top rows and the right-most column are coequalisers, (2) $a_1$ and $b_1$ have a common section, (3) $f_1$ and $g_1$ have a common section, (3) $f_2a_1{=}a_2f_1$, (4) $g_2b_1{=}b_2g_1$, (5) $h_2a_2{=}a_3h_1$ and (6) $h_2b_2{=}b_3h_1$; then $ch_2$ is a coequaliser of $f_2a_1{=}a_2f_1$ and $g_2b_1{=}b_2g_1$.
\end{lemma}
\noindent If $F:{\ca A_1}{\times}...{\times}{\ca A_n}{\rightarrow}{\ca B}$ is a functor which preserves connected colimits of a certain type, then it also preserves these colimits in each variable separately, because for a connected colimit, a cocone involving only identity arrows is a universal cocone. The most basic corollary of the $3{\times}3$-lemma says that the converse of this is true for reflexive coequalisers.
\begin{corollary}\label{cor:3by3}
Let $F:{\ca A_1}{\times}...{\times}{\ca A_n}{\rightarrow}{\ca B}$ be a functor. If $F$ preserves reflexive coequalisers in each variable separately then $F$ preserves reflexive coequalisers.
\end{corollary}
\noindent and this can be proved by induction on $n$ using the $3{\times}3$-lemma in much the same way as \cite{LkMonFinMon} lemma(1). The most well-known instance of this is
\begin{corollary}\label{cor:3by3-2}\cite{LkMonFinMon}
Let $\ca V$ be a biclosed monoidal category. Then the $n$-fold tensor product of reflexive coequalisers in $\ca V$ is again a reflexive coequaliser.
\end{corollary}
\noindent In particular note that by corollary(\ref{cor:3by3}) a multitensor $E$ preserves (some class of) reflexive coequalisers iff it preserves them in each variable separately.

Returning to our basic coequalisers an immediate consequence of the explicit description of $\Gamma E$ and corollary(\ref{cor:3by3}) is
\begin{lemma}\label{lem:reformulate-simplifying-conditions}
Let $E$ be a distributive multitensor on $V$ a cocomplete category, and $(X_i,x_i)$ a sequence of $E_1$-algebras. If $E$ preserves the basic coequalisers associated to all the subsequences of $(X_i,x_i)$, then for all $r \in \N$, $(\Gamma E)^r$ preserves the coequaliser (\ref{eq:monad-coeq2}).
\end{lemma}
\noindent and applying this lemma and corollary(\ref{cor:phi-shreik-simple}) gives
\begin{corollary}\label{cor:lifted-obj-simple}
Let $V$ be a locally presentable category, $\lambda$ a regular cardinal, $E$ a distributive $\lambda$-accessible multitensor on $V$ and $(X_i,x_i)$ a sequence of $E_1$-algebras. If $E$ preserves the basic coequalisers associated to all the subsequences of $(X_i,x_i)$, then one may take
\[ \opEpr\limits_i(X_i,x_i) = (\opEone\limits_i (X_i,x_i), a) \]
where the action $a$ is defined as the unique map such that $aE_1(q^{(0)})=q^{(0)}\sigma$.
\end{corollary}
\noindent Note in particular that when the sequence $(X_i,x_i)$ of $E_1$-algebras is of length $n=0$ or $n=1$, the associated basic coequaliser is absolute. In the $n=0$ case the basic coequaliser is constant at $E_0$, and when $n=1$ the basic coequaliser may be taken to be the canonical presentation of the given $E_1$-algebra. Thus in these cases it follows from corollary(\ref{cor:lifted-obj-simple}) that $E'_0=(E_0,\sigma)$ and $E_1'(X,x)=(X,x)$. Reformulating the explicit description of the unit in corollary(\ref{cor:vexp-simple}) one recovers the fact from our explicit descriptions, that the unit of $E'$ is the identity, which was of course true by construction.

To complete the task of giving a completely explicit description of the multitensor $E'$ we now turn to unpacking its substitution. So we assume that $E$ is a distributive $\lambda$-accessible multitensor on $V$ a locally presentable category, and fix an ordinal $m$ so that $|m| \geq \lambda$, so that $E'$ may be constructed as $E^{(m)}$ as in corollary(\ref{cor:lifted-obj}). By transfinite induction on $r$ we shall generate the following data:
\[ \sigma^{(r)}_{X_{ij},x_{ij}} : \opEr\limits_i(\opEm\limits_j(X_{ij}),x_{ij}) \to \opEm\limits_i(X_i,x_i) \]
and $\sigma^{(r{+}1)}_{X_{ij},x_{ij}}$ whenever $\con(X_{ij},x_{ij})=(X_i,x_i)$, such that
\[ \xygraph{!{0;(2.5,0):(0,.5)::} {\opE\limits_i\opEr\limits_j\opEm\limits_k}="tl" [r] {\opErpone\limits_{ij}\opEm\limits_k}="tr" [d] {\opEm\limits_{ijk}}="br" [l] {\opE\limits_i\opEm\limits_{jk}}="bl" "tl":"tr"^-{v^{(r)}E^{(m)}}:"br"^-{\sigma^{(r{+}1)}}:@{<-}"bl"^-{(q^{(m)})^{-1}v^{(m)}}:@{<-}"tl"^-{\opE\limits_i\sigma^{(r)}}} \]
commutes.
\\ \\
{\bf Initial step}. Define $\sigma^{(0)}$ to be the identity and $\sigma^{(1)}$ as the unique map such that $\sigma^{(1)}q^{(0)}=(q^{(m)})^{-1}v^{(m)}$ by the universal property of the coequaliser $q^{(0)}$.
\\ \\
{\bf Inductive step}. Define $\sigma^{(r{+}2)}$ as the unique map such that
\[ \sigma^{(r{+}2)}(v^{(r{+}1)}E^{(m)})=(q^{(m)})^{-1}v^{(m)}(\opE\limits_i\sigma^{(r{+}1)}) \]
using the universal property of $v^{(r{+}1)}$ as a coequaliser.
\\ \\
{\bf Limit step}. When $r$ is a limit ordinal define $\sigma^{(r)}$ as induced by the $\mu^{(s)}$ for $s<r$ and the universal property of $E^{(r)}$ as the colimit of the sequence of the $E^{(s)}$ for $s<r$. Then define $\sigma^{(r{+}1)}$ as the unique map such that
\[ \sigma^{(r{+}1)}(v^{(r)}E^{(m)})=(q^{(m)})^{-1}v^{(m)}(\opE\limits_i\sigma^{(r)}) \]
using the universal property of $v^{(r)}$ as a coequaliser.
\\ \\
The fact that the transfinite construction just specified was obtained from that for corollary(\ref{cor:induced-monad-very-explicit}), by taking $S=\Gamma E$ and looking at the homs, means that by corollaries (\ref{cor:induced-monad-very-explicit}) and (\ref{cor:vexp-simple}) one has
\begin{corollary}\label{cor:induced-substitution-very-explicit}
Let $V$ be a locally presentable category, $\lambda$ a regular cardinal, $E$ a distributive $\lambda$-accessible multitensor on $V$ and $(X_i,x_i)$ a sequence of $E_1$-algebras. Then one has
\[ \sigma'_{(X_i,x_i)} = \sigma^{(m)}_{(X_i,x_i)} \]
as an explicit description of the substitution of $E'$. If moreover $E$ preserves the basic coequalisers of all the subsequences of $(X_i,x_i)$, then one may take $\sigma^{(1)}_{(X_i,x_i)}$ as the explicit description of the substitution.
\end{corollary}
%

\subsection{Functoriality of lifting}\label{ssec:functoriality-lifting}
Recall \cite{Str72} \cite{LS00} that when $\ca K$ has Eilenberg-Moore objects, the one and 2-cells of the 2-category $\MND(\ca K)$ admit another description. Given monads $(V,T)$ and $(W,S)$ in $\ca K$, and writing $U^T:V^T{\rightarrow}V$ and $U^S:W^S{\rightarrow}W$ for the one-cell data of their respective Eilenberg-Moore objects, to give a monad functor $(H,\psi):(V,T){\rightarrow}(W,S)$, is to give $H$ and $\tilde{H}:V^T{\rightarrow}W^S$ such that $U^S\tilde{H}=HU^T$. This follows immediately from the universal property of Eilenberg-Moore objects. Similarly to give a monad 2-cell $\phi:(H_1,\psi_1){\rightarrow}(H_2,\psi_2)$ is to give $\phi:H_1{\rightarrow}H_2$ and $\tilde{\phi}:\tilde{H_1}{\rightarrow}\tilde{H_2}$ commuting with $U^T$ and $U^S$. Note that Eilenberg-Moore objects in $\CAT/\Set$ are computed as in $\CAT$, and we shall soon apply these observations to the case $\ca K = \CAT/\Set$.

Suppose we have a lax monoidal functor $(H,\psi):(V,E){\rightarrow}(W,F)$, that $V$ and $W$ are locally presentable, and that $E$ and $F$ are accessible. Then we obtain a commutative diagram
\[ \xygraph{!{0;(1.5,0):(0,.667)::} {\Enrich E}="tl" [r] {\ca GV^{E_1}}="tm" [r] {\ca GV}="tr" [d] {\ca GW}="br" [l] {\ca GW^{F_1}}="bm" [l] {\Enrich F}="bl" "tl":"tm":"tr" "bl":"bm":"br" "tl":"bl" "tm":"bm" "tr":"br"} \]
of forgetful functors in $\CAT/\Set$. Applying the previous paragraph to the left-most square gives a monad morphism $(\ca GV,\Gamma{E'}){\rightarrow}(\ca GW,\Gamma{F'})$, and then applying $\overline{(-)}$ to this gives the lax monoidal functor
\[ (\psi_1^*,\psi') : (V^{E_1},E'){\rightarrow}(W^{F_1},F') \]
between the induced lifted multitensors. Arguing similarly for monoidal transformations and monad 2-cells, one finds that the assignment $(V,E) \mapsto (V^{E_1},E')$ is 2-functorial.

Let $\varepsilon:E{\rightarrow}\ca T^{\times}_{{\leq}n}$ be an $n$-multitensor. In terms of the previous paragraph, this is the special case $V=W=\ca G^n\Set$, $H=\id$, $\psi{=}\varepsilon$. The lifted multitensor corresponding to $\ca T^{\times}_{{\leq}n}$ is just cartesian product for strict $n$-categories. One has a component of $\varepsilon'$ for each sequence $((X_1,x_1),...,(X_n,x_n))$ of strict $n$-categories, and since $\varepsilon_1^*:\Enrich n{\rightarrow}\Alg {E_1}$ as a right adjoint preserves products, this component may be regarded as a map
\begin{equation}\label{eq:mult->prod} \begin{array}{c} {\varepsilon'_{(X_i,x_i)} : {\opE\limits_i}' \varepsilon_1^*(X_i,x_i) \rightarrow \prod\limits_i \varepsilon_1^*(X_i,x_i)} \end{array} \end{equation}
of $E_1$-algebras. If in particular $E_1$ is itself $\ca T_{{\leq}n}$ and $\varepsilon_1{=}\id$, then these components of $\varepsilon'$ give a canonical comparison from the lifted tensor product $E'$ of $n$-categories to the cartesian product. For instance, when $E$ is the multitensor corresponding to the 3-operad for Gray categories, then $\varepsilon'$ gives the well-known comparison map from the Gray tensor product of 2-categories to the cartesian product, which we recall is actually a componentwise biequivalence.

Returning to the general situation, it is routine to unpack the assignment $(H,\psi) \mapsto (\psi_1^*,\psi')$ as in section(\ref{ssec:explicit-lifting}) and so obtain the following 1-cell counterpart of corollary(\ref{cor:lifted-obj-simple}).
\begin{corollary}\label{cor:free-lift-1cell}
Let $(H,\psi):(V,E){\rightarrow}(W,F)$ be a lax monoidal functor such that $V$ and $W$ are locally presentable, and $E$ and $F$ are accessible. Let $(X_1,...,X_n)$ be a sequence of objects of $V$. Then the component of $\psi'$ at the sequence \[ (E_1X_1,...,E_nX_n) \] of free $E_1$-algebras is just $\psi_{X_i}$.
\end{corollary}
%

\section{Contractibility}\label{sec:contractibility}

\subsection{Trivial Fibrations}
Let $V$ be a category and $\ca I$ a class of maps in $V$. Denote by $\ca I^{\uparrow}$ the class of maps in $V$ that have the right lifting property with respect to all the maps in $\ca I$. That is to say, $f:X{\rightarrow}Y$ is in $\ca I^{\uparrow}$ iff for every $i:S{\rightarrow}B$ in $\ca I$, $\alpha$ and $\beta$ such that the outside of
\[ \xygraph{{S}="tl" [r] {X}="tr" [d] {Y}="br" [l] {B}="bl" "tl" (:"tr"^-{\alpha}:"br"^{f},:"bl"_{i}(:"br"_-{\beta},:@{.>}"tr"|{\gamma}))} \]
commutes, then there is a $\gamma$ as indicated such that $f\gamma{=}\beta$ and $\gamma{i}=\alpha$. An $f \in \ca I^{\uparrow}$ is called a \emph{trivial $\ca I$-fibration}. The basic facts about $\ca I^{\uparrow}$ that we shall use are summarised in
\begin{lemma}\label{lem:basic-tf}
Let $V$ be a category, $\ca I$ a class of maps in $V$, $J$ a set and \[ (f_j:X_j{\rightarrow}Y_j \,\,\, | \,\,\, j \in J) \]
a family of maps in $V$.
\begin{enumerate}
\item  $\ca I^{\uparrow}$ is closed under composition and retracts.\label{tfib1}
\item  If $V$ has products and each of the $f_j$ is a trivial $\ca I$-fibration, then
\[ \begin{array}{c} {\prod\limits_{j} f_j : \prod\limits_j X_j \rightarrow \prod\limits_j Y_j} \end{array} \]
is also a trivial $\ca I$-fibration.\label{tfib2}
\item  The pullback of a trivial $\ca I$-fibration along any map is a trivial $\ca I$-fibration.\label{tfib3}
\item  If $V$ is extensive and $\coprod_jf_j$ is a trivial $\ca I$-fibration, then each of the $f_j$ is a trivial fibration.\label{tfib4}
\item  If $V$ is extensive, the codomains of maps in $\ca I$ are connected and each of the $f_j$ is a trivial $\ca I$-fibration, then $\coprod_jf_j$ is a trivial $\ca I$-fibration.\label{tfib5}
\end{enumerate}
\end{lemma}
\begin{proof}
(\ref{tfib1})-(\ref{tfib3}) is standard. If $V$ is extensive then the squares
\[ \xygraph{!{0;(1.5,0):(0,.666)::}
{X_j}="tl" [r] {\coprod_jX_j}="tr" [d] {\coprod_jY_j}="br" [l] {Y_j}="bl" "tl" (:"tr":"br"^{\coprod_jf_j},:"bl"_{f_j}:"br")} \]
whose horizontal arrows are the coproduct injections are pullbacks, and so (\ref{tfib4}) follows by the pullback stability of trivial $\ca I$-fibrations. As for (\ref{tfib5}) note that for $i:S{\rightarrow}B$ in $\ca I$, the connectedness of $B$ ensures that any square as indicated on the left
\[ \xygraph{{\xybox{\xygraph{!{0;(1.5,0):(0,.666)::} {S}="tl" [r] {\coprod_jX_j}="tr" [d] {\coprod_jY_j}="br" [l] {B}="bl" "tl" (:"tr":"br"^{\coprod_jf_j},:"bl"_{i}:"br")}}}
[r(4)]
{\xybox{\xygraph{!{0;(1.5,0):(0,.666)::} {S}="tl" [r] {X_j}="tr" [d] {Y_j}="br" [l] {B}="bl" "tl" (:"tr":"br"^{f_j},:"bl"_{i}:"br")}}}} \]
factors through a unique component as indicated on the right, enabling one to induce the desired filler.
\end{proof}
\begin{definition}
Let $F,G:W{\rightarrow}V$ be functors and $\ca I$ be a class of maps in $V$. A natural transformation $\phi:F{\implies}G$ is a \emph{trivial $\ca I$-fibration} when its components are trivial $\ca I$-fibrations.
\end{definition}
\noindent Note that since trivial $\ca I$-fibrations in $V$ are pullback stable, this reduces, in the case where $W$ has a terminal object $1$ and $\phi$ is cartesian, to the map $\phi_1:F1{\rightarrow}G1$ being a trivial $\ca I$-fibration.

Given a category $V$ with an initial object, and a class of maps $\ca I$ in $V$, we denote by $\ca I^+$ the class of maps in $\ca GV$ containing the maps{\footnotemark{\footnotetext{Recall that $0$ is the $V$-graph with one object whose only hom is initial, or in other words the representing object of the functor $\ca GV{\rightarrow}\Set$ which sends a $V$-graph to its set of objects.}}}
\[ \begin{array}{lccr} {\emptyset \rightarrow 0} &&& {(i) : (S) \rightarrow (B)} \end{array} \]
where $i \in \ca I$. The proof of the following lemma is trivial.
\begin{lemma}\label{lem:ind-tf}
Let $V$ be a category with an initial object and $\ca I$ a class of maps in $V$. Then $f:X{\rightarrow}Y$ is a trivial $\ca I^+$-fibration iff it is surjective on objects and all its hom maps are trivial $\ca I$-fibrations.
\end{lemma}
\noindent In particular starting with $V=\PSh {\G}$ the category of globular sets and $\ca I_{-1}$ the empty class of maps, one generates a sequence of classes of maps $\ca I_n$ of globular sets by induction on $n$ by the formula $\ca I_{n+1}=(\ca I_{n})^+$ since $\ca G(\PSh {\G})$ may be identified with $\PSh {\G}$, and moreover one has inclusions $\ca I_n \subset \ca I_{n+1}$. More explicitly, the set $\ca I_n$ consists of $(n+1)$ maps: for $0{\leq}k{\leq}n$ one has the inclusion $\partial{k} \hookrightarrow k$, where $k$ here denotes the representable globular set, that is the ``$k$-globe'',
and $\partial{k}$ is the $k$-globe with its unique $k$-cell removed. One defines $\ca I_{{\leq}\infty}$ to be the union of the $\ca I_n$'s.
Note that by definition $\ca I_{{\leq}\infty}=\ca I_{{\leq}\infty}^{+}$.

There is another version of the induction just described to produce, for each $n \in \N$, a class $\ca I_{{\leq}n}$ of maps of $\ca G^n(\Set)$. The set $\ca I_{{\leq}0}$ consists of the functions
\[ \begin{array}{lccr} {\emptyset \rightarrow 0} &&& {0+0 \rightarrow 0,} \end{array} \]
so $\ca I_{\leq 0}^{\uparrow}$ is the class of bijective functions. For $n \in \N$, $\ca I_{{\leq}n{+}1}=\ca I_{{\leq}n}^{+}$. As maps of globular sets, the class $\ca I_{\leq n}$ consists of all the maps of $\ca I_{n}$ together with the unique map $\partial{(n{+}1)}{\rightarrow}n$.
\begin{definition}
Let $0{\leq}n{\leq}\infty$. An $n$-operad{\footnotemark{\footnotetext{The monad $\ca T_{{\leq}\infty}$ on globular sets is usually just denoted as $\ca T$: it is the monad whose algebras are strict $\omega$-categories.}}} $\alpha:A{\rightarrow}\ca T_{{\leq}n}$ is \emph{contractible} when it is a trivial $\ca I_{{\leq}n}$-fibration. An $n$-multitensor $\varepsilon:E{\rightarrow}\ca T_{{\leq}n}^{\times}$ is \emph{contractible} when it is a trivial $\ca I_{{\leq}n}$-fibration.
\end{definition}
\noindent By the preceeding two lemmas, an $(n+1)$-operad $\alpha:A{\rightarrow}\ca T_{{\leq}n{+}1}$ over $\Set$ is contractible iff the hom maps of $\alpha_1$ are trivial $\ca I_{\leq n}$-fibrations.

\subsection{Contractible operads versus contractible multitensors}
As one would expect an $(n+1)$-operad over $\Set$ is contractible iff its associated $n$-multitensor is contractible. This fact has quite a general explanation. Recall the 2-functoriality of $\Gamma$ described in section(\ref{ssec:Gamma-2-functor}) and that of the lifting described in section(\ref{ssec:functoriality-lifting}).
\begin{proposition}\label{prop:contractible}
Let $(H,\psi):(V,E){\rightarrow}(W,F)$ be a lax monoidal functor between distributive lax monoidal categories, and $\ca I$ a class of maps in $W$. Suppose that $W$ is extensive, $H$ preserves coproducts and the codomains of maps in $\ca I$ are connected. Then the following statements are equivalent
\begin{itemize}
\item[(1)] $\psi$ is a trivial $\ca I$-fibration.
\item[(2)] $\Gamma\psi$ is a trivial $\ca I^+$-fibration.
\end{itemize}
and moreover when in addition $V$ and $W$ are locally presentable and $E$ and $F$ are accessible, these conditions are also equivalent to
\begin{itemize}
\item[(3)] The components of $U^F\psi'$ at sequences $(E_1X_1,...,E_1X_n)$ of free $E_1$-algebras are trivial $\ca I$-fibrations.
\end{itemize}
\end{proposition}
\begin{proof}
For each $X \in \ca GV$ the component $\{\Gamma\psi\}_X$ is the identity on objects and for $a,b \in X_0$, the corresponding hom map is obtained as the composite of
\[ \begin{array}{c} {\coprod\limits_{a{=}x_0,...,x_n{=}b} \psi : \coprod\limits_{x_0,...,x_n} \opF\limits_iHX(x_{i-1}x_i) \rightarrow \coprod\limits_{x_0,...,x_n} H\opE\limits_iX(x_{i-1}x_i)} \end{array} \]
and the canonical isomorphism that witnesses the fact that $H$ preserves coproducts. In particular note that for any sequence $(Z_1,...,Z_n)$ of objects of $V$, regarded as $V$-graph in the usual way, one has \[ \{\Gamma\psi\}_{(Z_1,...,Z_n)} = \psi_{Z_1,...,Z_n}.\]
Thus $(1){\iff}(2)$ follows from lemmas(\ref{lem:basic-tf}) and (\ref{lem:ind-tf}). $(2){\iff}(3)$ follows immediately from corollary(\ref{cor:free-lift-1cell}).
\end{proof}
\begin{corollary}\label{cor:contractible}
Let $0{\leq}n{\leq}\infty$, $\alpha:A{\rightarrow}\ca T_{{\leq}n{+}1}$ be an $n{+}1$-operad over $\Set$ and $\varepsilon:E{\rightarrow}\ca T^{\times}_{{\leq}n}$ be the corresponding $n$-multitensor. TFSAE:
\begin{enumerate}
\item  $\alpha:A{\rightarrow}\ca T_{{\leq}n{+}1}$ is contractible.
\item  $\varepsilon:E{\rightarrow}\ca T^{\times}_{{\leq}n}$ is contractible.
\item  The components of $\varepsilon'_{(X_i,x_i)}$ of section(\ref{ssec:functoriality-lifting})(\ref{eq:mult->prod}) are trivial $\ca I_{{\leq}n}$-fibrations of $n$-globular sets, when the $(X_i,x_i)$ are free strict $n$-categories.
\end{enumerate}
\end{corollary}
\begin{proof}
By induction one may easily establish that the codomains of the maps in any of the classes: $\ca I_n$, $\ca I_{{\leq}n}$, $\ca I_{{\leq}\infty}$ are connected so that proposition(\ref{prop:contractible}) may be applied.
\end{proof}
\begin{example}\label{ex:Gray-contractible}
Applying this last result to the 3-operad $G$ for Gray-categories, the contractibility of $G$ is a consequence of the fact that the canonical 2-functors from the Gray to the cartesian tensor product are identity-on-object biequivalences.
\end{example}

\subsection{Trimble \`{a} la Cheng II}\label{ssec:TCII}
Continuing the discussion from section(\ref{ssec:TCI}), we now explain why the operads which describe ``Trimble $n$-categories" are contractible. This result appears in \cite{ChengCompOp} as theorem(4.8) and exhibits Trimble $n$-categories as weak $n$-categories in the Batanin sense.

Let us denote by $\ca J$ the set of inclusions $S^{n{-}1}{\rightarrow}D^n$ of the $n$-sphere into the $n$ disk for $n \in \N$. As we remarked in section(\ref{ssec:enriched-graphs}) these may all be obtained by successively applying the reduced suspension functor $\sigma$ to the inclusion of the empty space into the point. As we recalled in section(\ref{ssec:TCI}), a basic ingredient of the Trimble definition is a version of the little intervals operad which acts on the path spaces of any space. A key property of this operad is that it is contractible -- a topological operad $A$ being contractible when for each $n$ the unique map $A_n \rightarrow 1$ is in $\ca J^{\uparrow}$. This is equivalent to saying that the cartesian multitensor map $A \rightarrow \prod$ is a trivial $\ca J$-fibration. A useful fact about the class trivial $\ca J$-fibrations is that it gets along with the construction of path-spaces in the sense of
\begin{lemma}\label{lem:tf-path-spaces}
If $f:X{\rightarrow}Y$ is a trivial $\ca J$-fibration then so is $f_{a,b}:X(a,b){\rightarrow}Y(fa,fb)$ for all $a,b \in X$.
\end{lemma}
\begin{proof}
To give a commutative square as on the left in
\[ \begin{array}{lccr} {\xybox{\xygraph{!{0;(2,0):(0,.5)::} {S^{n{-}1}}="tl" [r] {X(a,b)}="tr" [d] {Y(fa,fb)}="br" [l] {D^n}="bl" "tl"(:"tr":"br",:"bl":"br")}}} &&&
{\xybox{\xygraph{!{0;(2,0):(0,.5)::} {S^n}="tl" [r] {(a,X,b)}="tr" [d] {(fa,Y,fb)}="br" [l] {D^{n{+}1}}="bl" "tl"(:"tr":"br",:"bl":"br")}}} \end{array} \]
is the same as giving a commutative square in $\Top_{\bullet}$ as on the right in the previous display, by $\sigma \ladj h$. The square on the right admits a diagonal filler $D^{n{+}1} \rightarrow X$ since $f$ is a trivial $\ca J$-fibration, and thus so does the square on the left.
\end{proof}
We shall write $U_n:\Trimble n \rightarrow \ca G^n\Set$ for the forgetful functor for each $n$. The relationship between trivial fibrations of spaces and of globular sets is expressed in
\begin{proposition}\label{prop:pres-tf}
If $f:X{\rightarrow}Y$ is a trivial $\ca J$-fibration then $U_n\pi_n{f}$ is a trivial $\ca I_{\leq{n}}$-fibration.
\end{proposition}
\begin{proof}
We proceed by induction on $n$. Having the right lifting property with respect to the inclusions
\[ \begin{array}{lccr} {\emptyset \hookrightarrow 1} &&& {1{+}1=\partial{I} \hookrightarrow I} \end{array} \]
ensures that $f$ surjective and injective on path components, and thus is inverted by $\pi_0$. For the inductive step we assume that $U_n\pi_n$ sends trivial $\ca J$-fibrations to trivial $\ca I_{{\leq}n}$-fibrations and suppose that $f$ is a trivial $\ca J$-fibration. Then so are all the maps it induces between path spaces by lemma(\ref{lem:tf-path-spaces}). But from the inductive definition of $\Trimble {n{+}1}$ recalled in section(\ref{ssec:TCI}), we have $U_{n{+}1}\pi_{n{+}1}=\ca G(u_n\pi_n)P$ and so $U_{n{+}1}\pi_{n{+}1}(f)$ is a morphism of $(n{+}1)$-globular sets which is surjective on objects (as argued already in the $n=0$ case) and whose hom maps are trivial $\ca I_{{\leq}n}$-fibrations by induction. Thus the result follows by lemma(\ref{lem:ind-tf}).
\end{proof}
In section(\ref{ssec:TCI}) we exhibited $\Trimble n$ as the algebras of an $n$-operad by a straight-forward application of two abstract results -- propositions(\ref{prop:monad-trimble}) and (\ref{prop:operad-trimble}). We now provide a third such result relating to contractibility.
\begin{proposition}\label{prop:tf-formal}
Given the data and hypotheses of proposition(\ref{prop:operad-trimble}): $V$ is a lextensive category, $T$ a cartesian and coproduct preserving monad on $V$, $\psi:S{\rightarrow}T$ a $T$-operad and $(E,\varepsilon)$ a non-symmetric operad in $V^S$. Suppose furthermore that a class $\ca I$ of maps of $V$ is given, and that the non-symmetric operad $\alpha:E \rightarrow \prod$ and the $T$-operad $\psi$ are trivial $\ca I$ fibrations. Then the $\Gamma(T^{\times})$-operad
\[ \Gamma(E)\ca G(S) \rightarrow \Gamma(T^{\times}) \]
of proposition(\ref{prop:operad-trimble}) is a trivial $\ca I^{+}$-fibration.
\end{proposition}
\begin{proof}
By definition this monad morphism may be written as the composite $(\Gamma(\psi^{\times}))(\Gamma(\alpha)\ca G(S))$. Since $\psi$ is a trivial $\ca I$-fibration so is $\psi^{\times}$ by lemma(\ref{lem:basic-tf}), and thus $\Gamma(\psi^{\times})$ is a trivial $\ca I^{+}$-fibration by proposition(\ref{prop:contractible}). Since $\alpha$ is a trivial $\ca I$-fibration, $\Gamma(\alpha)$ is a trivial $\ca I^{+}$-fibration again by proposition(\ref{prop:contractible}), and so the result follows since trivial fibrations compose.
\end{proof}
Starting with a contractible topological operad $A$ which acts on path spaces, proposition(\ref{prop:pres-tf}) ensures that $U_n\pi_nA$ will be a contractible non-symmetric operad of $n$-globular sets. Then proposition(\ref{prop:tf-formal}) may be applied to give, by induction on $n$, the contractibility of the $n$-operad defining Trimble $n$-categories.

\section{Acknowledgements}\label{sec:Acknowledgements}
We would like to acknowledge Clemens Berger, Richard Garner, Andr\'{e} Joyal, Steve Lack, Joachim Kock, Jean-Louis Loday, Paul-Andr\'{e} Melli\`{e}s, Ross Street and Dima Tamarkin for interesting discussions on the substance of this paper. We would also like to acknowledge the Centre de Recerca Matem\`{a}tica in Barcelona for the generous hospitality and stimulating environment provided during the thematic year 2007-2008 on Homotopy Structures in Geometry and Algebra. The first author would like to acknowledge the financial support on different stages of this project of the Scott Russell Johnson Memorial Foundation, the Australian Research Council (grant No.DP0558372) and L'Universit\'{e} Paris 13. The second author would like to acknowledge the support of the ANR grant no. ANR-07-BLAN-0142.
The third author would like to acknowledge the laboratory PPS (Preuves Programmes Syst\`{e}mes) in Paris, the Max Planck Institute in Bonn and the Institut des Hautes \'Etudes Scientifique in Bures sur Yvette for the excellent working conditions he enjoyed during this project, as well as Macquarie University for the hospitality he enjoyed in August of 2008.

\appendix

\section{Coequalisers in categories of algebras}\label{sec:coequalisers}
In these appendices we review some of the transfinite constructions in monad theory that we used in section(\ref{sec:lift-mult}). An earlier reference for these matters is \cite{Kel80}. However due to the technical nature of this material, and our need for its details when we come to making our constructions explicit, we feel that it is appropriate to give a rather thorough account of this background.

Let $T$ be a monad on a category $V$ that has filtered colimits and coequalisers and let
\[ \xygraph{!{0;(2,0):} {(A,a)}="l" [r] {(B,b)}="r" "l":@<-1ex>"r"_-{g} "l":@<1ex>"r"^-{f}} \]
be morphisms in $V^T$. We shall now construct morphisms
\[ \begin{array}{lcccr} {v_n : TQ_n \rightarrow Q_{n{+}1}} && {q_n:Q_n \to Q_{n{+}1}} && {q_{{<}n}:B \to Q_n} \end{array} \]
starting with $Q_0=B$ by transfinite induction on $n$, such that for $n$ large enough $q_{<{n}}$ is the coequaliser of $f$ and $g$ in $V^T$ when $T$ is accessible. The initial stages of this construction are described in the following diagram.
\[ \xygraph{!{0;(1.5,0):(0,.666)::} {T^2A}="tA" [r] {T^2B}="tB" [r] {T^2Q_1}="t0" [r] {T^2Q_2}="t1" [r] {T^2Q_3}="t2" [r] {T^2Q_4}="t3" [r] {...}="t4"
"tA" [d] {TA}="mA" [r] {TB}="mB" [r] {TQ_1}="m0" [r] {TQ_2}="m1" [r] {TQ_3}="m2" [r] {TQ_4}="m3" [r] {...}="m4"
"mA" [d] {A}="bA" [r] {B}="bB" [r] {Q_1}="b0" [r] {Q_2}="b1" [r] {Q_3}="b2" [r] {Q_4}="b3" [r] {...}="b4"
"tA":@<-1ex>"tB"_-{T^2g} "tA":@<1ex>"tB"^-{T^2f}:"t0"^-{T^2q_0}:"t1"^-{T^2q_1}:"t2"^-{T^2q_2}:"t3"^-{T^2q_3}:"t4"
"mA":@<-1ex>"mB"_-{Tg} "mA":@<1ex>"mB"^-{Tf}:"m0"^-{Tq_0}:"m1"^-{Tq_1}:"m2"^-{Tq_2}:"m3"^-{Tq_3}:"m4"
"bA":@<-1ex>"bB"_-{g} "bA":@<1ex>"bB"^-{f}:"b0"_-{q_0}:"b1"_-{q_1}:"b2"_-{q_2}:"b3"_-{q_3}:"b4"
"tA":@<-1ex>"mA"_-{\mu}:@<-1ex>"bA"_-{a} "tA":@<1ex>"mA"^-{Ta}:@<1ex>@{<-}"bA"^-{\eta}
"tB":@<-1ex>"mB"_-{\mu}:@<-1ex>"bB"_-{b} "tB":@<1ex>"mB"^-{Tb}:@<1ex>@{<-}"bB"^-{\eta}
"t0":"m0"^-{\mu}:@{<-}"b0"^-{\eta} "t1":"m1"^-{\mu}:@{<-}"b1"^-{\eta} "t2":"m2"^-{\mu}:@{<-}"b2"^-{\eta} "t3":"m3"^-{\mu}:@{<-}"b3"^-{\eta}
"tB":"m0"^-{Tv_0} "mB":"b0"^-{v_0} "t0":"m1"^-{Tv_1} "m0":"b1"^-{v_1} "t1":"m2"^-{Tv_2} "m1":"b2"^-{v_2} "t2":"m3"^-{Tv_3} "m2":"b3"^-{v_3}} \]
{\bf Initial step}. Define $q_{{<}0}$ to be the identity, $q_0$ to be the coequaliser of $f$ and $g$, $q_{<{1}}=q_0$ and $v_0=q_0b$. Note also that $q_0=v_0\eta_B$.
\\ \\
{\bf Inductive step}. Assuming that $v_n$, $q_n$ and $q_{<{n{+}1}}$ are given, we define $v_{n{+}1}$ to be the coequaliser of $T(q_n)\mu$ and $Tv_n$, $q_{n{+}1}=v_{n{+}1}\eta$ and $q_{<{n{+}2}}=q_{n{+}1}q_{<{n{+}1}}$. One may easily verify that $q_{n{+}1}v_n=v_{n{+}1}T(q_n)$, and that $v_1$ could equally well have been defined as the coequaliser of $\eta{v_0}$ and $Tq_0$.
\\ \\
{\bf Limit step}. Define $Q_n$ as the colimit of the sequence given by the objects $Q_m$ and morphisms $q_m$ for $m < n$, and $q_{<{n}}$ for the component of the universal cocone at $m=0$.
\[ \xygraph{!{0;(3,0):(0,.333)::} {\colim_{m{<}n} T^2Q_m}="tl" [r] {\colim_{m{<}n} TQ_m}="tm" [r] {\colim_{m{<}n} Q_m}="tr" [d] {Q_n}="br" [l] {TQ_n}="bm" [l] {T^2Q_n}="bl" "tl":@<1ex>"tm"^-{\mu_{<{n}}}:@<1ex>"tr"^-{v_{<{n}}} "tl":@<-1ex>"tm"_-{(Tv)_{<{n}}}:@<-1ex>@{<-}"tr"_-{\eta_{<{n}}} "bl":"bm"_-{\mu}:@{<-}"br"_-{\eta} "tl":"bl"_{o_{n,2}} "tm":"bm"^{o_{n,1}} "tr":@{=}"br"} \]
We write $o_{n,1}$ and $o_{n,2}$ for the obstruction maps measuring the extent to which $T$ and $T^2$ preserve the colimit defining $Q_n$. We write $\mu_{<{n}}$, $(Tv)_{<{n}}$, $v_{<{n}}$ and $\eta_{<{n}}$ for the maps induced by the $\mu_{Q_m}$, $Tv_m$, $v_m$ and $\eta_{Q_m}$ for $m < n$ respectively. The equations
\[ \begin{array}{ccccccc} {{\mu}o_{n,2}=o_{n,1}\mu_{<{n}}} && {\eta=o_{n,1}\eta_{<{n}}} && {v_{<{n}}(Tv)_{<{n}}=v_{<{n}}\mu_{<{n}}} && {v_{<{n}}\eta_{<{n}}=\id} \end{array} \]
follow easily from the definitions. Define $v_n$ as the coequaliser of $o_{n,1}\mu_{<{n}}$ and $o_{n,1}(Tv)_{<{n}}$, $q_n=v_n\eta$ and $q_{<{n{+}1}}=q_nq_{<{n}}$.
\\ \\
{\bf Stabilisation}. We say that the sequence \emph{stabilises at $n$} when $q_n$ and $q_{n{+}1}$ are isomorphisms. In the case $n=0$ one may easily show that stabilisation is equivalent to just $q_0$ being an isomorphism, which is the same as saying that $f=g$.
\begin{lemma}\label{stable-limit}
If $n$ is a limit ordinal and $o_{n,1}$ and $o_{n,2}$ are invertible, then the sequence stabilises at $n$.
\end{lemma}
\begin{proof}
Let us write $q_{m,n}:Q_m \to Q_n$,
\[ \begin{array}{lccr} {q'_{m,n}:TQ_m \to \colim_{m<n}TQ_m} &&& {q''_{m,n}:T^2Q_m \to \colim_{m<n}T^2Q_m} \end{array} \]
for the colimit cocones. First we contemplate the diagram
\[ \xygraph{!{0;(3,0):(0,.333)::} {\colim_{m{<}n} T^2Q_m}="tl" [r] {\colim_{m{<}n} TQ_m}="tm" [r] {\colim_{m{<}n} Q_m}="tr" [d] {Q_n}="br" [l] {TQ_n}="bm" [l] {T^2Q_n}="bl" [d] {T^2Q_{n{+}1}}="bbl" [r] {TQ_{n{+}1}}="bbm" [r] {Q_{n{+}1}}="bbr" [d] {Q_{n{+}2}}="bbbr"
"tl":@<1ex>"tm"^-{\mu_{<{n}}}:@<1ex>"tr"^-{v_{<{n}}} "tl":@<-1ex>"tm"_-{(Tv)_{<{n}}}:@<-1ex>@{<-}"tr"_-{\eta_{<{n}}} "bl":"bm"_-{\mu}:@{<-}"br"_-{\eta} "bbl":"bbm"_-{\mu}:@{<-}"bbr"_-{\eta} "tl":"bl"_{o_{n,2}} "tm":"bm"^{o_{n,1}} "tr":@{=}"br" "bl":"bbl"_{T^2q_n} "bm":"bbm"^{Tq_n} "br":"bbr"^{q_n}:"bbbr"^{q_{n{+}1}} "bl":"bbm"_-{Tv_n}:"bbbr"_-{v_{n{+}1}} "bm":"bbr"_-{v_n}} \]
and in general one has
\begin{equation}\label{eq:succ-limit} T(q_n)o_{n,1}(Tv)_{{<}n} = T(v_n)o_{n,2}. \end{equation}
To prove this note that from the definitions of $q_m$ and $q_n$ and the naturality of the $q_{m,n}$ in $m$, one may show easily that $v_nT(q_{m,n})=q_nq_{m{+}1,n}v_m$, and from this last equation and all the definitions it is easy to show that
\[ T(q_n)o_{n,1}(Tv)_{{<}n}q''_{m,n} = T(v_n)o_{n,2}q''_{m,n} \]
for all $m < n$ from which (\ref{eq:succ-limit}) follows.

Suppose that $o_{n,1}$ and $o_{n,2}$ are isomorphisms. Then define $q'_{n}:Q_{n{+}1} \to Q_n$ as the unique map such that $q'_nv_no_{n,1}=q_nv_{{<}n}$. It follows easily that $q'_n=q^{-1}_n$. From (\ref{eq:succ-limit}) and the invertibility of $o_{n,2}$ it follows easily that $v_n\mu=v_nT(q_n^{-1})T(v_n)$ and so there is a unique $q'_{n{+}1}$ such that $q'_{n{+}1}v_{n{+}2}=v_nT(q_n^{-1})$, from which it follows easily that $q'_{n{+}1}=q_{n{+}1}^{-1}$.
\end{proof}
\begin{lemma}\label{really-stable}
If the sequence stabilises at $n$ then it stabilises at any $m \geq n$, and moreover one has an isomorphism of sequences between the given sequence $(Q_m,q_m)$ and the following one:
\[ \xygraph{{Q_0}="p1" [r] {...}="p2" [r] {Q_n}="p3" [r] {Q_n}="p4" [r] {...}="p5" "p1":"p2"^-{q_0}:"p3":"p4"^-{\id}:"p5"^-{\id}} \]
\end{lemma}
\begin{proof}
We show for $m \geq n$ that $q_m$ and $q_{m{+}1}$ are isomorphisms, and provide the component isomorphisms $i_m:Q_m \to Q_n$ of the required isomorphism of sequences, by transfinite induction on $m$. We define $i_m$ to be the identity when $m \leq n$. In the initial step $m=n$, $q_m$ and $q_{m{+}1}$ are isomorphisms by hypothesis and we define $i_{n{+}1}=q_n$. In the inductive step when $m \geq n$ is a non-limit ordinal, we must show that $q_{m{+}2}$ is an isomorphism and define $i_{m{+}2}=q_{m{+}1}i_{m{+}1}$. The key point is that
\begin{equation}\label{eq:key} v_{m{+}1}\mu = v_{m{+}1}T(q^{-1}_{m{+}1})T(v_{m{+}1}) \end{equation}
because with this equation in hand one defines $q'_{m{+}2}:Q_{m{+}3} \to Q_{m{+}2}$ as the unique morphism satisfying $q'_{m{+}2}v_{m{+}2}T(q_{m{+}1})=v_{m{+}1}$ using the universal property of $v_{m{+}2}$, and then it is routine to verify that $q'_{m{+}2}=q^{-1}_{m{+}2}$. So for the inductive step it remains to verify (\ref{eq:key}). But we have
\[ \begin{array}{rll} v_{m{+}1}{\mu}T^2(q_m) &=& v_{m{+}1}T(q_m)\mu = v_{m{+}1}T(v_m) = v_{m{+}1}T(q^{-1}_{m{+}1})T(q_{m{+}1}v_m) \\
&=& v_{m{+}1}T(q^{-1}_{m{+}1})T(v_{m{+}1})T^2(q_m) \end{array} \]
and so (\ref{eq:key}) follows since $q_m$ is an isomorphism. In the case where $m$ is a limit ordinal, we have stabilisation at $m'$ established whenever $n \leq m' < m$ by the induction hypothesis. Thus the colimit defining $Q_m$ is absolute (ie preserved by all functors) since its defining sequence from the position $n$ onwards consists only of isomorphisms. Thus $q_m$ and $q_{m{+}1}$ are isomorphisms by lemma(\ref{stable-limit}). By induction, the previously constructed $i_{m'}$'s provide a cocone on the defining diagram of $Q_m$ with vertex $Q_n$, thus one induces the isomorphism $i_m$ compatible with the earlier $i_{m'}$'s and defines $i_{m{+}1}=q_mi_m$.
\end{proof}
\begin{lemma}\label{coeq-when-stable}
If the sequence stabilises at $n$ then $(Q_n,q_n^{-1}v_n)$ is a $T$-algebra and
\[ q_{<{n}}:(B,b) \to (Q_n,q_n^{-1}v_n) \]
is the coequaliser of $f$ and $g$ in $V^T$.
\end{lemma}
\begin{proof}
The unit law for $(Q_n,q_n^{-1}v_n)$ is immediate from the definition of $q_n$ and the associative law is the commutativity of the outside of the diagram on the left
\[ \xygraph{{\xybox{\xygraph{!{0;(1.5,0):(0,.7)::} {T^2Q_n}="tl" [r(2)] {TQ_n}="tr" [d] {Q_{n{+}1}}="mr" [d] {Q_n}="br" [l] {Q_{n{+}1}}="bm" [l] {TQ_n}="bl" [u] {TQ_{n{+}1}}="ml" "bm" [u(.75)] {Q_{n{+}2}}="mm" [u(.75)] {TQ_{n{+}1}}="tm"
"tl":"tr"^-{\mu}:"mr"^-{v_n}:"br"^-{q_n^{-1}}:@{<-}"bm"^-{q_n^{-1}}:@{<-}"bl"^-{v_n}:@{<-}"ml"^-{Tq_n^{-1}}:@{<-}"tl"^-{Tv_n}
"tr":"tm"^-{Tq_n}:"mm"_-{v_{n{+}1}}:"bm"_-{q_{n{+}1}^{-1}} "ml":"mm"_-{v_{n{+}1}}:"mr"^-{q_{n{+}1}^{-1}}}}}
[r(5)u(.1)]
{\xybox{\xygraph{{TB}="tl" [r(2)] {TQ_n}="tr" [d] {Q_{n{+}1}}="m" [d] {Q_n}="br" [l(2)] {B}="bl" "tl":"tr"^-{Tq_{{<}n}}:"m"^-{v_n}:"br"^-{q_n^{-1}}:@{<-}"bl"^-{q_{{<}n}}(:@{<-}"tl"^-{b},:"m"^-{q_{{<}n{+}1}}) "m":@{}"tl"|{(I)})}}}} \]
the regions of which evidently commute. The commutativity of the outside diagram on the right exhibits $q_{{<}n}$ as a $T$-algebra map, and this follows immediately from the commutativity of the region labelled (I).

The equational form of (I) says $q_{{<}n{+}1}b=v_nT(q_{{<}n})$ and we now proceed to prove this by transfinite induction on $n$. The case $n=0$ is just the statement $v_0=q_0b$. The inductive step comes out of the calculation
\[ q_{{<}n{+}2}b = q_{n{+}1}q_{{<}n{+}1}b = q_{n{+}1}v_nT(q_{{<}n}) = v_{n{+}1}T(q_nq_{{<}n}) = v_{n{+}1}T(q_{{<}n{+}1}) \]
and since $Tq_{{<}n}=o_{n,1}q'_{0,n}$. The case where $n$ is a limit ordinal is the commutativity of the outside of
\[ \xygraph{!{0;(1.4,0):} {TB}="tb" ([ur] {B}="b" [r(3)] {Q_{n{+}1}}="qnp1",[dr] {TB}="tb2" [r(3)] {\col \, TQ_m}="ctq" [ur] {TQ_n}="tqn")
"tb"(:"b"^-{b}:"qnp1"^-{q_{{<}n{+}1}},:"tb2"_-{\id}|-{}="mla":@/_{1pc}/"ctq"_-{q'_{0,n}}|-{}="mma":"tqn"_-{o_{n,1}}|-{}="mra":"qnp1"_-{v_n} "b":@{}"ctq"|*{TB}="m")
"m" (:@{}"b"|(.6)*{Q_n}="ptl",:@{}"qnp1"|(.6)*{TQ_n}="ptr",:@{}"mra"|(.6)*{\col \, TQ_m}="pmr",:@{}"mma"|(.6)*{\col \, T^2Q_m}="pb",:@{}"mla"|(.6)*{T^2B}="pml")
"b":"ptl"_(.6){q_{{<}n}}:"ptr"^-{\eta}:"qnp1"^(.4){v_n}
"tb":"pml"^-{\eta}:"m"^-{Tb}(:"ptr"^(.4){Tq_{{<}n}},:"pmr"^-{q'_{0,n}})
"pml":"pb"^-{q''_{0,n}}:"pmr"^-{(Tv)_{{<}n}}:"ptr"^-{o_{n,1}}
"pml":"tb2"^-{\mu} "pb":"ctq"^-{\mu_{{<}n}}} \]
the regions of which evidently commute. Thus $q_{{<}n}$ is indeed a $T$-algebra map.

To see that $q_{{<}n}$ is a coequaliser let $h:(B,b) \to (C,c)$ such that $hf=hg$. For each ordinal $m$ we construct $h_m : Q_m \to C$ such that $h_{m{+}1}v_m=cT(h_m)$ for all $m$ by transfinite induction on $m$. When $m=0$ we define $h_0=h$ and $h_1$ as unique such that $h_1q_0=h$. The equation $h_1v_0=cT(h)$ is easily verified. For the inductive step we note that the commutativity of
\[ \xygraph{{T^2Q_n}="ml" [ur] {TQ_n}="tl" [r(2)] {TC}="tr" [dr] {C}="mr" [dl] {TC}="br" [l(2)] {TQ_{n{+}1}}="bl" [ur] {T^2C}="m"
"ml":"tl"^-{\mu}:"tr"^-{Th_n}:"mr"^-{c}:@{<-}"br"^-{c}:@{<-}"bl"^-{Th_{n{+}1}}:@{<-}"ml"^-{Tv_n} "m"(:@{<-}"ml"_-{T^2h_n},:"tr"^-{\mu},:"br"^-{Tc})} \]
and the universal property of $v_{m{+}1}$ ensures there is a unique $h_{m{+}2}$ such that $h_{m{+}2}v_{m{+}1}=cT(h_{m{+}1})$. When $m$ is a limit ordinal it follows from all the definitions that
\[ cT(h_m)o_{m,1}\mu_{{<}m}q''_{m',m} = cT(h_m)o_{m,1}(Tv)_{{<}m}q''_{m',m} \]
for all $m' < m$, and so $cT(h_m)o_{m,1}\mu_{{<}m} = cT(h_m)o_{m,1}(Tv)_{{<}m}$ and so by the universal property of $v_{m{+}1}$ there is a unique $h_{m{+}1}$ such that $h_{m{+}1}v_m=cT(h_m)$. The sequence of $h_m$'s just constructed is clearly unique such that $h_0=h$ and $h_{m{+}1}v_m=cT(h_m)$. It follows immediately that $h_n$ is a $T$-algebra map, and that $h_nq_{{<}n}=h$. Conversely given $h':Q_n \to C$ such that $h'q_{{<}n}=h$, one constructs $h'_m:Q_m \to C$ as $h'_{m}=h'q_{m,n}$, and it follows easily that $h'_0=h$, $h'_{m{+}1}v_m=cT(h'_m)$ and $h'_n=h'$ whence $h'_m=h_m$ and so $h=h'$.
\end{proof}
\noindent From these results we recover the usual theorem on the construction of coequalisers of algebras of accessible monads.
\begin{theorem}\label{thm:coeq-Talg}
Let $V$ be a category with filtered colimits and coequalisers, $T$ be a monad on $V$ and
\[ \xygraph{!{0;(2,0):} {(A,a)}="l" [r] {(B,b)}="r" "l":@<-1ex>"r"_-{g} "l":@<1ex>"r"^-{f}} \]
be morphisms in $V^T$. If $T$ is $\lambda$-accessible for some regular cardinal $\lambda$, then $q_{<{n}}$ as constructed above is the coequaliser of $f$ and $g$ in $V^T$ for any ordinal $n$ such that $|n| \geq \lambda$.
\end{theorem}
\begin{proof}
Take the smallest such ordinal $n$ -- it is necessarily a limit ordinal, and $T$ and $T^2$ by hypothesis preserve the defining colimit of $Q_n$. Thus by lemmas(\ref{stable-limit}) and (\ref{coeq-when-stable}) the result follows in this case, and in general by lemmas(\ref{really-stable}) and (\ref{coeq-when-stable}).
\end{proof}
Finally we mention the well-known special case when the above transfinite construction is particularly simple, that will be worth remembering.
\begin{proposition}\label{prop:simple-coeq}
Let $V$ be a category with filtered colimits and coequalisers, $T$ be a monad on $V$ and
\[ \xygraph{!{0;(2,0):} {(A,a)}="l" [r] {(B,b)}="r" "l":@<-1ex>"r"_-{g} "l":@<1ex>"r"^-{f}} \]
be morphisms in $V^T$. If $T$ and $T^2$ preserve the coequaliser of $f$ and $g$ in $V$, then the sequence $(Q_n,q_n)$ stabilises at $1$. Denoting by $w:TQ_1 \to Q_1$ the unique map such that $wT(q_0)=q_0b$, $q_0:(B,b) \to (Q_1,w)$ is the coequaliser of $f$ and $g$ in $V^T$.
\end{proposition}
\begin{proof}
Refer to the diagram in $V$ above that describes the first few steps of the construction of $(Q_n,q_n)$. Since $q_0$ and $T^2q_0$ are epimorphisms, the $T$-algebra axioms for $(Q_1,w)$ follow from those for $(B,b)$, and $q_0$ is a $T$-algebra map by definition. Thus $w$ is the coequaliser in $V$ of $\mu_{Q_1}$ and $Tw$, and since $T^2q_0$ is an epimorphism it is also the coequaliser of $\mu_{Q_1}T^2(q_0)$ and $T(w)T^2(q_0)=Tv_0$, but so is $v_1$, and so $q_1$ is the canonical isomorphism between them. To see that $q_2$ is also invertible, apply the same argument with the composite $q_1q_0$ in place of $q_0$. The result now follows by lemma(\ref{coeq-when-stable}).
\end{proof}
%

\section{Monads induced by monad morphisms}\label{sec:Dubuc}
Suppose that $V$ is locally presentable, $(M,\eta^M,\mu^M)$ and $(S,\eta^S,\mu^S)$ are monads on $V$, and $\phi:M{\rightarrow}S$ is a morphism of monads. Then one has the obvious forgetful functor $\phi^*:V^S{\rightarrow}V^M$ and when $S$ is accessible, $\phi^*$ has a left adjoint which we denote as $\phi_!$. The general fact responsible for the existence of $\phi_!$, and which in fact gives a formula for it in terms of coequalisers in $V^S$, is the Dubuc adjoint triangle theorem \cite{Dubuc}: for an algebra $(X,x:MX{\rightarrow}X)$ of $M$, one has the reflexive coequaliser
\[ \xygraph{!{0;(3,0):}
{(SMX,\mu^S_{MX})}="l" [r] {(SX,\mu^S_X)}="m" [r] {\phi_!(X,x)}="r"
"l":@<2ex>"m"^-{\mu^S_XS(\phi_X)}:"l"|-{S\eta^M_X}:@<-2ex>"m"_-{Sx}:"r"^-{q_{(X,x)}}} \]
in $V^S$. Putting this together with section(\ref{sec:coequalisers}) an explicit description of the composite $U^S\phi_!$ is given as follows. We construct morphisms
\[ \begin{array}{lccr} {v_{n,X,x} : SQ_n(X,x) \rightarrow Q_{n{+}1}(X,x)} &&& {q_{n,X,x}:Q_n(X,x) \to Q_{n{+}1}(X,x)} \end{array} \]
\[ \begin{array}{c} {q_{{<}n,X,x}:SX \to Q_n(X,x)}  \end{array} \]
starting with $Q_0(X,x)=SX$ by transfinite induction on $n$.
\\ \\
{\bf Initial step}. Define $q_{{<}0}$ to be the identity, $q_0$ to be the coequaliser of $\mu^S(S\phi)$ and $Sx$, $q_{<{1}}=q_0$ and $v_0=q_0b$. Note also that $q_0=v_0\eta^S$.
\\ \\
{\bf Inductive step}. Assuming that $v_n$, $q_n$ and $q_{<{n{+}1}}$ are given, we define $v_{n{+}1}$ to be the coequaliser of $S(q_n)(\mu^SQ_n)$ and $Sv_n$, $q_{n{+}1}=v_{n{+}1}(\eta^SQ_{n{+}1})$ and $q_{<{n{+}2}}=q_{n{+}1}q_{<{n{+}1}}$.
\\ \\
{\bf Limit step}. Define $Q_n(X,x)$ as the colimit of the sequence given by the objects $Q_m(X,x)$ and morphisms $q_m$ for $m < n$, and $q_{<{n}}$ for the component of the universal cocone at $m=0$.
\[ \xygraph{!{0;(3,0):(0,.333)::} {\colim_{m{<}n} S^2Q_m}="tl" [r] {\colim_{m{<}n} SQ_m}="tm" [r] {\colim_{m{<}n} Q_m}="tr" [d] {Q_n}="br" [l] {SQ_n}="bm" [l] {S^2Q_n}="bl" "tl":@<1ex>"tm"^-{\mu_{<{n}}}:@<1ex>"tr"^-{v_{<{n}}} "tl":@<-1ex>"tm"_-{(Sv)_{<{n}}}:@<-1ex>@{<-}"tr"_-{\eta_{<{n}}} "bl":"bm"_-{\mu}:@{<-}"br"_-{\eta} "tl":"bl"_{o_{n,2}} "tm":"bm"^{o_{n,1}} "tr":@{=}"br"} \]
We write $o_{n,1}$ and $o_{n,2}$ for the obstruction maps measuring the extent to which $S$ and $S^2$ preserve the colimit defining $Q_n(X,x)$. We write $\mu^S_{<{n}}$, $(Sv)_{<{n}}$, $v_{<{n}}$ and $\eta^S_{<{n}}$ for the maps induced by the $\mu^SQ_m$, $Sv_m$, $v_m$ and $\eta^SQ_m$ for $m < n$ respectively. Define $v_n$ as the coequaliser of $o_{n,1}\mu_{<{n}}$ and $o_{n,1}(Sv)_{<{n}}$, $q_n=v_n({\eta^S}Q_n)$ and $q_{<{n{+}1}}=q_nq_{<{n}}$.
\\ \\
Instantiating theorem(\ref{thm:coeq-Talg}) to the present situation gives
\begin{corollary}\label{cor:explicit-phi-shreik}
Suppose that $V$ is a locally presentable category, $M$ and $S$ are monads on $V$, $\phi:M{\rightarrow}S$ is a morphism of monads, and $(X,x)$ is an $M$-algebra. If moreover $S$ is $\lambda$-accessible for some regular cardinal $\lambda$, then for any ordinal $n$ such that $|n| \geq \lambda$
one may take
\[ \begin{array}{lccr} {\phi_!(X,x)=(Q_n(X,x),q_n^{-1}v_n)} &&& {q_{<n} : (SX,\mu_X) \to (Q_n(X,x),q_n^{-1}v_n)} \end{array} \]
as an explicit definition of $\phi_!(X,x)$ and the associated coequalising map in $V^S$ coming from the Dubuc adjoint triangle theorem.
\end{corollary}
\noindent and instantiating proposition(\ref{prop:simple-coeq}) to the present situation gives
\begin{corollary}\label{cor:phi-shreik-simple}
Suppose that under the hypotheses of corollary(\ref{cor:explicit-phi-shreik}) that $S$ and $S^2$ preserve the coequaliser of $\mu_{X}^{S}S(\phi_X)$ and $Sx$ in $V$. Then the sequence $(Q_n,q_n)$ stabilises at $1$, and writing $w:SQ_1 \to Q_1$ for the unique map such that $wS(q_0)=q_0\mu_{X}^{S}$, one may take
\[ \begin{array}{lccr} {\phi_!(X,x)=(Q_1(X,x),w)} &&& {q_0 : (SX,\mu_X) \to (Q_1(X,x),w)} \end{array} \]
as an explicit definition of $\phi_!(X,x)$ and the associated coequalising map in $V^S$.
\end{corollary}
\begin{remark}
Here is a degenerate situation in which corollary(\ref{cor:phi-shreik-simple}) applies. Since $U^M\phi^*=U^S$ we have $\phi_!F^M \iso F^U$, but another way to view this isomorphism as arising is to apply the corollary in the case where $(X,x)$ is a free $M$-algebra, say $(X,x)=(MZ,\mu^M_Z)$, for in this case one has the dotted arrows in
\[ \xygraph{!{0;(2,0):} {SM^2Z}="l" [r] {SMZ}="m" [r] {SZ}="r" "l":@<1.5ex>"m"^-{(\mu^SM)(S{\phi}M)} "l":@<-1.5ex>"m"^-{S\mu^M} "l":@{<.}@<-3.5ex>"m"_-{SM\eta^M} "m":"r"^-{\mu^SS(\phi)} "m":@{<.}@<-2ex>"r"_-{S\eta^M}} \]
exhibiting $\mu_Z^SS(\phi_Z)$ as a split coequaliser, and thus absolute.
\end{remark}
Let us denote by $(T,\eta^T,\mu^T)$ the monad on $V^M$ induced by the adjunction $\phi_! \ladj \phi^*$. While a completely explicit description of this monad is unnecessary for the proof of theorem(\ref{thm:lift-mult}), we will require such a description in section(\ref{ssec:explicit-lifting}) when we wish to give an explicit description of the ``lifted'' multitensors that this theorem provides for us. Let $(X,x)$ be in $V^M$, suppose $S$ is $\lambda$-accessible and fix an ordinal $n$ such that $|n| \geq \lambda$. Then by corollary(\ref{cor:explicit-phi-shreik}) one may take
\[ \begin{array}{lccr} {T(X,x) = (Q_n(X,x),a(X,x)\phi_{Q_n(X,x)})} &&& {a(X,x) = (q^{-1}_n)_{Q_n(X,x)} (v_{n})_{Q_n(X,x)}} \end{array} \]
as the definition of the endofunctor $T$. Note that $(Q_n(X,x),a(X,x))$ is just a more refined notation for $\phi_!(X,x)$. Referring to the diagram
\[ \xygraph{!{0;(2,0):(0,.6)::} {M^2X}="tl" [r] {MX}="tm" [r] {X}="tr" [d] {Q_n}="br" [l] {SX}="bm" [l] {SMX}="bl" "tl":@<-1ex>"tm"_-{Mx} "tl":@<1ex>"tm"^-{\mu^M_X}:"tr"^-{x} "bl":@<-1ex>"bm"_-{Sx} "bl":@<1ex>"bm"^-{\mu^S_XS(\phi_X)}:"br"_-{q_{{<}n}} "tl":"bl"_{\phi_{MX}} "tm":"bm"^{\phi_X} "tr":@{.>}"br"^{\eta^T_{(X,x)}}} \]
one may define the underlying map in $V$ of $\eta^T_{(X,x)}$ as the unique map making the square on the right commute. This makes sense since the top row is a coequaliser in $V$. Via the evident $M$-algebra structures on each of the objects in this diagram, one may in fact interpret the whole diagram in $V^M$ with the top row now being the canonical presentation coequaliser for $(X,x)$, and this is why $\eta^T_{(X,x)}$ is an $M$-algebra map. The proof that $\eta^T_{(X,x)}$ possesses the universal property of the unit of $\phi_! \ladj \phi^*$ is straight forward and left to the reader. As for $\mu^T_{(X,x)}$, it is induced from the following situation in $V^S$:
\[ \xygraph{{(SMQ_n,\mu^S)}="l" [r(2.5)] {(SQ_n,\mu^S)}="m" [r(4)] {(Q_n(Q_n,a\phi),a(Q_n,a\phi))}="r" [dl] {(Q_n(X,x),a(X,x))}="b" "l":@<-1ex>"m"_-{\mu^SS(\phi)} "l":@<1ex>"m"^-{S(a(X,x)\phi)}(:"r"^-{(q_{{<}n})_{(Q_n,a\phi)}}:@{.>}"b"^(.35){\mu^T_{(X,x)}},:"b"_{a(X,x)})} \]
Since by definition $\mu^T_{(X,x)}$ underlies an $S$-algebra map, to finish the proof that our definition really does describe the multiplication of $T$, it suffices by the universal property of $\eta^T$ to show that $\mu^T_{(X,x)}\eta^T_{T(X,x)}$ is the identity, and this is easily achieved using the defining diagrams of $\mu^T$ and $\eta^T$ together.

The data of $T$ is still not quite explicit enough for our purposes. What remains to be done is to describe $\eta^T$ and (especially) $\mu^T$ in terms of the transfinite data that gives $Q_n(X,x)$. So we shall for each ordinal $m$ provide
\[ \begin{array}{lccr} {\eta^{(m{+}1)}_{(X,x)} : X \to Q_{m{+}1}(X,x)} &&& {\mu^{(m)}_{(X,x)} : Q_m(Q_n(X,x),a(X,x)\phi) \to Q_n(X,x)} \end{array} \]
and $\mu^{(m{+}1)}_{(X,x)}$ in $V$ such that $\mu^{(m{+}1)}v_m=a(X,x)S(\mu^{(m)})$, by transfinite induction on $m$.
\\ \\
{\bf Initial step}. Define $\mu^{(0)}_{(X,x)}$ to be the identity, and $\eta^{(1)}_{(X,x)}$ and $\mu^{(1)}_{(X,x)}$ as the unique morphisms such that
\[ \begin{array}{lccr} {\eta^{(1)}_{(X,x)}x=(q_0)_{(X,x)}\phi_X} &&& {\mu^{(1)}_{(X,x)}(q_0)_{(Q_n,a\phi)}=a(X,x)} \end{array} \]
by the universal properties of $x$ and $q_0$ (as the evident coequalisers) respectively.
\\ \\
{\bf Inductive step}. Define $\eta^{(m{+}2)}=q_{m{+}1}\eta^{(m{+}1)}$ and $\mu^{(m{+}2)}$ as the unique map satisfying $\mu^{(m{+}2)}v_{m{+}1}=a(X,x)S(\mu^{(m{+}1)})$ using the universal property of $v_{m{+}1}$ as a coequaliser.
\\ \\
{\bf Limit step}. When $m$ is a limit ordinal define $\eta^{(m)}_{(X,x)}$ and $\mu^{(m)}_{(X,x)}$ as the maps induced by the $\eta^{(r)}$ and $\mu^{(r)}$ for $r<m$ and the universal property of $Q_m(X,x)$ as the colimit of the sequence of the $Q_r$ for $r<m$. Then define $\eta^{(m{+}1)}=q_m\eta^{(m)}$ and $\mu^{(m{+}2)}$ as the unique map satisfying $\mu^{(m{+}2)}v_{m{+}1}=a(X,x)S(\mu^{(m{+}1)})$ using the universal property of $v_{m{+}1}$ as a coequaliser.
\\ \\
The fact that the induction just given was obtained by unpacking the descriptions of $\eta^T$ and $\mu^T$ of the previous paragraph in terms of the transfinite construction of the endofunctor $T$ (ie the $Q_m(X,x)$), is expressed by
\begin{corollary}\label{cor:induced-monad-very-explicit}
Suppose that $V$ is a locally presentable category, $M$ and $S$ are monads on $V$, $\phi:M{\rightarrow}S$ is a morphism of monads, and $(X,x)$ is an $M$-algebra. If moreover $S$ is $\lambda$-accessible for some regular cardinal $\lambda$, then for any ordinal $n$ such that $|n| \geq \lambda$
one may take
\[ \begin{array}{lcccr} {T(X,x) = (Q_n(X,x),a(X,x)\phi_{Q_n(X,x)})} && {\eta^T_{(X,x)}=\eta^{(n)}_{(X,x)}} && {\mu^T_{(X,x)}=\mu^{(n)}_{(X,x)}} \end{array} \]
as constructed above as an explicit description underlying endofunctor, unit and multiplication of the monad generated by the adjunction $\phi_! \ladj \phi^*$.
\end{corollary}
\noindent and the simplification coming from proposition(\ref{prop:simple-coeq}) gives
\begin{corollary}\label{cor:vexp-simple}
Under the hypotheses of corollary(\ref{cor:induced-monad-very-explicit}), if for $(X,x) \in V^M$, $S$ and $S^2$ preserve the coequaliser of $\mu_{X}^{S}S(\phi_X)$ and $Sx$ in $V$, then one may take
\[ \begin{array}{lcccr} {T(X,x) = (Q_1(X,x),w\phi)} && {\eta^T_{(X,x)}=\eta^{(1)}_{(X,x)}} && {\mu^T_{(X,x)}=\mu^{(1)}_{(X,x)}} \end{array} \]
with $w$ as constructed in corollary(\ref{cor:phi-shreik-simple}).
\end{corollary}

\end{document}